\documentclass[12pt]{article}
\usepackage{amsmath}
\usepackage{amssymb}
\numberwithin{equation}{section}

\title{ Integrals over Grassmannians and
 Random permutations }

\author{
M. Adler\thanks{ Department of Mathematics, Brandeis
University, Waltham, Mass 02454, USA. E-mail:
adler@math.brandeis.edu.  The support of a National
Science Foundation grant \# DMS-98-4-50790 is
gratefully acknowledged.}~~~~~~ P. van
Moerbeke\thanks{ Department of Mathematics,
Universit\'e de Louvain, 1348 Louvain-la-Neuve,
Belgium and Brandeis University, Waltham, Mass 02454,
USA. E-mail: vanmoerbeke@geom.ucl.ac.be and
@math.brandeis.edu. The support of a National Science
Foundation grant \# DMS-98-4-50790, a Nato, a FNRS and
a Francqui Foundation grant is gratefully
acknowledged.}}

\date{October 8, 2001}

\catcode`\@=11
\let\c@equation=\relax
\newcounter{equation}[subsection]

\catcode`\@=12

\newcommand{\MAT}[1]{\left(\begin{array}{*#1c}}
\newcommand{\mat}{\end{array}\right)}
\newcommand{\qed}{\leavevmode\unskip\nobreak\penalty200\hskip2pt\null
\nobreak\hfill\rule{1.1ex}{1.1ex}%\parfillskip=0pt
\medbreak
}

\newcommand{\rg}{\rightarrow}
\newcommand{\lrg}{\longrightarrow}

\newcommand{\LR}{{\cal L}}

\newcommand{\BC}{{\mathbb C}}

\newcommand{\BF}{{\mathbb F}}
\newcommand{\BH}{{\mathbb H}}

\newcommand{\BY}{{\mathbb Y}}

\newcommand{\Sg}{\Sigma}
\newcommand{\iy}{\infty}
\newcommand{\pl}{\partial}
\newcommand{\al}{\alpha}

\newenvironment{proof}
{\medskip\noindent{\it Proof:\/} }{\qed}
\newenvironment
        {remark}{\medskip\noindent\underline{\it Remark:\/} }{\medbreak}

\newcommand{\la}{\langle}
\newcommand{\ra}{\rangle}

\newcommand{\dt}{\delta}
\newcommand{\Dt}{\Delta}
\newcommand{\sg}{\sigma}
\newcommand{\BR}{{\mathbb R}}
\newcommand{\lb}{\lambda}

\newcommand{\BJ}{{\mathbb J}}

\newcommand{\diag}{\operatorname{diag}}

\def\be#1\ee{\begin{equation}#1\end{equation}}
\def\bea#1\eea{\begin{eqnarray}#1\end{eqnarray}}
\def\bean#1\eean{\begin{eqnarray*}#1\end{eqnarray*}}

\newcommand{\Tr}{\operatorname{\rm Tr}}

\newtheorem{definition}{Definition}[section]
\newtheorem{theorem}[definition]{Theorem}
\newtheorem{lemma}[definition]{Lemma}
\newtheorem{corollary}[definition]{Corollary}
\newtheorem{proposition}[definition]{Proposition}

\makeatletter
\def\ps@X{\let\@mkboth\@gobbletwo
        \def\@oddhead{\tt Adler-van\,Moerbeke:%
        Grassmannian\hfil %\today
        October 8, 2001
        \hfil\S\thesection,
        p.\thepage}
        \def\@oddfoot{\rm\hfil\thepage\hfil}
        \let\@evenhead\@oddhead
        \let\@evenfoot\@oddfoot}
\makeatother

\catcode `!=11

\newdimen\squaresize
\newdimen\thickness
\newdimen\Thickness
\newdimen\ll! \newdimen \uu! \newdimen\dd! \newdimen \rr! \newdimen
\temp!

\def\sq!#1#2#3#4#5{%
\ll!=#1 \uu!=#2 \dd!=#3 \rr!=#4
\setbox0=\hbox{%
 \temp!=\squaresize\advance\temp! by .5\uu!
 \rlap{\kern -.5\ll!
 \vbox{\hrule height \temp! width#1 depth .5\dd!}}%
 \temp!=\squaresize\advance\temp! by -.5\uu!
 \rlap{\raise\temp!
 \vbox{\hrule height #2 width \squaresize}}%
 \rlap{\raise -.5\dd!
 \vbox{\hrule height #3 width \squaresize}}%
 \temp!=\squaresize\advance\temp! by .5\uu!
 \rlap{\kern \squaresize \kern-.5\rr!
 \vbox{\hrule height \temp! width#4 depth .5\dd!}}%
 \rlap{\kern .5\squaresize\raise .5\squaresize
 \vbox to 0pt{\vss\hbox to 0pt{\hss $#5$\hss}\vss}}%
}%end of \hbox
 \ht0=0pt \dp0=0pt \box0
}%end of \sq!

\def\vsq!#1#2#3#4#5\endvsq!{\vbox to \squaresize{\hrule
width\squaresize height 0pt%
\vss\sq!{#1}{#2}{#3}{#4}{#5}}}

\newdimen \LL! \newdimen \UU! \newdimen \DD! \newdimen \RR!

\def\vvsq!{\futurelet\next\vvvsq!}
\def\vvvsq!{\relax
  \ifx     \next l\LL!=\Thickness \let\continue=\skipnexttoken!
  \else\ifx\next u\UU!=\Thickness \let\continue=\skipnexttoken!
  \else\ifx\next d\DD!=\Thickness \let\continue=\skipnexttoken!
  \else\ifx\next r\RR!=\Thickness \let\continue=\skipnexttoken!
  \else\def\continue{\vsq!\LL!\UU!\DD!\RR!}%
  \fi\fi\fi\fi
  \continue}

\def\skipnexttoken!#1{\vvsq!}

\def\place#1#2#3{\vbox to 0pt{\vss
\rlap{\kern#1\squaresize
  \raise#2\squaresize\hbox{$#3$}}
\vss}}

\def\Young#1{\LL!=\thickness \UU!=\thickness \DD! = \thickness \RR! =
\thickness \vbox{\smallskip\offinterlineskip
\halign{&\vvsq! ## \endvsq!\cr #1}}}

\catcode `!=12

\begin{document}
\maketitle

\begin{abstract}
In testing the independence of two Gaussian
populations, one computes the distribution of the
sample canonical correlation coefficients, given that
the actual correlation is zero.  The ``Laplace
transform" of this distribution is not only an
integral over the Grassmannian of p-dimensional planes
in complex $n$-space, but is also related to a
generalized hypergeometric function. Such integrals
are solutions of Painlev\'e-like equations. They also
have expansions, related to random words of length
$\ell$ formed with an alphabet of $p$ letters. Given
that each letter appears in the word, the maximal
length of the disjoint union of $p$ increasing
subsequences of the word clearly equals $\ell$. But
the maximal length of the disjoint union of $p-1$
increasing subsequences leads to a non-trivial
distribution. It is precisely this probability which
appears in the expansion above.
\end{abstract}

%\newpage

\tableofcontents

\vspace{1cm}

\subsubsection*{Random words, longest increasing
sequences and mean hook lengths}

Consider the set of words
  $$
  \pi\in S_{\ell}^p:=
\left\{\mbox{words $\pi$ of length $\ell$, built from
an alphabet $\{1,...,p\}$}\right\},
 $$
  with the uniform probability distribution
\be
 P^{\ell,p}(\pi)=\frac{1}{p^{\ell}}.
  \ee
 The RSK correspondence (see section 2.1) between words and
 pairs of semi-standard and standard tableaux
 induces a probability measure on partitions
  \be
  \lb \in \BY_{\ell}=\{\mbox{partitions
 $ \lb \in \BY$ of weight $|\lb|=\ell$} \},
  \ee
 given by
 \be
  P^{\ell,p}(\lb)=
        \frac{f^{\lb}~s_{\lb}(1^p)}{p^{\ell}},
        ~~|\lb|=\ell,
 \ee
 where $s_{\lb}$ is the Schur polynomial associated
 with the partition $\lb$,
 \be
  1^p=(\overbrace{1,\ldots,1}^p,0,0,\ldots)~\mbox{ and }~ f^{\lb}=\# \{\mbox{standard tableaux of shape
$\lb$}\},
 \ee
with\footnote{$\lb^{\top}$ is the dual partition, i.e.,
  obtained by flipping the Young diagram $\lb$ about
  its diagonal. So, $\lb_1^{\top}$ is the length of the
  first column of $\lb$.
 }
$$(\mbox{ support }P^{\ell,p})\subseteq
 \BY_{\ell}^{(p)}:=\{\lb \in \BY_{\ell},~ \mbox{such that}~
\lb_1^{\top}\leq p \}.$$
 A subsequence $\sigma$ of the word $\pi$ is {\em weakly
 $k$-increasing}, if it can be written as
 \be
 \sigma=\sigma_1\cup \sigma_2\cup \ldots \cup
 \sigma_k,
  \ee
 where $\sigma_i$ are disjoint weakly increasing subsequences of the word
 $\pi$, i.e., possibly with repetitions.
  The length of the longest increasing/decreasing
  subsequences is closely related to the shape of the
  associated partition, via the RSK correspondence
  :
 \bea
 d_1(\pi)&=& \left\{\begin{array}{l} \mbox{length of
the longest {\em strictly}}\\\mbox{ decreasing
 subsequence of $\pi$}
\end{array}
\right\}=\lambda^{\top}_1
 \nonumber\\
i_k(\pi)&=&
 \left\{\begin{array}{l}
\mbox{length of the
 longest {\em weakly}}\\ \mbox{ $k$-increasing
 subsequence of $\pi$}
\end{array}
\right\}= \lambda_1+\ldots+\lambda_k
 \nonumber\\
 \eea

   For integer $0\leq p < n$, consider
   the fixed rectangular Young diagram of width $n-p$,
\be\mu=(n-p)^p:=(\overbrace{n-p,n-p,\ldots,n-p}^p).\ee
Consider a word $\pi \in S_{\ell}^p$. Then the
statement $d_1(\pi)=p$ implies, in particular, that
all letters of the alphabet $\{1,...,p\}$ are
represented in $\pi$; then automatically
$i_p(\pi)=\ell$. The theorem below deals with the
first non-trivial quantity $i_{p-1}(\pi)$, given that
$d_1(\pi)=p$ .
 Using the standard
notation, defined for a general parameter $\beta>0$,
 \be
  (a)_{\lb}^{(1/\beta)}:=
\prod_i(a+\beta(1-i))_{\lb_i}  %~~\beta=1/2
 ,
\mbox{with $(x)_n:=x(x+1)\ldots (x+n-1),~x_0=1$}, \ee
we now state Theorem 0.1, which will be established
in section 3.2 (note that here symbol (0.0.8) is used for
 $\beta=1$):

$$
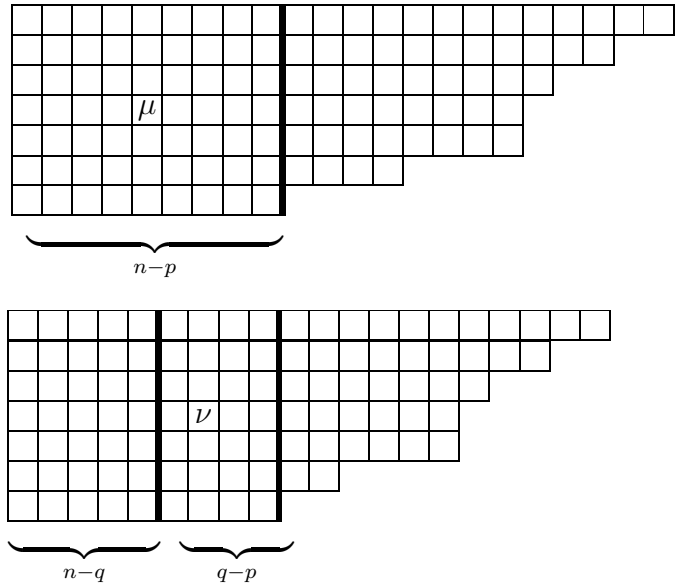
\squaresize .4cm \thickness .01cm \Thickness .07cm
\Young{
 &&&&&&&&r&&&&&&&&&&&&&\cr
 &&&&&&&&r&&&&&&&&&&&\cr
 &&&&&&&&r&&&&&&&&&\cr
 &&&&\mu&&&&r&&&&&&&&\cr
 &&&&&&&&r&&&&&&&&             \cr
 &&&&&&&&r&&&&                 \cr
 &&&&&&&&r                   \cr }
 % \right.
 $$

\vspace{-1cm}

 $$\hspace{.2cm}\underbrace{\hspace{3.4cm}}_{n-p}
%\hspace{.01cm} ~~~
%~\underbrace{\hspace{1.5cm}}_{q-p}
\hspace{5.2cm}$$

$$ \squaresize .4cm \thickness .01cm \Thickness .07cm
\Young{
 &&&&r&&&&r&&&&&&&&&&&\cr
 &&&&r&&&&r&&&&&&&&&\cr
 &&&&r&&&&r&&&&&&&\cr
 &&&&r&&\nu&&r&&&&&&\cr
 &&&&r&&&&r&&&&&&             \cr
 &&&&r&&&&r&&                 \cr
 &&&&r&&&&r                   \cr }\hspace{.9cm}$$

\vspace{-1cm}

 $$\hspace{.4cm}\underbrace{\hspace{2.0cm}}_{n-q}
%\hspace{.01cm} ~~~
~\underbrace{\hspace{1.5cm}}_{q-p} \hspace{5.5cm}$$

%\newpage

\begin{theorem} Given the probability (0.0.1) and (0.0.3)
, the following holds ($h^{\kappa}$ denotes the
product of hook lengths over all boxes of the
partition $\kappa$): ($0\leq p < n$)
   \bea
 P^{\ell,p}\left( \lb\supseteq\mu
   \right)&=&
   P^{\ell,p}\left(\pi \in S_{\ell}^p~\Bigr|~{
 {d_1(\pi)=p
 \mbox{ and} }\atop {  i_{p-1}(\pi)\leq
 \ell-n+p} } \right) \nonumber \\
 &=&
  \frac{\ell
!}{p^{\ell}}\prod_1^p\frac{(p-i)!}{(n-i)!}\sum_{\kappa
\in \BY_{\ell -p(n-p)} \atop \kappa_1^{\top}\leq
p}\frac{1}{(h^{\kappa})^2}
 \frac{((p)_{\kappa})^2 }{(n)_{\kappa}}
%\prod^p_{i=1}
%\frac{\left((p-i+1)_{\kappa_i}\right)^2}
%{(n-i+1)_{\kappa_i}}
.
%
%\nonumber\\
   \eea
   More generally, for fixed $p\leq q <n$,
   the mathematical expectation
   (with regard to $P^{\ell,p}$)
   of the hook length of $\lb$,
   emanating from the vertical strip
   $\nu$
   (of width $q-p$),
   % of the Young diagrams $\lb$, averaged
   %with regard to $P^{\ell,p}$,
   equals:

\be
 E^{\ell,p}\left(I_{\{\lb\supseteq\mu\}}
 (\lb)
   \prod_{{(i,j)\in \lb}
 \atop {n-q<j\leq n-p}}
  h_{(i,j)}^{\lb}
   \right)=
    \frac{\ell
!}{p^{\ell}}
 \prod_1^p\frac{(q-i)!}{(n-i)!}
  \sum_{\kappa
\in \BY_{\ell -p(n-p)} \atop \kappa_1^{\top}\leq
p}\frac{1}{(h^{\kappa})^2}
 \frac{(p)_{\kappa}(q)_{\kappa}}{(n)_{\kappa}}
%\prod^p_{i=1}
%\frac{(p-i+1)_{\kappa_i}(q-i+1)_{\kappa_i}}
% {(n-i+1)_{\kappa_i}}
. \ee

\end{theorem}

%\newpage

\subsubsection*{Generating function for the
mathematical expectation of the
hook length, integrals over Grassmannians
 and Painlev\'e V
  } Theorem 0.2 below involves an {\em integral over
 the
 Grassmannian} \be
  Gr(p, \BC^n)
 = \frac{U(n)}{U(p)\times U(n-p)}=:G/K
 %&&&&
 \ee
  of $p$-dimensional planes in $\BC^n$ through the
  origin
 and Haar measure $d\mu(Z)$ on $Gr(p, \BC^n)$, expressed
 in the parametrizing coordinate $Z$ of
  \be
 \mbox{Affine~}Gr(p,\BC^n)=\left\{\mbox{span}
 \MAT{1}I_p\\Z\mat
\Big| Z:%=Z(A)
 = A_{21}A_{11}^{-1}
,~A\in G\right\},
 \ee
where $A\in G$ is represented in block form:

$\hspace*{39mm}\stackrel{p}{\longleftrightarrow}
\quad\stackrel{n-p}{\longleftrightarrow}$

\vspace{-.7cm}

\be
 \left(\begin{array}{cc}A_{11}&A_{12}\\
A_{21}&A_{22}\end{array}\right)
 \begin{array}{l}\updownarrow p\\ \updownarrow n-p\end{array}
\in G.
 \ee

 Section 1 will be devoted to the geometry of
 $Gr(p,\BF^n)$, where $\BF=\BC,~\BR$ or $\BH$, and to
 the study of integrals over $Gr(p,\BF^n)$.

 We also need {\em Jack polynomials},
  which are the unique symmetric functions
 orthogonal with respect to a certain
 $\al$-dependent inner-product $\la,\ra_{\al}$, such
 that
 \be
  \la J_{\lb}^{(\al)},J_{\mu}^{(\al)}\ra
=\dt_{\mu\lb}j_{\lb}^{(\al)},
 \ee
 with
  $$
j_{\lb}^{(\al)}=\prod_{(i,j)\in\lb}
 \left(\lb^{\top}_j-i+\al(\lb_{i}-j
+1)\right)\left(\lb^{\top}_j-i+1+\al(\lb_i-j)\right)
 .$$
Facts about Jack polynomials relevant for this project
will be discussed in section 2.3.

Finally, {\em generalized hypergeometric functions}
${}_2F_1^{(\al)}$ are defined by: ($p,q,n\in \BC,~
x=(x_1,x_2,\ldots)$)
\be
{}_2F_1^{(\al)}(p,q;n;x): =\sum_{\kappa \in \BY}
 \frac{(p)^{(\al)}_{\kappa}(q)^{(\al)}_{\kappa}}{(n)^{(\al)}_{\kappa}}
~ \al^{|\kappa|}\frac{J_{\kappa}^{(\al)}(x)}
 {j_{\kappa}^{(\al)}}.
 \ee
In particular,
 \be
{}_2F_1^{(1)}(p,q;n;x): =\sum_{\kappa \in \BY}
 \frac{(p)^{(1)}_{\kappa}(q)^{(1)}_{\kappa}}{(n)^{(1)}_{\kappa}h^{\kappa}}
~s_{\kappa}(x).
 \ee
The following theorem will be established in section
4.2; it is strongly motivated by certain integrals
appearing in the context of testing statistical
independence of Gaussian populations, as will be
explained in the next paragraph.
\begin{theorem}
For fixed $p\leq q\leq n/2$, the generating function
for the mathematical expectation of the hook length
(0.0.10) over a strip, with
regard to the probability (0.0.3), is given
by\footnote{Remember $(\mbox{ support }P^{\ell,p})
 \subseteq\{\lb \in \BY,~ \mbox{such that}~
|\lb|=\ell,~\lb_1^{\top}\leq p \}.$}

 \bea
\lefteqn{\hspace{-2cm}
 \prod_{i=1}^p \frac{(n-i)!}{(q-i)!}
 ~ x^{-p(n-p)}
% \frac{\prod_{i=1}^p i!~(n-q-i)!}{x^{p(n-p)}}
  \sum_{\ell\geq p(n-p)}
  \frac{(px)^{\ell}}{\ell!}
 E^{\ell,p}\left( I_{\{\lambda \supseteq \mu \}}(\lb)
  \prod_{{(i,j)\in \lb}
 \atop {n-q<j\leq n-p}}
  h_{(i,j)}^{\lb}
   \right)
   }
  \nonumber \\
   &=&
%    ({\tilde c_{n,q,p}})^{-1}
 (c^{(1)}_{n,q,p})^{-1}  \int_{Gr(p,\BC^n)}e^{x\Tr(I+Z^{\dag}Z)^{-1}}\det(Z^{\dag
}Z)^{-(q-p)} d\mu(Z)  \nonumber\\
 &=&
  ~{}_2F_1^{(1)}(p,q;n;y)
 \Bigr|_{\sum_{\ell}y_{\ell}^i= \delta_{1i}x}
  \nonumber \\
  &=&
    \exp {\displaystyle{\int_0^x\frac{u(y)-p(n-p)+py}{y}
    dy}
    }
 \eea
       where
$u(x)$
  is the unique solution to the initial value problem:
 \be
 \left\{\begin{array}{l} \displaystyle{ x^2
u^{\prime\prime\prime}+xu^{\prime\prime}
 +6x{u^{\prime} }^2-4uu^{\prime}+4Qu^{\prime}
 -2Q^{\prime}u+2R=0 } \\\hspace{10cm}
  \mbox{\bf (Painlev\'e V)} \\
\displaystyle{\mbox{with}~~  u_{}(x)=p(n-p)-
 \frac{p(n-q)}{n}x+\ldots+
+  a_{n+1}x^{n+1}+O(x^{n+1})+\ldots ,~\mbox{near}~ x=0
.}
\end{array}
 \right.
 \ee
  with $a_{n+1}$ specified
  by the hypergeometric function in (0.0.17) and
  \bea
 4Q&=&
 -x^2+2(n+2(p-q))x-(n-2p)^2 \nonumber\\
 2R&=&
  p(p-q)(x+n-2p).
  \eea
\end{theorem}

\remark The constant $c^{(1)}_{n,q,p}$ in (0.0.17) is
the one below for $\beta =1$:
  \be
c:=c^{(\beta)}_{n,q,p}:=\prod^p_{i=1}\frac{\Gamma(i{\beta}+1)
 \Gamma(\beta(n-q-i+1)) \Gamma(\beta(q-i+1))}
{\Gamma(\beta +1) \Gamma(\beta (n-i+1) )}. \ee

  \subsubsection*{Testing Statistical Independence
 of Gaussian Populations }

To summarize section 5, consider $p+q$ normally distributed random variables
$(X_1,...,X_p)^{\top}$ and $(Y_1,...,Y_q)^{\top}$
 ($p\leq q$) with mean zero and
covariance matrix $\Sigma$. According to Hotelling
(\cite{Hotelling})(see also Muirhead \cite{Muirhead},
p106), $(X_1,...,X_p)^{\top}$ and
$(Y_1,...,Y_q)^{\top}$ can be replaced by a linearly
transform of the  $X$'s and $Y$'s, so that the
covariance matrix takes on the canonical form:

\smallbreak

$\hspace*{44mm}\stackrel{p}{\longleftrightarrow}
\quad\stackrel{q}{\longleftrightarrow}$

\vspace{-.8cm}

 $$\Sigma=\mbox{cov}\MAT{1}X\\Y\mat =
 \left(\begin{array}{cc}\Sigma_{11}&\Sigma_{12}\\
\Sigma_{12}^{\top}&\Sigma_{22}\end{array}\right)
 \begin{array}{l}\updownarrow p\\ \updownarrow q\end{array}
\longrightarrow \Sigma_{{can}}=\MAT{2}I_{p}&P
\\P^{\top}&I_q\mat,
$$
 where
 $$ \hspace{1cm}{{q}\atop{\longleftarrow\longrightarrow}}$$

 \vspace{-.3cm}

  $$
 P=\left(
\begin{array}{lllllll|l}
\rho_1& & & & & & & \\
 &\ddots& & O& & & &\\
 & &\rho_k& & & & &O\\
& & &\rho_{k+1}& & & &  \\ & O& & & &\ddots& &\\ & & &
& & & \rho_{p}&
\end{array}
\right){\Bigg\updownarrow} p,\quad
k=\mbox{\,rank\,}\Sigma_{12},\\
%\hspace*{3cm}{{\longleftarrow\longrightarrow}\atop{p}}
$$ is the $p\times q$-matrix of canonical correlation
coefficients:
  $$
1\geq\rho_1\geq\rho_2\geq\ldots\geq\rho_k
>0,~\rho_{k+1}= ...=\rho_p=0.
%\quad\mbox{(canonical
%correlation coefficients),}
 $$

% From here on, we may take $\Sigma=\Sigma_{can}$.
 The $n$ ($n\geq p+q$) independent samples
$(x_{11},\ldots,x_{1p},y_{11},\ldots,y_{1q})^{\top},\ldots$,
$(x_{n1},\ldots,x_{np},y_{n1},\ldots,y_{nq})^{\top}$,
arising from observing $\MAT{1}X\\Y\mat$ lead to a
matrix $\MAT{1}x\\y\mat$ of size $(p+q,n)$, having a
normal distribution with correlation $\Sigma$.  The
roots $r_1^2,\ldots, r_p^2$ (sample canonical
correlation coefficients) of the equation
  $$
\det(xy^{\top}(yy^{\top})^{-1}yx^{\top}-r^2xx^{\top})=
0 $$ are the estimators (maximum likelihood
estimators) of the canonical correlation coefficients
 $\rho_1^2,\ldots, \rho_p^2$.

In testing the {\em null hypothesis}, versus the {\em alternative
 hypothesis},

$$   H_0: ~~\rho_1^2=\ldots =\rho_p^2=0
~~~\mbox{versus} ~~H_a: ~~(\rho_1^2,\ldots,
\rho_p^2)\neq 0
 $$
 one needs the joint density of
the $r^2_i$, given that $\rho_1^2=\ldots =\rho_p^2=0$;
namely, up to a $(q,p,n)$-dependent normalizing
constant, the density is given by
 \be %\pi^{p^2/2}c_{n,p,q}
% \prod_{1\leq i<j\leq p}(r_i^2-r_j^2)
  |\Delta_p(r^2)|^{\beta}
  \prod^p_{i=1}(r^2_i)^{\beta(q-p+1)-1}
   (1-r_i^2)^{\beta(n-q-p+1)-1}
dr_i^2
 \ee
for $\beta=1/2$. This formula generalizes to the
formula above, upon considering random variables
$(X_1,...,X_p)^{\top}$ and $(Y_1,...,Y_q)^{\top}$,
with values in the complex $\BC$ ($\beta=1$) and the
quaternions $\BH$ ($\beta=2$).

\subsubsection*{Expectation of the ratio of Jack polynomials, integrals
over Grassmannians, sample canonical correlation
coefficients and PDE's }

%{\bf Probabilities on partitions, integrals over
%Grassmannians and PDE's }

Consider now a Poissonized probability on partitions $\lb
\in \BY$, which depends on a parameter $x>0$:
 \bea
P_{x,p}(\lb)&=&e^{-\beta px}\frac{x^{|\lb|}
J^{(1/\beta)}_{\lb}(1^p)}{
 j_{\lb}^{(1/\beta)}}=
 e^{-\beta px}
  \frac{(\beta px)^{|\lb|}}{|\lb|!}
 P^{\ell,p}(\lb)
 ,~~~\lb \in \BY,\nonumber\\
\eea
 with $P^{\ell,p}$ generalizing probability measure
  (0.0.3), as we shall see in section 3.1,
  $$
   P^{\ell,p}(\lb)=
  \frac{J^{(1/\beta)}_{\lb}(1^p)}{j_{\lb}^{(1/\beta)}}
  \frac{|\lb|!}{(\beta
  p)^{|\lb|}}~~\mbox{with}~~|\lb|=\ell,
  $$

 This probability has its support
  on $\lb^{\top}_1 \leq p$. Many of these probability
  distributions on partitions have been introduced and
  extensively studied by Borodin, Kerov, Okounkov,
  Olshanski and Vershik (see \cite{Borodin1,
  Borodin2,Borodin3,Vershik and Kerov,Kerov})
 The following statement involves an integral
 suggested again by the statistical
theory mentioned earlier.

\begin{theorem}
For fixed $p\leq q\leq n/2$, the following holds
($\beta=1/2,1,2$)
  \bea
I^{(\beta)}_{n,p,q}&=&
 ce^{\beta
px}E_{x,p}\left(\frac{1}{\beta^{|\lb|}}\frac{J_{\lb}^{(1/\beta)}(1^q)}
{J_{\lb}^{(1/\beta)}(1^n)}%I_{\{\lb_1^{\top}\leq p\}}(\lb)
 \right)
 \nonumber \\&=&
  \int_{Gr(p,\BF^n)}e^{x\Tr(I+Z^{\dag}Z)^{-1}}\det(Z^{\dag
}Z)^{-\beta(q-p)} d\mu(Z)  \nonumber \\
 &=&\int_{[0,1]^p}e^{x
\sum^p_1z_i}
|\Delta_p(z)|^{2\beta}\prod_1^p
z_i^{\beta(q-p+1)-1}(1-z_i)^{\beta(n-q-p+1)-1}dz_i
 \nonumber\\
 &=&
 c~
   {}_2F_1^{(1/\beta)}(\beta p,\beta q;\beta n;y)
    \Bigr|_{\sum_{\ell}y_{\ell}^i= \frac{x}{\beta}\delta_{1i}}
      \nonumber\\  %&&\\
&=&  c \exp \int_0^x v(y)dy.
 \eea
 where $c=c^{(\beta)}_{n,q,p}$ is as in (0.0.20) and where%\footnote{  }
\begin{itemize}
  \item
   $d\mu(Z)$ is Haar measure on
  the space $Gr(p,\BF^n)$ of $p$-planes in $\BF^n$, where
  $\BF:= \BC, \BR$ or the quaternions $\BH$.
  \item
  The integral over $[0,1]^p$,
 appearing in (0.0.23), is the ``Laplace transform" of the distribution of
 the sample canonical correlation coefficients (0.0.21).
 This integral is a H\"ankel determinant for $\beta=1$ and a
  Pfaffian for $\beta=1/2$ and $2$.

   \item

$v(y):= v_{n,p,q}^{(\beta)}(y)$ and $I_p:=
 I^{(\beta)}_{n,p,q}$ satisfies the differential
equation(define $\delta_1^{\beta}:=1$ for $\beta=1$
and $:=0$ otherwise):
\bea && \hspace{-1cm}
 4\left(y^3
v^{\prime\prime\prime}+6y^3{v^{\prime}
}^2+(1+\delta^{\beta}_{1})(2 y^2
 v^{\prime\prime}
 +4 y^2vv^{\prime}
  + yv^2)\right)-yP_0v^{\prime}+P_1v
  +P_2  \nonumber\\
%  &&
% -yP_0H^{\prime}+P_1H
%  +P_2\\
&&\nonumber\\ &&%-yP_0H^{\prime}+P_1H
   =
 \left\{ \begin{array}{ll}
    0,~~& \mbox{for}~~ \beta =1,
    \mbox{\bf (\bf Painlev\'e V)} \nonumber\\
    \nonumber\\
   \displaystyle{  \frac{3}{16}\frac{p(p-1)}{(p+1)(p+2)}y^3~\frac{I_{p-2} I_{p+2}}{I_{p}^2}
} ,~~& \mbox{for}~~ \beta =1/2, \\
        \\
  \displaystyle{  \frac{3}{16^2}\frac{p}{p+1}y^3~\frac{I_{p-1}I_{p+1}}{I_{p}^2}
},~~& \mbox{for}~~ \beta =2,
 \end{array}
  \right.
  \eea
with $P_0$ quadratic and $P_1$, $P_2$ linear
polynomials in y, with coefficients depending on $n$
and
  \be
  r=pq,~~~~s=n-2p-2q.
   \ee
\end{itemize}

\end{theorem}

This statement will be established in section 4.1 and
the differential equation part in section 6.1. As a
by-product, we show incidentally that the multivariate
hypergeometric function ${}_2F_1^{(1)}( p, q; n;y)$
expressed in the
    $it_i=\Sigma_{k\geq 1}
y^i_k$-variables are $\tau$-functions for the
KP-hierarchy; but also that the function
${}_2F_1^{(1)}( p, q; n;y)$ properly restricted is a
solution of Painlev\'e V.
 For related questions, see Orlov and Sherbin
(\cite{Orlov1,Orlov2}). Section 7 gives new
differential equations for the spectrum of Wishart
matrices and for the sample canonical correlations of
Gaussian populations.

\underline{\sl Acknowledgment}: The authors thank
Professors I. Gessel, S. Helgason, B. Lian, G.
Schwarz, R. Stanley and C.-L. Terng for useful advice,
especially regarding section 1.

\newpage

\section{Integrals over Grassmannians}

Consider the
 Grassmannian $Gr(p, \BF^n)$ of $p$-planes through the origin in
 $\BF^n$,
   where $\BF= \BC,~\BR$ or $\BH$ ($=$ quaternions).
Let $G$ be the group of
matrices $A$, with entries in $\BF$, such that
$A^{-1}=A^{\dag}$, with\footnote{Given $a_i\in \BR$, we define
for $a=a_0+ia_1
\in \BC$, $\bar a=a_0-ia_1$
and for $a=a_0+a_1e_1+a_2e_2+a_3e_3\in \BH$, $\bar a:=
a_0-a_1e_1-a_2e_2-a_3e_3$.}
$A^{\dag}:=\bar A^{\top}$.
  Matrices
$A\in G$ will be represented as block matrices

$\hspace*{39mm}\stackrel{p}{\longleftrightarrow}
\quad\stackrel{n-p}{\longleftrightarrow}$

\vspace{-.5cm}

\be
 \left(\begin{array}{cc}A_{11}&A_{12}\\
A_{21}&A_{22}\end{array}\right)
 \begin{array}{l}\updownarrow p\\ \updownarrow n-p\end{array}
\in G.
 \ee

   The main statement of this section is theorem
1.1, where it is assumed, without loss of generality,
that $n\geq 2p$. The values of $\beta$ are related to
$\BC,\BR$ and $\BH$, as follows:

\medbreak

\hspace{1cm}{\small
\begin{tabular}{ll}

$Gr(p,\BC^n):$  &
 $\beta=1$

%\\ &
 \\

$Gr(p,\BR^n):$ & $\beta=1/2$
%  \\ &
 \\

 $Gr(p,\BH^n):$ &$\beta=2$

\\

\end{tabular}
}

The geometry of the symmetric spaces $G /K  $ and
$K\backslash G /K $ has been studied by Helgason
\cite{Helgason1,Helgason2}. In his recent Princeton
thesis, Due\~nez (\cite{Duenez}) has systematically
studied integrals over symmetric spaces. Explicit
information on this subject is not readily available
in the literature; therefore we explain the theory in
the Grassmannian case and the useful aspects for our
purposes.

%Sarnak \cite{Sarnak} in his
%MSRI-lectures suggested considering

\begin{theorem} Consider the two parametrization
 of $Gr(p,\BF^n)$
\be
 \mbox{Affine~}Gr(p,\BF^n)=\left\{\mbox{span}
 \MAT{1}I_p\\Z\mat
\Big| Z:=Z(A)=A_{21} A_{11}^{-1},~A\in G\right\},
 \ee
 and its invariant measures $d\mu(Z)$. Then,
 for $p\leq q\leq n/2$, we have
  \bea
    \lefteqn{\int_{Gr(p,\BF^n)}e^{x\Tr(I+Z^{\dag}Z)^{-1}}\det(Z^{\dag
}Z)^{-\beta(q-p)} d\mu(Z)} \nonumber\\
 &=& \int_{[0,1]^p}e^{x%\displaystyle{
\sum^p_1z_i}
%}
|\Delta_p(z)|^{2\beta}\prod_1^p
z_i^{\beta(q-p+1)-1}(1-z_i)^{\beta(n-q-p+1)-1}dz_i.
 \nonumber\\ \eea
This is the ``{\em Fourier transform}" of the joint
density of the sample canonical correlations
$(z_1,\ldots,z_p)=(r_1^2,\ldots,r^2_p)$, for $\beta
=1/2,1,2$, for the real, complex and quaternionic
cases, given that the canonical correlation
coefficients are zero. (see (0.0.21) and section 5)
\end{theorem}

%\noindent
Consider the following block matrix
 \bean
 I_{p,q}&=&
  \MAT{2}I_p&O
  \\
O&-I_q\mat.
 \eean

\begin{theorem} An alternative description for
$Gr(p,\BF^n)$ is given by

 \be
 Gr(p,\BF^n)\simeq {\cal S}:= \left\{M=AI_{p,n-p}A^{-1}
 I_{p,n-p}~\Big| ~A\in G\right\}
 \ee
The $n\times n$ matrices $M$ have a $n-2p$-dimensional
eigenspace corresponding to the eigenvalue $1$, so
that $M$ can be decomposed into $M=M_0\otimes M_1$
(with $M_1$ corresponding to the $1$-eigenspace), with
$Tr (M-M_0)=n-2p$. Then we have
 \bea
    \lefteqn{ e^{px/2} \int_{\cal S}e^{x Tr M_0/4}\det
  \left( \frac{I-M_0}{I+M_0} \right)^{\beta
  (p-q)}d\mu(M)}\nonumber\\
 &=& \int_{[0,1]^p}e^{x%\displaystyle{
\sum^p_1z_i}
%}
|\Delta_p(z)|^{2\beta}\prod_1^p
z_i^{\beta(q-p+1)-1}(1-z_i)^{\beta(n-q-p+1)-1}dz_i.
 \nonumber\\ \eea

 \end{theorem}
\begin{theorem} Considering the parametrization
\be
 T_{Id}Gr(p,\BF^n)=\{Z~\Bigr|~ \mbox{arbitrary
  $(n-p)\times p$ matrix,
 with $Z_{ij}\in\BF$}\},
 \ee
we have %of $T_{id}Gr(p,\BF^n)$
 $$ \int_{   { {Z\in
T_{id}Gr(p,\BF^n)} \atop {  \mbox{with
spectrum~}(Z^{\dag}Z)\in E  } }  } e^{-x\Tr Z^{\dag}Z}
\det(Z^{\dag}Z)^{\beta p}d\mu(Z),\qquad E\subset\BR^+
$$
 \be = \int_{E^p}|\Delta _p(u)|^{2\beta} \prod^p_1e^{-x
u_i}u_i^{\beta(n-p+1)-1}du_i. \ee
 For
$\beta =1/2$ and $x=\frac{1}{2\lb}$, this is an
integral of the joint density (Wishart density)(see
section 5 and Muirhead \cite{Muirhead},
 p.107) of the eigenvalues
$u_1,\ldots,u_p$ of the matrix $A=Z^{\dag}Z$, where
$Z$ is a $n\times p$ matrix ($p\leq n$), with Gaussian
density centered at $0$ and covariance $\lb I_p$,
namely the density
 $$
c_{n,p}(2\pi\lb)^{-np/2}e^{-\frac{1}{2\lb}\Tr
Z^{\dag}Z}\prod_{{{1\leq i\leq n}\atop{1\leq j\leq
p}}}dz_{ij}. $$
\end{theorem}

%\newpage

\begin{proposition}
Then
  \bea
Gr(p,\BC^n)&= & \frac{U(n)}{U(p)\times U(n-p)}=:G/K
 \nonumber\\
 Gr(p,\BR^n)&= & \frac{SO(n)}{SO(p)\times
SO(n-p)}=:G/K   \nonumber\\
 Gr(p,\BH^n)&= & \frac{Sp(n)}{Sp(p)\times
Sp(n-p)}=:G/K, \eea
%  with $K=K_1\times K_2$. Also, ($Z$
%is a $(n-m)\times m$ matrix)
and the affine part can be parametrized
as follows:
 $$
 \mbox{Affine}~ Gr(p,\BF^n)=
 \left\{ \mbox{span }\left( {I_p}\atop {Z} \right)~ \Bigr|~
 Z=Z(A)=A_{21}A_{11}^{-1}
    ,~ A \in G \right\}
 $$
 and
 \bea
 K\backslash G /K=
 \left\{   \mbox{span }\left( {I_p}\atop {Z} \right)~
   \left|~
 \begin{array}{l} Z
    =\left(
\begin{array}{lll}
-\tan \theta_1& &~~O      \\
 &\ddots&      \\
O & &-\tan \theta_p    \\ \hline
  & O_{n-2p,p}&
\end{array}
\right),~\\ \\ \hspace{1.6cm} \mbox{with}~~0\leq
\theta _i\leq {\pi}
\end{array} \right.
 \right\}; \nonumber\\&& \eea
also\footnote{Setting Id $=\left( {I_p}\atop {O}
\right)$.},
 \be
 T_{Id}Gr(p,\BF^n)=\{Z~\Bigr|~ \mbox{arbitrary $(n-p)\times p$ matrix,
 with $Z_{ij}\in\BF$}\}.
 \ee
Setting, for the respective cases of $Gr(p,\BF^n)$ and
$T_{Id}Gr(p,\BF^n)$,
 \bean
%z_i&:=&\frac{1}{1+\tan^2\theta_i}=\cos^2\theta_i,~~
\{y_1,\ldots,y_p\}&=&
 \mbox{spectrum}~\left(\frac{1-Z^{\dag}Z}{1+Z^{\dag}Z}
 \right)=
( \cos
2\theta_1,\ldots,\cos2\theta_p)~~
%1\leq i\leq p
\\
  \{u_1,\ldots,u_p\}&:=&\mbox{spectrum}~(Z^{\dag}Z)
  \eean
with $-1\leq
 y_i\leq 1,$ and $0\leq
 u_i < \iy,$
Haar measure on
 $Gr(p, \BF^n)$ and $T_{Id}Gr(p, \BF^n)$ reads,
 setting $k=n-2p$, ({\bf Weyl integration formulae})

\bea \mbox{$Gr(p, \BF^n)$}:&&~~d\mu(Z)=
 |\Dt_p(y)|^{2 \beta} \prod_1^p
(1-y_i)^{\beta k +(\beta-1) }
 (1+y_i)^{\beta-1}dy_idK\nonumber\\
 \mbox{$T_{Id}Gr(p, \BF^n)$}:&&~~d\mu(Z)=
 |\Dt_p(u)|^{2 \beta} \prod_1^p
u_i^{\beta k +(\beta -1)}
 du_idK,
\eea
leading to the table:

\vspace{.8cm}

\hspace{-1cm}{\small
\begin{tabular}{l|l|l}
$G/K$  &
 induced measure $d\mu$ on $K\backslash G /K  $
 & $d\mu$ on $T_{Id}G/K  $ \\
 \hline
 \\
$Gr(p,\BC^n)$  &
 $ \Dt_p(y)^2 \prod_1^p
(1-y_i)^{k} dy_i $ &
 $ \Dt_p(u)^2 \prod_1^p
u_i^{k} du_i $
\\ & & \\

$Gr(p,\BR^n)$  &
 $|\Dt_p(y)| \prod_1^p
(1-y_i)^{\frac{1}{2}(k-1)}
 (1+y_i)^{-\frac{1}{2}}dy_i $
  & $|\Dt_p(u)| \prod_1^p
u_i^{\frac{1}{2}(k-1)}
 du_i $ \\

& & \\

 $Gr(p,\BH^n)$ &

  %$\frac{f^{(2\lb^{\top})}}{(2|\lb|-1)!!}$&
  $\Dt_p(y)^4 \prod_1^p
(1-y_i)^{2k+1}
 (1+y_i)dy_i $
 & $\Dt_p(u)^4 \prod_1^p
u_i^{2k+1}
 du_i $\\

\end{tabular}
} $$ \mbox{Table 1}$$ %\vspace{.5cm}

In the other description (1.0.4) of $Gr(p,\BF^n)$,
given by
  $$
 Gr(p,\BF^n)\simeq {\cal S}:= \left\{M=AI_{p,n-p}A^{-1}
 I_{p,n-p}~\Big| ~A\in G\right\}
 $$
an appropriate left action of $B\in K$ on $A\in G$,
amounting to conjugation in ${\cal S}$, leads to the
matrix in the torus ${\frak A}\subset G$,
 \bean
(BA)I_{p,n-p}(BA)^{-1}
 I_{p,n-p}&=&
  B(AI_{p,n-p}A^{-1}
 I_{p,n-p})B^{-1}\\ &&\\
  &=&\left(\begin{array}{cccc}
 \Re e^{2i\Theta_p} & & \Im e^{2i\Theta_p} &  \\
                  & &                  & O\\
 -\Im e^{2i\Theta_p} && \Re e^{2i\Theta_p} & \\
  & O & & I_{n-2p} \\
\end{array}
\right) \in {\frak A},
 \eean
  with $\Theta_p:=
\diag(\theta_1,\ldots,\theta_p)$.

\end{proposition}

\medskip\noindent{\it Proof of Proposition 1.4:\/}
The $p$ columns of the $n\times p$ matrix
$A\MAT{1}I_p\\O_{}\mat$, with $A\in G,$ and
 $O:=O_{n-p,p}$ (a zero matrix of size $(n-p,p)$) span a $p$-dimensional plane in
$\BF^n$, so that $$
Gr(p,\BF^n)=\left\{\mbox{span}~A\MAT{1}I_p\\O\mat,
\mbox{~with~}A\in G\right\}. $$
  The {\bf right action} of $B\in K$ on $A\in G$
   acts on the $p$-plane $A\MAT{1}I_p\\O\mat$ as:
$$ AB\MAT{1}I_p\\O\mat = A\MAT{2}B_{11}&O\\
O&B_{22}\mat\MAT{1}I_p\\O\mat=A\MAT{1}B_{11}\\O\mat
=A\MAT{1}I_p\\O\mat B_{11} $$ and therefore it has no
effect on that plane $$ \mbox{span}~AB\MAT{1}I_p\\O\mat
=~\mbox{span}~A\MAT{1}I_p\\O\mat B_{11}=
 ~\mbox{span}~A\MAT{1}I_p\\O\mat %B_{11}
  , $$ since
multiplication to the right by $B_{11}$ merely
replaces the $p$ columns of $A\MAT{1}I_p\\O\mat$ by
$p$ linear combination.
 Then the $n\times p$
matrix $V:=A\MAT{1}I_p\\O\mat$ satisfies
\be
V^{\dag}V=(I_p\,\,\,O)~ A^{\dag}A\MAT{1}I_p\\O\mat = (I_p\,\,\,O)~
I_n\MAT{1}I_p\\O\mat =I_p.
\ee

{\em Conversely}, we show that
 \be
   \mbox{span}~A_1\MAT{1}I_p\\O\mat  =
\mbox{span}~A_2\MAT{1}I_p\\O\mat
 \ee
 implies $A^{-1}_2A_1 \in K=K_1\times K_2$
Indeed, (1.0.13) holds if and only if
\be
A_1\MAT{1}I_p\\O\mat = A_2\MAT{1}I_p\\O\mat g,
\mbox{~with an invertible $p\times p$ matrix $g$}. \ee
Then we prove $g\in K_1$. Indeed, from (1.0.13), the
matrices $V_i:=A_i\MAT{1}I_p\\O\mat $ satisfy
$V^{\dag}V=I_p$ and so $V_1=V_2g$ implies
$g^{\dag}g=(V_2g)^{\dag}V_2g=V_1^{\dag}V_1=I_p$.
Multiplying (1.0.14) to the left with $A_2^{-1}$
yields $$ A_2^{-1}A_1\MAT{1}I_p\\O\mat
=\MAT{1}g\\O\mat\mbox{ ~and thus~} G \ni
~A_2^{-1}A_1=\MAT{2}g&\ast\\ O&h\mat ; $$ the fact
that the latter matrix is in $G$, implies $\ast =0$
and $g\in K_1$, $h\in K_2$. This means that
$A_2^{-1}A_1\in K$ and so $Gr(p,\BF^n)\cong G/K$.

{\bf To describe $\mbox{Affine}~ Gr(p,\BF^n)$}, notice
that a plane $A\MAT{1}I_p\\O\mat$ for which $\det
A_{11}\neq 0$, can be expressed as $$
G/K\ni~\mbox{span}~A\MAT{1}I_p\\O\mat=
 ~\mbox{span}\MAT{1}A_{11}\\A_{21}\mat=
 ~\mbox{span} \MAT{1}A_{11}\\A_{21}\mat A^{-1}_{11}=~\mbox{span}
\MAT{1}I_p\\Z(A)\mat , $$ where
$Z(A):=A_{21}A_{11}^{-1}$ is a $(n-p)\times p$ matrix.
Also notice that $Z(A)$ is unchanged upon multiplying
$A$ to the right with $B\in K$.

The {\bf left action} of $B\in K$ on $A\in G$ has the
following effect on
\be
 \mbox{Affine~}G/K\lrg\mbox{~Affine~}
 G/K:\MAT{1}I_p\\Z\mat\curvearrowright
\MAT{1}I_p\\B_{22}ZB^{-1}_{11}\mat ,
 \ee
  because in $$
BA\MAT{1}I_p\\O\mat = \MAT{2}B_{11}&O\\
O&B_{22}\mat\MAT{1}A_{11}\\A_{21}\mat=
\MAT{1}B_{11}A_{11}\\ B_{22}A_{21}\mat, $$ we have
$Z(BA)=(B_{22}\,\,\,A_{21})(B_{11}
 \,\,\,A_{11})^{-1}=B_{22}Z(A)B_{11}^{-1}$.
Picking arbitrary matrices $B_{11}\in K_1$, $B_{22}\in
K_2$, the $(n-p)\times p$ ($n\geq 2p$) matrix $Z(A)$
can be ``diagonalized", namely
\be
 Z(BA)=B_{22}Z(A)B_{11}^{-1}=\left(
\begin{array}{ccc}
\al_1 & &O      \\
 &\ddots&      \\
O & &\al_p   \\
\hline
  & O&
\end{array}
\right).\ee

Now we use the fact that the $p\times p$ matrix
$Z^{\dag}Z$ is ``self-adjoint" and positive
definite\footnote{since
$v^{\dag}Z^{\dag}Zv=(Zv)^{\dag}Zv=
\displaystyle{\sum_1^p}|(Zv)_i|^2>0$, for $v\in
\BF^p \backslash 0$.}, and so, setting $\al_i=:-\tan\theta_i$,
 \be
(Z(BA))^{\dag}Z(BA)=\left(
\begin{array}{ccc}
\al_1^2 & &O      \\
 &\ddots&      \\
O & &\al_p^2
\end{array}
\right)=\left(
\begin{array}{ccc}
\tan^2\theta_1 & &O      \\
 &\ddots&      \\
O & &\tan^2\theta_p
\end{array}
\right).
\ee
Therefore, by the left action of $K$ on $G$,
 the $p$-plane in $\BF^n$ can be
represented by the span of the columns of the
following matrix, which by taking the linear
combination of the columns, each multiplied with
$\cos\theta_i$ reads
 $$ \mbox{span}\left(
\begin{array}{ccc}
&I_p& \\
\hline
-\tan ~\theta_1 & &O      \\
 &\ddots&      \\
O & &-\tan~\theta_p\\
\hline
&O_{n-2p,p}&
\end{array}
\right)=\mbox{span}\left(
\begin{array}{ccc}
\cos\theta_1& & O\\
&\ddots& \\
O& &\cos\theta_p\\
\hline
-\sin \theta_1 & &O      \\
 &\ddots&      \\
O & &-\sin \theta_p\\
\hline
&O_{n-2p,p}&
\end{array}
\right),
$$
and so
$$
 K\backslash G /K\simeq
\left\{
%\left( {I_m}\atop {Z} \right)~ \Bigr|~ Z=
%
\mbox{span }
 \MAT{1}\Re e^{i \Theta_p}\\ -\Im e^{i
\Theta_p}
\\O_{n-2p,p}\mat
 ~,~0 \leq \theta_i \leq {\pi}  \right\}
 . $$

Similarly, $Z\in T_{Id}Gr(p,\BF^n)$ can be
diagonalized by the action (1.0.16) on
$T_{Id}Gr(p,\BF^n)$; i.e.,
\be
 Z \mapsto B_{22}Z B_{11}^{-1}=\left(
\begin{array}{ccc}
z_1 & &O      \\
 &\ddots&      \\
O & &z_p   \\
\hline
  & O&
\end{array}
\right)\in {\frak a}^+,\ee where the $z_i$ are
linearized versions of the $\tan \theta_i$'s and where
${\frak a}^+$ is a fixed Weyl chamber in the Cartan of
${\frak p}.$

By the Weyl integration formula, the measure induced
on $K\backslash G/K$, via the Haar measure and the
embedding $G/K\hookrightarrow G$,  is given by (in the
compact case) (see Helgason \cite{Helgason1}, p.~188)
 $$
 d\mu(H)=\prod_{\al\in\Sg^+}
  |\sin\al(iH)|^{m_{\al}}d{H},
 ~~~
~H\in {\frak
a}_*, $$
  where $g=\mathfrak{k}+\mathfrak{p}$, with compact
  real form $u=\mathfrak{k}+\mathfrak{p}_*,~
  \mathfrak{p}_*=i
  \mathfrak{p}$ and
 with ${\frak a}_*$ a maximal abelian subspace of
$\mathfrak{p}$,
where $\Sg^+$ is the set of roots having
positive values on the fixed Weyl chamber ${\frak
a}^+$ of ${\frak
a}_*$ and
where the root space $g_{\al}$ has
dimension $m_{\al}$ for any restricted root $\al$: $$
g_{\al}:=\{X\in g\mid ~ [H,X]=\al(H)X,\mbox{~for all
$H\in{\frak a}\}$},
 $$
  and we also have the induced
measure on $T_{Id}(K\backslash G/K)$
(\cite{Helgason1}, p.~195)
$$ d\mu=\prod_{\al\in\Sg^+}\al(H)^{m_{\al}}d{ H},\qquad
H\in{\frak a}^+. $$ The roots and multiplicities
$m_{\al}$ are as follows ($n\geq 2p$):

\medbreak

$$
\begin{array}{c|c|c|c}
Gr(p,\BC^n)&Gr(p,\BR^n)&Gr(p,\BH^n)&\al\in\Sg^+\\ \hline
2&1&4& i(\varepsilon_j+\varepsilon_k) \\
2&1&4& i(\varepsilon_j-\varepsilon_k)\\ 1&0&3&%\\
2i\varepsilon_{\ell}\\ 2(n-2p)&n-2p&4(n-2p)&%\\
i{\varepsilon_{\ell}~}
\end{array}
$$ with $1\leq j<k\leq p$, $1\leq\ell\leq p$, yielding
Table 1, upon setting
\be
 H=i(\theta_1,\ldots,\theta_p),
  0\leq\theta_i\leq\pi,
 \ee
e.g., we check table 1 for $Gr(p,\BH^n)$.
Setting $k=n-2p$ and $y_j=\cos 2\theta_j%=2z_j-1
$, which is very natural in view of (1.0.20):

\medbreak \noindent \underline{ for $K\backslash
G/K$:} \bea
d\mu&=&\prod_{\al\in\Sg^+}|\sin\al(iH)|^{m_{\al}}dH
\nonumber\\ &=&\prod_{1\leq j<k\leq p}|\sin
(\theta_j-\theta_k)\sin
 (\theta_j+\theta_k)|^4\prod^p_{j=1}|
 \sin 2\theta_j|^3|
\sin^2 {\theta_j}|^{2(n-2p)}d\theta_j\nonumber\\
&=&\prod_{1\leq j<k\leq p}|\frac{1}{2}(\cos
2\theta_j-\cos 2\theta_k)|^4
 \prod^p_{j=1}(1-\cos^2 2\theta_j)
\left(\frac{1-\cos 2\theta_j}{2}\right)^{2k}\frac{1}{2}
 d\cos 2\theta_j\nonumber\\
&=&2^{-
p(2(p+k)-1)}\Delta(y)^4\prod^p_{j=1}(1+y_j)(1-y_j)^{2k+1}dy_j
\eea

\medbreak \noindent \underline{for $T_{id}K\backslash
G/K$:}
 \bean
d\mu&=&\prod_{\al\in\Sg^+}\al(H)^{m_{\al}}dH\\
&=&c\prod_{1\leq j<k\leq
p}(v_j-v_k)^4(v_j+v_k)^4\prod^p_{j=1}v_j^{4k+3}dv_j\\
&=&c\prod_{1\leq j<k\leq p} (v_j^2-v^2_k)^4
\prod^p_{j=1}(v^2_j)^{2k+1}\frac{1}{2}dv^2_{j}\\
&=&c2^{-p}\Delta^4(u)\prod^p_{j=1}u_j^{2k+1}du_j,
  \eean   \vspace{-1.9cm}\be \ee
setting $u_j=v_j^2$, where in the above we made the
identification (1.0.6) of Proposition 1.3 and put $Z$
in the normal form (1.0.16), so that $Z^{\dag}Z= \diag
(u_1,\ldots , u_p).$

To describe $G/K$ in a second way, remember
 $g=\mathfrak{k}+\mathfrak{p}$, with
 $\mathfrak{k},\mathfrak{p}$ the $\pm$ eigenspaces
 of a lie algebra involution $\sigma$. The latter lifts
 to the group as an involution $\sigma$,
 which commutes with inversion, i.e., $(g^{\sigma})^{-1}
 = (g^{-1})^{\sigma}$. Use $\sigma$ to define
  the following embedding
 \be {\iota }: G \hookrightarrow G: g\mapsto \iota
 (g)= g(g^{\sigma})^{-1},\ee which induces a natural injective
  map
 $$ {\iota }: G/K \hookrightarrow G: g\mapsto \iota
 (g)%= g(g^{\sigma})^{-1}
 .
  $$
   Indeed, $\iota (g_1)=\iota (g_2)$
   is equivalent to
   $g_1(g_1^{\sigma})^{-1}=g_2(g_2^{\sigma})^{-1}$,
   which amounts to
   $ g_2^{-1}g_1=(g_2^{\sigma})^{-1}g_1^{\sigma}
   = (g_2^{-1})^{\sigma}g_1^{\sigma}
   = (g_2^{-1}g_1)^{\sigma}$, meaning that
   $g_2^{-1}g_1\in K$.

From the polar decomposition $G=K{\cal A}K$, we have
that every $g \in G$ can be decomposed into
 $$
g=k_1^{-1} A(\theta) k_2^{-1}, ~~k_1,k_2 \in K,
~A(\theta) \in {\cal A},
 $$
 with $\theta$ the torus coordinates, such that
 $A(\theta) A(\theta')=A(\theta+\theta'), ~A(0)=I$ and
$A^{-1}(\theta)=A(-\theta)$. Since $A=\mbox{exp}~ a$,
with $a \in i\mathfrak p$, we have
$A(\theta)^{\sigma}=A^{-1}(\theta)$. From (1.0.22) the
torus $\cal A$ embeds into $\cal S$ as follows:
 $$
   \iota (A(\theta))=
    A(\theta) ( A^{\sigma}(\theta))^{-1}=
    A(\theta) A(\theta)=A(2\theta).
    $$

Moreover the polar decomposition $ g=k_1^{-1}
A(\theta) k_2^{-1} $ yields conjugation in $\cal S$ by
$k_1$ :
\bea
  \iota(g)
    &=& k_1^{-1}( k_1 g k_2 k_2^{-1} ( g^{\sigma})^{-1}
    k_1^{-1})  k_1  \nonumber\\
   &=& k_1^{-1}( k_1 g k_2 (k^{\sigma}_2)^{-1}
    ( g^{\sigma})^{-1}
    (k^{\sigma}_1)^{-1})  k_1  \nonumber\\
    &=& k_1^{-1}\left( k_1 g k_2
    ((k_1 g k_2)^{\sigma})^{-1}\right)
       k_1
    \nonumber\\
    &=& k_1^{-1} \iota (k_1 g k_2)    k_1
    \nonumber\\
    &=& k_1^{-1} \iota (A(\theta))    k_1
    \nonumber\\
    &=& k_1^{-1}  A(2\theta)   k_1
    \eea

 Specializing to the Grassmannian case, we have the involution $$ \sigma:G
\rightarrow G:  g\mapsto   g^{\sigma}:=
  I_{p,n-p}~ g~ I_{p,n-p}
$$

with
\bean
K&=&\{\mbox{fixed points in $G$ of the involution
$\sigma$} \}\\
&=&
 \{ g \in G \mbox{ such that}~ g~ I_{p,n-p}=
 I_{p,n-p}~g\}.
\eean

Setting
$$
\Theta_p:= \diag (\theta_1,\ldots, \theta_p),
$$
we have that a maximal abelian subspace of the Lie
algebra is given by
$$
{\frak a} =
 \left\{
a(\theta) =
 \left(\begin{array}{cccc}
 O_p &  & \Theta_p  & \\
  &&& O\\
 -\Theta_p  & & O_p  &  \\
  & O & & I_{n-2p} \\
  \end{array}
\right),~~ 0\leq \theta < 2\pi
\right\}
 $$
and by exponentiation, we find the torus
$${\frak A} =\left\{
A(\theta) =
 \left(\begin{array}{cccc}
 \Re e^{i\Theta_p} & & \Im e^{i\Theta_p} &  \\
                  & &                  & O\\
 -\Im e^{i\Theta_p} && \Re e^{i\Theta_p} & \\
  & O & & I_{n-2p} \\
\end{array}
\right),~~ 0\leq \theta < 2\pi \right\} \in G. $$
%with $A(\theta) A(\theta')=A(\theta+\theta')$ and
%$A(\theta)^{\sigma}=A(-\theta)$.
The spectrum of
$A(\theta)$ is easily seen to be
\be
(e^{i\theta_1},\ldots, e^{i\theta_p},
e^{-i\theta_1},\dots, e^{-i\theta_p},
 \underbrace{1,\ldots,1}_{n-2p})
.\ee

  To connect with the previous description, given $g \in G$,
  we pick $k_1,k_2$ such that $k_1 gk_2=A(\theta)$ and by the
  previous discussions, we have
\bean
 Z^{\dag}(g)Z(g)
 & =&
Z^{\dag}(k_1gk_2)Z(k_1gk_2)  \\
 &=&
Z^{\dag}(A(\theta))Z(A(\theta)) \\
  &=& \mbox{diag}~(\tan^2 \theta_1,\ldots, \tan^2
 \theta_p).
 \eean
  The embedding ${\frak A}\hookrightarrow Gr(p,\BF^n)$ is
then given by $${\frak A}\hookrightarrow Gr(p,\BF^n):
A(\theta)\mapsto A(\theta) \MAT{1}I_p\\O_{n-p,p}\mat=
\MAT{1}\Re e^{i \Theta_p}\\ -\Im e^{i \Theta_p}
\\O_{n-2p,p}\mat,
$$
 and so the $\theta_i$'s in the two discussions are
 identical.
\qed

\bigbreak

\noindent{\it Proof of Theorems 1.1, 1.2 and 1.3:}
 In order to compute integral (1.0.3), recall in the above description
of $K\backslash G/K$,
$$H=i(\theta_1,\ldots,\theta_p)%=2i(\lb_1,\ldots,\lb_p).
.$$ Since the integrand is invariant under the left
action of $K$ on $G$, which induces on $Z$ the map
$Z\rg B_{22}ZB^{-1}_{11}$, use (1.0.17), from which it
follows
 that, upon using $\tan^2\theta_i=
  \frac{1-\cos2\theta_i}{1+\cos2\theta_i}=\frac{1-y_i}{1+y_i}$,
and $\cos^2\theta_i=(1+y_i)/2$,% we may set $$
%Z^{\dag}Z=\MAT{3}\tan^2\lb_1& &O\\
% &\ddots& \\
%O& &\tan^2\lb_p\mat, \mbox{~where
%$\tan^2\lb_i=\frac{1-z_i}{z}$}$$ and thus
we have
$$
 Z^{\dag}Z=\diag
 (\tan^2\theta_1,\ldots,\tan^2\theta_p)
 =\diag \left( \frac{1-y_1}{1+y_1},\ldots,
  \frac{1-y_p}{1+y_p} \right)
 $$
and
$$
(I+Z^{\dag}Z)^{-1}=\left(
\begin{array}{ccc}
\cos^2\theta_1& & O\\ &\ddots& \\ O& &\cos^2\theta_p
\end{array}
\right) =\frac{1}{2}\diag \left( 1+y_1,\ldots, 1+y_p
\right). $$
 Hence, setting $z_i=\cos^2\theta_i$ in the last identity below, we have
 $y_i=2z_i-1$, using (1.0.11) and picking an
 appropriate normalizing constant $c$:

\medbreak

$\displaystyle{\int_{Gr(p,\BF^n)}e^{x\Tr(I+Z^{\dag}Z)^{-1}}\det
(Z^{\dag}Z)^{-\beta(q-p)}d\mu(Z)}$
 \bean
&=& c \int_{[-1,1]^p}e^{\frac{x}{2}
 \sum_1^p(1+y_i)}|\Delta_p(y)|^{2\beta}\prod^p_
{i=1} \left(\frac{1-y_i}{1+y_i}\right)^{\beta(p-q)}
(1-y_i)^{\beta (n-2p)}
 (1-y_i^2)^{\beta-1}dy_i
% (1-z_i)^{\beta(n-2p)+\beta
%-1} z_i^{\beta -1}dz_i
\\
 & & \\
&=&c\int_{[0,1]^p}e^{x
\sum_1^pz_i}|\Delta_p(z)|^{2\beta}\prod_1^p
(1-z_i)^{\beta(n-p-q+1) -1}z_i^{\beta(q-p+1)-1}dz_i,
 \\ &&
\hspace{7cm}~~\mbox{upon setting} ~ y_i=2z_i-1.
\eean
 Finally, using the Weyl integration formula (1.0.21) on
the tangent space and identifying the
$(z_1,\ldots,z_p)$ in $(1.0.18)$ with a point in
${\frak a}^+$, and setting $u_i=z_i^2\geq 0$, $1\leq
i\leq p$, we find $$ \int_{{{Z\in
T_{id}Gr(p,\BF^n)}\atop{\mbox{\tiny with
spectrum~}(Z^{\dag}Z)\in E}}}e^{-x\Tr Z^{\dag}Z}
\det(Z^{\dag}Z)^{\beta p}d\mu(Z),$$ $$=
\int_{E^p}|\Delta _p(u)|^{2\beta} \prod^p_1e^{-x
u_i}u_i^{\beta(n-p+1)-1}du_i, $$ establishing Theorem
1.3.

Finally, to prove Theorem 1.2, notice that, according
to (1.0.23), matrices $M$ in ${\cal S}$ can be
diagonalized to matrices $A(2\theta)$, with spectrum
\be
(e^{2i\theta_1},\ldots, e^{2i\theta_p},
e^{-2i\theta_1},\dots, e^{-2i\theta_p},
 \underbrace{1,\ldots,1}_{n-2p})
.\ee
So $M$ decomposes into $M=M_0\otimes M_1$, with $M_1$
being the $n-2p$-dimensional eigenspace corresponding
 to the eigenvalue $1$ and so
$$
\prod_1^p\frac{1-y_k}{1+y_k}=\prod_1^p
 \frac{1-\cos 2\theta_k}
 {1+\cos 2\theta_k}=
 \prod_1^p\frac{(1-e^{2i\theta_k})(1-e^{-2i\theta_k})}
 {(1+e^{2i\theta_k}) (1+e^{-2i\theta_k})}
 =\det \frac{I-M_0}{I+M_0}.
$$
\qed

%\newpage

%\newpage

\section{Jack polynomials}

\subsection{Young diagrams and Schur polynomials}

 Standard references to this
subject are MacDonald, Sagan, Stanley, Stanton and
White \cite{MacDonald,Sagan, Stanley, Stanton}. To set
the notation, we remind the reader of a few basic
facts.

\begin{itemize}
  \item A A {\em partition } of
  $n=|\lambda|:=\lambda_1+...+\lambda_{\ell}$ (with
  $n=|\lambda|$ called the weight) is represented by
  a {\em Young diagram} $\lambda_1\geq
  \lambda_2 \geq ...\geq \lambda_{\ell}\geq 0$.
 A {\em dual Young diagram} $ \lambda^{\top}=(
  \lambda^{\top}_1 \geq \lambda^{\top}_2 \geq...)$ is the
  diagram obtained by flipping the diagram $\lambda$
  about its diagonal.

  \item A {\em semi-standard Young tableau} of shape $\lambda$
  %$\lambda_1\geq  ...\geq \lambda_{\ell}$
  is an array of positive
   integers $a_{ij}$ placed at $(i,j)$ in the Young diagram $\lambda$, which are
   non-decreasing from left to right {\em and} strictly
   increasing from top to bottom.

  \item A {\em standard Young tableau} of shape $\lambda$ is an array of
   integers $1,...,n$ placed in the Young diagram, which are
   strictly
   increasing from left to right {\em and} from top to bottom.

 \item The {\em Schur polynomial $s_{\lambda}$} associated with a Young
  diagram $\lambda$ is a symmetric function in the
  variables $x_1,x_2,...$, (finite or infinite), where
  $n=|\lb|$ and
  defined by (for notation $f^{\lb}$, see the next
  point)
  $$
  s_{\lambda}(x_1,x_2,...)=\sum_{{\{a_{i,j}\}~
  \mbox{\tiny{semi-standard}} }\atop
  {\mbox{\tiny{tableaux
  }}\lambda}}~\prod_{(i,j)\in \lb}x_{a_{i,j}}
  = f^{\lb}~ x_1\ldots x_n+...$$

\item The {\em hook length} of the $i,j$th box is defined by
$h^{\lb}_{ij}:=\lb_i + \lb^{\top}_j-i-j+1$. Also
define
  \bea h^{\lb}&:=&\prod_{(i,j)\in
\lb}h^{\lb}_{ij} \nonumber\\ &=& \frac{\prod_1^m
(m+\lb_i-i)!}{\Dt_m(m+\lb_1-1,\ldots,m+\lb_m-m)},~~\mbox{for}~
m \geq \lb_1^{\top}. \eea

%\newpage

\item   The {\em number of standard Young
   tableaux} of a given shape $\lambda=
   (\lambda_1\geq...\geq\lambda_m)$ is given
   by
  \bea f^{\lb}&=&\# \{\mbox{standard tableaux of shape
$\lb$}\}\nonumber\\
&=&  |\lb|!~s_{\lb}(x)\Bigr|_{\sum_{k}x_k=1 \atop
 \sum_{k}x^i_k=0~\mbox{{\tiny for $i\geq 2$}}}
\nonumber\\
&=&\mbox{coefficient of $x_1x_2\ldots x_n$ in
$s_{\lb}(x)$} \nonumber
\\&=&\frac{|\lb
|!}{h^{\lb}}~=~ |\lb |
!\det\left(\frac{1}{(\lb_{i}-i+j)!}\right)\quad
\nonumber\\
%|\lb|=k
  &=&  |\lb|! ~
  \frac{\Delta_{m}(m+\lb_1-1,\ldots,m
 +\lb_{m}-m)}
{\displaystyle{\prod_1^{m}} (m+\lb_i-i)!},
\quad\mbox{for $m\geq\lb^{\top}_{1}$}.\nonumber \\
% &=& |\lb |
%! \prod_{1\leq i< j \leq m}(h_i-h_j)\prod_1^m
%\frac{1}{h_i!},~~\mbox{with}~h_i:=\lambda_i-i+m,~m:=\hat
%\lambda_1 \nonumber
 . \eea

\item   The {\em number of semi-standard Young
   tableaux} of a given shape $\lambda$, with numbers $1$ to $k$
for $k\geq 1$:
  \bea
 \lefteqn{ \# \left\{\begin{array}{l}\mbox{semi-standard
tableaux of shape $\lb$} \\ \mbox{filled with numbers
from $1$ to $k$}
\end{array}
\right\}}\nonumber\\
 &=&
s_{\lb}(\overbrace{1,\ldots,1}^k,0,0,\ldots)
\nonumber\\
 &=&
\prod_{(i,j)\in\lb} \frac{j-i+k}{h^{\lb}_{i,j}}
\nonumber\\
&=& \left\{ \begin{array}{l}
  \displaystyle{
\frac{\Delta_k(k+\lb_1-1,\ldots,k+\lb_k-k)}
{\displaystyle{\prod^{k-1}_{i=1}}i!},~~~\mbox{when $ k
\geq \lb_1^{\top}$,}}\\
\\
 \displaystyle{   0},~~~\mbox{when $ k
< \lb_1^{\top}$,}
 \end{array} \right.\nonumber\\
 \nonumber\\
 \eea
 using the fact that
\be
\prod_{(i,j)\in\lb} (j-i+k)=\frac {\prod_{i=1}^k
(k+\lb_i-i)!}{\prod_1^{k-1}i!}. \ee

  \item

\noindent{\em Robinson-Schensted-Knuth
correspondence}: Given
\bean S_n&=&\mbox{group of permutations of $\{1,\ldots
n \}$}\\ S_n^k&=&
 \{\mbox{words of length $n$ built from the set
 $\{1,\ldots,k\}$}\} \eean
the following 1-1 correspondences hold:
 $$
S_n\longrightarrow\left\{\begin{array}{l} (P,Q),
\mbox{where $P$ and $Q$ are two standard Young
tableaux}\\ \mbox{of same shape $\lb$, with $|\lb|=n$
and taken from $\{1,\ldots,n\}$ }
\end{array}
\right\} $$

 $$
S_n^k\longrightarrow\left\{\begin{array}{l} (P,Q),
\mbox{where $P$ and $Q$ have same shape $\lb$, with
$|\lb|=n$}\\ \mbox{$P$ is semi-standard, filled with
numbers from $1$ to $k$, and} \\ \mbox{$Q$ is
standard, filled with numbers from $1$ to $n$}
\end{array}
\right\} $$

It follows that for given $n$ and $k$, we have

\bea \sum_{{\lb\mbox{\tiny with}}\atop {
|\lb|=n}}\left(f^{\lb}\right)^2&=&n!  \\
  \sum_{{\lb\mbox{\tiny with}}\atop {
|\lb|=n}}f^{\lb}~ s_{\lb}(1^k)
% \prod_{(i,j)\in\lb}\frac{j-i+k}{h^{\lb}_{i,j}})
&=&k^n
 \eea

\item  {\em Increasing and decreasing sequences}

According to Greene \cite{Greene}, given a word
$\pi\in S_n^k$, mapped, via the RSK correspondence,
into $(P,Q)$ of shape
$\lambda=(\lambda_1,\ldots,\lambda_{\ell})$, then for
any $k$, \bean
 \lambda_1+\ldots+\lambda_k&=&
 \left\{\begin{array}{l}
\mbox{length of the
 longest {\em weakly}}\\ \mbox{ $k$-increasing subsequence }
\end{array}
\right\}\\
\lambda^{\top}_1+\ldots+\lambda^{\top}_k&=&
\left\{\begin{array}{l} \mbox{length of the
 longest {\em strictly}}\\\mbox{ $k$-decreasing subsequence}
\end{array}
\right\}
  \eean
  \end{itemize}

\subsection{Some useful formulae on hook length}

Remembering the notation (0.0.8), we have the
following statement:

 $$ \squaresize .4cm \thickness
.01cm \Thickness .07cm \Young{
 &&&&&&&&r&&&&&&&&&&&&&\cr
 &&&&&&&&r&&&&&&&&&&&\cr
 &&&&&&&&r&&&&&&&&&\cr
 &&&&\mu&&&&r&&&&&&&&\cr
 &&&&&&&&r&&&&&&&&             \cr
 &&&&&&&&r&&&&                 \cr
 &&&&&&&&r                   \cr }
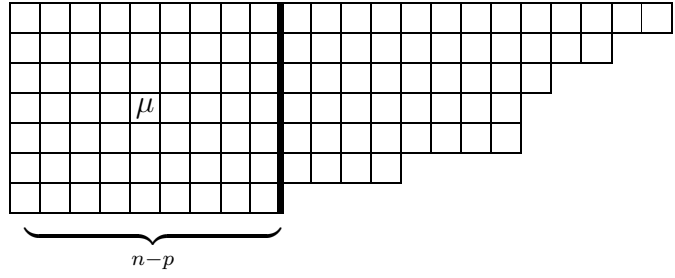
$$

\vspace{-1cm}

 $$\hspace{.2cm}\underbrace{\hspace{3.4cm}}_{n-p}
%\hspace{.01cm} ~~~
%~\underbrace{\hspace{1.5cm}}_{q-p}
\hspace{5.2cm}$$

\begin{lemma} Given a partition $\lb\supseteq
\mu :=(n-p)^p$ with $\lb_1^{\top}=p$, $0\leq p<n$,
then\footnote{Here $(a)_{\lb}:=
(a)_{\lb}^{(1)}=\prod_i (a+1-i)_{\lb_i}$.}
 \bea
h^{\lb}&=&h^{\mu}h^{\lb\backslash\mu}
  \frac{(n)_{\lb\backslash\mu}}{(p)_{\lb\backslash\mu}}
%\prod^p_{i=1}\frac{(n-i+1)_{(\lb\backslash\mu
%)_i}}{ (p-i+1)_{(\lb\backslash\mu)_i}}
,\\
h^{\mu}&=&\frac{\displaystyle{\prod^p_{1}}(n-i)!}{\displaystyle{\prod^{p-1}_{1}}i!
}\mbox{~~and~~}s_{\mu}(1^p)=1, \\
\frac{s_{\lb}(1^p)}{h^{\lb}}&=&
 \frac{1}{h^{\mu}(h^{\lb\backslash\mu})^2}
\frac{((p)_{\lb\backslash\mu})^2}{(n)_{\lb\backslash\mu}}
% \prod^p_{i
%=1}\frac{\left((p-i+1)_{(\lb\backslash\mu)_i}
%\right)^2}{(n-i+1)_{(\lb\backslash\mu)_i}}
 . \eea
\end{lemma}

\proof Setting $\kappa =\lb\backslash\mu$ \bean
\frac{h^{\lb}}{h^{\lb\backslash\mu}}
 &=&\prod_{(i,j)\in\mu}h^{\lb}_{(i,j)}
 \\
&=&\prod^p_{i=1}\prod^{n-p}_{j=1}(n+1+\kappa_i-i-j)
 \\
&=&\prod^p_{i=1}(n+\kappa_i-i)
\ldots(n+1+\kappa_i-i-n+p)
\\
&=&h^{\mu}\prod^p_{i=1}\frac{(p+1+\kappa_i-i)\ldots(n+\kappa_i-i)}{(p+1-i)\ldots
(n-i)}\\
&=&
h^{\mu}\prod^p_{i=1}
 \frac{(n-i+1)_{\kappa_i}}{(p-i+1)_{\kappa_i}}
\eean

\bea h^{\mu}&=&[(n-1)(n-2)\ldots
p][(n-2)\ldots(p-1)]\ldots[(n-p)\ldots 1]\nonumber\\ &
&\nonumber\\ &=&\frac{(n-1)!(n-2)!\ldots
(n-p)!}{(p-1)!(p-2)!\ldots 1!}\nonumber\\ &
&\nonumber\\
&=&\frac{\displaystyle{\prod^p_1}(n-i)!}{\displaystyle{\prod^{p-1}_1}i!}.
\eea Then $$
s_{\mu}(1^p)=\frac{\displaystyle{\prod_{(i,j)\in\mu}}
 (j-i+p)}{h^{\mu}}=
\frac{\displaystyle{\prod_1^p}(n-i)!}{h^{\mu}
 \displaystyle{\prod_1^{p-1}}i!}=1,
  \quad\mbox{using (2.1.4).}
   $$
%\newpage

    Finally, using in the second identity below formula (2.1.4) and in the fourth identity
    $\prod^p_1
(n-i+\kappa_i)!
 =
 \prod^p_1(n-i)!\prod^p_1(n-i+1)_{\kappa_i}$, one
 computes
  \bean
\frac{s_{\lb}(1^p)}{h^{\lb}}
 &=&
   \frac{\displaystyle{\prod_{(i,j)\in\lb}}(j-i+p)}{
(h^{\lb})^2}
 \\
 &=&\frac{1}{(h^{\lb})^2}
 \frac{\displaystyle{\prod_{i=1}^p}
(p+\lb_i-i)!}{\displaystyle{\prod_1^{p-1}}i!}
 \\
&=&\frac{1}{(h^{\lb})^2}
 \frac{\displaystyle{\prod_{i=1}^p}
(n-i+\kappa_i)!}{\displaystyle{\prod_1^{p-1}}i!},
 \quad\mbox{since $\lb_i=n-p+\kappa_i$}
%  \\
\eean \bean
&=&\frac{1}{(h^{\lb})^2}\frac{\displaystyle{\prod_{1}^p}
(n-i)!}{\displaystyle{\prod_1^{p-1}}i!}\prod^p_{i=1}
(n-i+1)_{\kappa_i}
  \\
&=&\frac{h^{\mu}}{(h^{\mu}h^{\lb\backslash\mu})^2}\prod^p_{i=1}
 \frac{\left((p-i+1)_{({\lb\backslash\mu})_i}\right)^2}
  {(n -i+1)_{({\lb\backslash\mu})_i}},\quad\mbox{using
(2.2.1) and (2.2.2)} \eean ending the proof of Lemma
2.1.\qed

%\newpage

\subsection{Jack polynomials}

Define symmetric polynomials
\bean
p_{\lb}(x_1,x_2,\ldots )&:=& p_{\lb_1}p_{\lb_2}\cdots
 =\sum_i x_i^{\lb_1}
  \sum_i x_i^{\lb_2} \cdots
  \\
  m_{\lb}(x_1,x_2,\ldots )&:=&
  \sum_{\mbox{\small permutations}\atop \mbox{of
  $x_1,...,x_n$}} x_1^{\lb_1} \ldots x_n^{\lb_n},
  \eean
and the dominance ordering between partitions:
$$\mu \leq \lb ~~\mbox{means} ~:~~~\sum_1^{\ell}\mu_i\leq
 \sum_1^{\ell}\lb_i,~~\mbox{for all $\ell$} .$$
 Given that $\lb$ has $m_i=m_i(\lb)$ parts equal to $i$,
define the inner-product $\langle \, , \,\rangle$ on
the vector space of all symmetric functions of bounded
degree (to be explained)
 $$ \la p_{\lb},p_{\mu}\ra
=\delta_{\lb\mu}(1^{m_1}2^{m_2}\ldots)m_1!m_2!\ldots
 \al^{ \lb_1^{\top}}.
$$

\noindent Jack polynomials are the unique symmetric functions
$J_{\lb}^{(\al)}$ satisfying

\medbreak

(i)\quad $\la J_{\lb}^{(\al)},J_{\mu}^{(\al)}\ra =0$, if $\lb\neq\mu$
,
\medbreak

(ii)\quad
$J_{\lb}^{(\al)}=\sum_{\mu\leq
\lb}v_{\lb\mu}(\al)m_{\mu}$,

\medbreak

(iii)\quad If $\vert\lb \vert =n$, then
$$
J_{\lb}^{(\al)}=n!x_1\ldots x_n+\ldots
$$

It follows that
 $$v_{\lb\lb}=
  \prod_{(i,j)\in\lb}
% \left(\lb^{\top}_j-i+\al(\lb_{i}-j+1)\right)
\left(\lb^{\top}_j-i+1+\al(\lb_i-j)\right).
  $$

\subsubsection*{Special cases:}

{\em Jack polynomials} for $\al=1$ are proportional to Schur polynomials,
namely
$$
J_{\lb}^{(1)}=h^{\lb}s_{\lb}.
$$

\noindent{\em Zonal polynomials} are given by
 $$Z_{\lb}^{( \beta)}=J_{\lb}^{(1/\beta)}~,~~\mbox{
with $\beta = 1,1/2,2 $}.$$ They have the remarkable
property that for $G=O(n), ~U(n)$ or $U(n,\BH):=
\{g~|~ g\bar g^{\top}=I\}$,
 \be
\int_{G}J^{(\alpha)}_{\lb}(\sigma k\tau
k^{-1})dk=\frac{J^{(\alpha)}_{\lb}(\sigma)
 J^{(\alpha)}_{\lb}(\tau)}
{J^{(\alpha)}_{\lb}(1^n)}  , \ee
 for all
$$
\begin{array}{lll}
\sigma, \tau\in  \Sigma&= \{ \mbox{real symmetric
matrices}\}
           & ~\mbox{for} ~\alpha=2  \\
       &= \{ \mbox{ Hermitian matrices}\}
            & ~\mbox{for} ~\alpha=1  \\
       &= \{ \mbox{quaternionic matrices, with}~
             \sigma =\bar \sigma^{\top}\}
            & ~\mbox{for} ~\alpha=1/2. \\
  \end{array}
 $$
The function $J^{(\alpha)}_{\lb}(\tau)$ is a symmetric
function of the (real) spectrum of $\tau$.

\subsubsection*{Orthogonality:}
\be
\la J_{\lb}^{(\al)},J_{\mu}^{(\al)}\ra
=\dt_{\mu\lb}j_{\lb}^{(\al)}, \ee where
 \bea
j_{\lb}^{(\al)}&=&\prod_{(i,j)\in\lb}
 \left(\lb^{\top}_j-i+\al(\lb_{i}-j
+1)\right)\left(\lb^{\top}_j-i+1+\al(\lb_i-j)\right)
\nonumber\\&&\nonumber\\ &=&\left\{
 \begin{array}{ll}
  (h^{\lb})^2&\mbox{~for~}\al=1 \\
  h^{2\lb} & \mbox{~for~}\al=2\\
   h^{2\lb^{\top}}/2^{2|\lb|}&\mbox{~for~}\al=1/2.
  \end{array}\right.
   \eea

\subsubsection*{Special values:}
For arbitrary $n$, we have
\bea
J_{\lb}^{(\al)}(1^n)&=&\prod_{(i,j)\in\lb}
 \left(n-(i-1)+\al(j-1)
\right),
\mbox{~where~}1^n=(\overbrace{1,\ldots,1}^n,0,0,\ldots)
 \nonumber\\
&=&\al^{|\lb|}\prod_{(i,j)\in\lb}\left(\frac{1}{\al}
(n-i+1)+j-1\right)=\al^{|\lb|} \left( \frac{n}{\al}
\right)^{(\al)}_{\lb} \nonumber\\ &=&
   \left\{ \begin{array}{l}
\displaystyle{ \al^{|\lb|}\prod_{i=1}^m
 \frac{\Gamma\left(\frac{1}{\al}(n-i+1)+
\lb_i\right)}{\Gamma\left(\frac{1}{\al}(n-i+1)\right)}
 > 0~,} ~~\mbox{for all $m\geq \lb_1^{\top}$, if $n \geq\lb_1^{\top}$},\\
   \\
 0  ~, ~~\mbox{if $n <\lb_1^{\top}$}
 \end{array} \right.
% \nonumber\\ & &
  \eea
    and so, for $\al=1$,
  \bea
s_{\lb}(1^n)&=&\frac{1}{h^{\lb}}J_{\lb}^{1}(1^n)
 \nonumber\\
&=&
 %\frac{1}{h^{\lb}}\prod^m_{i=1}\frac{(m-i+\lb_i)!}{(m-i)!}
\left\{\begin{array}{l}\displaystyle{\frac{1}{h^{\lb}}
 \prod^m_{i=1}\frac{(n-i+\lb_i)!}
{ (n-i)!}\geq 0~,~~\mbox{~for all~}m \geq \lb^{\top}_1},~
\mbox{~if~}n \geq \lb^{\top}_1,\\
\\
0 ~,~~\mbox{~if~}n < \lb^{\top}_1.
 \end{array}  \right.
   \eea
  The last identity in (2.3.4) is
  obtained by taking the product over the $i$th
  row of $\lb$ and using
  $(x+n)(x+n-1)\ldots x=\frac{\Gamma
  (x+n+1)}{\Gamma(x)}.$ When the Gamma functions blow
  up, the formulas must be understood as
  limits. Also

  \bean
  \left.J_{\lb}(x)
   \right|_{\sum_{\ell}x_{\ell}^i=\delta_{1i}u}
   &=& u^{|\lb|}
   \\
   \left.s_{\lb}(x)
   \right|_{\sum_{\ell}x_{\ell}^i=\delta_{1i}u}
   &=& \frac{u^{|\lb|}}{h^{\lb}}.
   \eean

\subsubsection*{Expansion of $(x_1+x_2+\ldots)^n$:}

\be
\frac{(x_1+x_2+\ldots)^n}{n!}=\al^n\sum_{\vert\lb\vert
=n} \frac{J_{\lb}^{(\al)}(x)}{j_{\lb}^{(\al)}}. \ee
 Then also
 \be (x_1+x_2+\ldots)^n=
 \sum_{\vert\lb\vert
=n}  C^{(\al)}_{\lb}(x)
 ~,~\mbox{with}~~
C^{(\al)}_{\lb}(x)=\frac{|\lb|!~
  \al^{|\lb|}}{j_{\lb}^{(\al)}} J_{\lb}^{(\al)}(x).
 \ee

\subsubsection*{Cauchy identity:}

\be
  \prod_{(i,j)\geq
1}(1-x_iy_j)^{-1/\al}=\sum_{\lb\in\BY}\frac{J_{\lb}^{(\al)}(x)
J_{\lb}^{(\al)}(y)}{j_{\lb}^{(\al)}};
 \ee
  in particular,
for $\al=1$
 \be \prod_{i,j\geq
1}(1-x_iy_j)^{-1}=\sum_{\lb\in\BY}s_{\lb}(x)s_{\lb}(y).
\ee

%\newpage

\subsubsection*{Hypergeometric functions}

Generalized hypergeometric functions ${}_2F_1^{(\al)}$
are defined by:
\be
{}_2F_1^{(\al)}(p,q;n;x): =\sum_{\kappa \in \BY}
 \frac{(p)^{(\al)}_{\kappa}(q)^{(\al)}_{\kappa}}
  {(n)^{(\al)}_{\kappa}}
~
\al^{|\kappa|}\frac{J_{\kappa}^{(\al)}(x)}
 {j_{\kappa}^{(\al)}}.
 \ee
For $\al=1$, using
$J_{\kappa}^{(1)}=h^{\kappa}s_{\kappa}$
 and $j_{\kappa}^{(1)}=(h^{\kappa})^2$, we have

\be
{}_2F_1^{(1)}(p,q;n;x): =\sum_{\kappa \in \BY}
 \frac{(p)^{(1)}_{\kappa}(q)^{(1)}_{\kappa}}
  {(n)^{(1)}_{\kappa}}
~  \frac{s_{\kappa}(x)}
 {h^{\kappa}},
 \ee
 and so, upon restriction,
\be
   {}_2F_1^{(1)}(p,q;n;x)
 \Bigr|_{\sum_{\ell}x_{\ell}^i= \delta_{1i}u}
=
 \sum_{\kappa \in \BY}u^{|\kappa|}
 \frac{(p)^{(1)}_{\kappa}(q)^{(1)}_{\kappa}}{({h^{\kappa}})^2
  (n)^{(1)}_{\kappa}}
 \ee
\subsubsection*{Generalized Selberg formula:}
 Kaneko~\cite{Kaneko} computes
 the following integrals, subjected to the condition that
 $a,b>\beta (p-1)$ (see also MacDonald~\cite{MacDonald}
  and
 Kadell~\cite{Kadell}):

\bea \lefteqn{
\int_{[0,1]^p}J_{\lb}^{(1/\beta)}(x)|\Delta_p(x)\Bigr|^{2\beta}
\prod^p_{i=1}(1-x_i)^{a-\beta (p-1)-1}x_i^{b-\beta
(p-1)-1}dx_i} \nonumber\\
&=&J_{\lb}^{(1/\beta)}(1^p)\prod^p_{i=1}
\frac{\Gamma(i{\beta}+1)~\Gamma(a+\beta(1-
i))~\Gamma(\lb_i +b+\beta(1-i))}{\Gamma(\beta+1)~
\Gamma(\lb_i+a+b+\beta(1-i))}.\nonumber\\
 \eea
  Setting
$u=(u_1,\ldots,u_m)$, we have the following
representation in terms of the hypergeometric
function:
   \bea \lefteqn{
\int_{[0,1]^p}\prod_{1\leq i\leq p \atop 1\leq k\leq
m}(1-x_iu_k)^{-\beta} |\Delta_p(x)\Bigr|^{2\beta}
\prod^p_{i=1}(1-x_i)^{a-\beta (p-1)-1}x_i^{b-\beta
(p-1)-1}dx_i} \nonumber\\
  &=&
 {}_2F_1^{(1/\beta)}(\beta p,b,a+b;u)
 \prod^p_{i=1}
\frac{\Gamma(i{\beta}+1)~\Gamma(a+\beta(1-
i))~\Gamma(b+\beta(1-i))}{\Gamma(\beta+1)~
\Gamma(a+b+\beta(1-i))}.\nonumber\\ \eea

%\newpage

\section{Probability measures on partitions}

\subsection{Probability measure on the set $\BY$ of
all partitions}

 In view of formula (2.3.8), define
the (not necessarily positive) probability measure on
the space $\BY$ of Young diagrams, depending on $x,y$
and $\al$: (see \cite{Borodin1,
  Borodin2,Borodin3,Vershik and Kerov,Kerov})
 \be
%p_{x,y}
P(\lb):=\frac{J_{\lb}^{(\al)}(x)J_{\lb}^{(\al)}(y)}{j_{\lb}^{(\al)}
\displaystyle{\prod_{i,j}}(1-x_iy_j)^{-1/\al}},~~~~\lb
\in \BY.
 \ee
   In particular, evaluating $P(\lb)$
along the locus
 \be
\LR=\left\{\begin{array}{l}\mbox{for}~x=
(x_1,x_2,\ldots)\mbox{~such that~}
\displaystyle{\sum_{\ell}}x_{\ell}^i=
\delta_{1i}u,\mbox{~we have~}J_{\lb}^{(\al)}
 (x)=u^{|\lb|}\\
 %\\
\mbox{for}~y=(\overbrace{1,1,\ldots,1}^p,0,0,\ldots)=1^p,
\mbox{~we have~}J_{\lb}^{(\al)}(1^p)\mbox{~as
in~}(2.3.4),
\end{array}\right.
 \ee
and using\footnote{$%H(z):=
\displaystyle{\prod_{i\geq
1}(1-x_iz)^{-1}=\exp\left(\sum^{\iy}_{\ell=1}\frac{z^{\ell}}{\ell}\sum_{i\geq
1}x_i^{\ell}\right)}$.}

$$ \prod_{i,j\geq 1}(1-x_iy_j)^{-1}=\prod_{j=
1}^p\prod_{i\geq 1}(1-x_iy_j)^{-1}
\Big\vert_{y_j=1}=\prod_{i\geq 1}(1-x_i)^{-p}=%(H(1))^p=
e^{pu}, $$ we obtain the genuine ($\geq 0$)
probability measure for $u>0$ on the space $\BY$,
depending on $u\in \BR$ and the integer $p>0$,
 \bea
  P_{u,p}(\lb):=P(\lb)\Bigl|_{\LR}
   &=&
 e^{-pu/\alpha}u^{|\lb|}
 \frac{
J_{\lb}^{(\alpha)}(1^p)}{j_{\lb}^{(\alpha)}}
 \nonumber \\
 &=&
  e^{-pu/\alpha}%u^{|\lb|}
  \frac{(pu/\alpha)^{|\lb|}}{|\lb|!}
  P^{\ell,p}(\lb),~\mbox{with}~~\ell:=|\lb|, %\in \BY
 \eea
 with
$$ (\mbox{support}~P_{u,p})\subseteq  \BY^{(p)}
:=\{\lb \in \BY \mbox{ such that}~ \lb_1^{\top} \leq
p\}.
$$
Probability (3.1.3) can be viewed as a Poissonized probability of
 \be
P^{\ell,p}(\lb)= \frac{
J_{\lb}^{(\alpha) }(1^p)\ell!}{
 j_{\lb}^{(\alpha)}(p/\alpha)^{\ell}
 },~~\mbox{for}~\lb\in \BY^{(p)}_{\ell}
 %=\{\lb~\mbox{with}~|\lb|=\ell\}
   .
 \ee
Probability $P_{u,p}(\lb)$ is $\geq 0$, because, from
(2.3.4),
 $J^{(\al)}_{\lb}(1^p) > 0$, for $p\geq \lb_1^{\top}$
  and $=0$ otherwise.

 In
particular, setting $\al=1$ and using
$J_{\lb}^{(1)}=h^{\lb}s_{\lb}$ and
$j_{\lb}^{(1)}=(h^{\lb})^2$, (3.1.3) leads to
  \be
P_{u,p}(\lb)\Bigr|_{\al=1}=\frac{u^{|\lb|}J_{\lb}^{(1)}(1^p)}
 {e^{pu}j^{(1)}_{\lb}}=
\frac{u^{|\lb|}s_{\lb}(1^p)}{e^{pu}h^{\lb}}
   = e^{-pu}\frac{(up)^{|\lb|}}{|\lb|!}
    P^{\ell,p}(\lb)\Bigr|_{\al=1}
%   \frac{f^{\lb}s_{\lb}(1^p)}{p^{|\lb|}}
%
 ~,~~~\lb
\in \BY^{(p)},
 \ee
 where by (2.1.2)
 \be
  P^{\ell,p}(\lb)\Bigr|_{\al=1}=
   \frac{f^{\lb}s_{\lb}(1^p)}{p^{|\lb|}}
   ,~~~\lb \in \BY^{(p)}_{\ell}.\ee
%Note that in both cases, the probability has support
%on the set $ \{\lb \in \BY~|~ \lb_1^{\top}\leq p\}.$
Probability (3.1.6) will be considered next.

%\newpage

\subsection{Probability measures on the set
$\BY_{\ell}$ of partitions of $\ell$ and random words}

From (3.1.6), setting in this section
$P^{\ell,p}:=
 P^{\ell,p}\Bigr|_{\al=1}$,
   \bean
  P^{\ell,p}(\lb)&=&
        \frac{f^{\lb}~s_{\lb}(1^p)}{p^{|\lb|}}
%\frac{1}{Z_{\ell,m}}
% f^{\lb}
%   \prod_{(i,j)\in\lb}\frac{j-i+m}{h^{\lambda}(i,j)}
    ,~~\mbox{for}~ \lb \in\BY^{(p)}_{\ell}
%    =\{\lb~\mbox{with}~|\lb|=\ell\}
    \nonumber\\
    &=&\frac{\ell !}{p^{\ell}}
     \Delta_{p}(h)^2 \frac{1}{\prod_1^{p-1}i!  \prod_1^p
 h_i!}
,~ \mbox{by (2.1.2) and (2.1.3),}
  \\
   &=& \mbox{probability on Young diagrams
    $\lb\in \BY^{(p)}_{\ell}$ coming from the
  } \\ & & \mbox{ uniform distribution on
  $S_{\ell}^p$, via the RSK correspondence,
  }\eean
 \vspace{-.7cm}
   \be  \ee
where $h:=(h_1,...,h_p)$ with $h_i:=p+\lb_i-i$, and
where $$S_{\ell}^p:= \left\{\mbox{words of length
$\ell$, built from an
 alphabet $\{1,...,p\}$}\right\},$$
with $|S_{\ell}^p|=p^{\ell}$.
% As before, the support
%of the  probability $P^{\ell,p}$ is given by $ \{\lb
%\in \BY_{\ell}~|~ \lb_1^{\top}\leq p\}.$
As already pointed out, this is a probability, firstly
because of the connection with the word problem,
secondly because of (3.1.5). As already pointed out,
this probability was considered in \cite{Borodin1,
  Borodin2,Borodin3,Vershik and Kerov,Kerov} and also
  in the context of random words, by Tracy and Widom
  (\cite{TW4}).

\begin{proposition} Given rectangular Young diagrams
 $\mu =(n-p)^p
\supseteq\mu' =(n-q)^p$, with $p\leq q <n$, and
$\ell\geq p(n-p)$, the expectation equals

\bigbreak

$E^{\ell,p}\displaystyle{\left(I_{\{\lb\supseteq\mu\}}
 (\lb)\frac{h^{\lb\backslash\mu'}}
  {h^{\lb\backslash\mu}}\right)}$
\be
=\frac{\ell
!}{p^{\ell}}\prod_1^p\frac{(q-i)!}{(n-i)!}\sum_{\kappa
\in \BY_{\ell -p(n-p)} \atop \kappa_1^{\top}\leq
p}\frac{1}{(h^{\kappa})^2}\prod^p_{i=1}
\frac{(p-i+1)_{\kappa_i}(q-i+1)_{\kappa_i}}{(n-i+1)_{\kappa_i}}.
\ee
\end{proposition}

\proof From (2.2.1) with $\lb \mapsto \lb\backslash
\mu'$, (so $ n\mapsto q,~p\mapsto p$)
 we have, upon setting $\kappa
=\lb\backslash\mu$,
\be
\prod_{(i,j)\in\mu\backslash\mu'}h^{\lb}_{(i,j)}
 =\frac{h^{\lb\backslash\mu'}}{
h^{\lb\backslash\mu}}=h^{\mu\backslash\mu'}\prod^p_{i=1}
 \frac{(q-i+1)
_{\kappa_i}}{(p-i+1)_{\kappa_i}}. \ee
 Then combining (3.2.3) and (2.2.3), $$
\frac{s_{\lb}(1^p)}{h^{\lb}}\frac{h^{\lb\backslash\mu'}}{h^{\lb\backslash\mu}}=
\frac{h^{\mu\backslash\mu'}}{h^{\mu}}\frac{1}{(h^{\lb\backslash\mu})^2}
\prod^p_{i=1}\frac{(q-i+1)_{\kappa_i}(p-i+1)_{\kappa_i}}{
(n-i+1)_{\kappa_i}}. $$
 In particular, setting
$\lb=\mu$, and using (2.2.2),
  $$
\frac{s_{\lb}(1^p)}{h^{\lb}}
 \frac{h^{\lb\backslash\mu'}}{h^{\lb\backslash\mu}}
 \Bigg|_{\lb=\mu}=
\frac{h^{\mu\backslash\mu'}}{h^{\mu}}=\prod_1^p\frac{(q-i)!}{(n-i)!}
,$$
 and therefore
 \bea
\frac{s_{\lb}(1^p)}{h^{\lb}}\frac{h^{\lb\backslash\mu'}}
{h^{\lb\backslash\mu}}
 &= &
 \left(\frac{s_{\lb}(1^p)}{h^{\lb}}
\frac{h^{\lb\backslash\mu'}}
 {h^{\lb\backslash\mu}}\right)
 \Bigg|_{\lb=\mu}\frac{1}{(h^{\lb\backslash\mu})^2}
  \prod_{i=1}^p
  \frac{(p-i+1)_{\kappa_i}(q-i+1)_{\kappa_i}}
   {(n-i+1)_{\kappa_i}}\nonumber\\ && \nonumber\\
&=&\prod_{i=1}^p\frac{(q-i)!}{(n-i)!}\frac{1}{(h^{\lb\backslash\mu})^2}
\prod_{i=1}^p\frac{(p-i+1)_{\kappa_i}(q-i+1)_{\kappa_i}}
{(n-i+1)_{\kappa_i}}. \eea
 Then, taking the
expectation $E^{\ell,p}$ with regard to the
probability measure $P^{\ell,p}$, defined in (3.2.1),
\bean
 \lefteqn{E^{\ell,p}\left(I_{\{\lb\supseteq\mu\}}(\lb)\prod_{n-q<j\leq
n-p} h^{\lb}_{(i,j)}\right)}\\
&=&\sum_{\lb \in
\BY_{\ell}  \atop
\lb\supseteq\mu}\frac{f^{\lb}s_{\lb}(1^p)}{p^{\ell}}\frac{h^{\lb\backslash\mu'}}
{ h^{\lb\backslash\mu}}\\ &=&\frac{\ell
!}{p^{\ell}}\sum_{\lb \in \BY_{\ell}
\atop\lb\supseteq\mu} \frac{s_{\lb}(1^p)}{h^{\lb}}
 \frac{h^{\lb\backslash\mu'}}{h^{\lb\backslash\mu}},
  \mbox{~~using~~}f^{\lb}=\frac{|\lb|!}{h^{\lb}}
   \mbox{~~as in (2.1.2),}\\ &=&\frac{\ell
!}{p^{\ell}}\prod_1^p\frac{(q-i)!}{(n-i)!}\sum_{\lb
\in \BY_{\ell}  \atop
\lb\supseteq\mu}\frac{1}{(h^{\lb\backslash\mu})^2}
 \prod^p_{i=1}\frac{(p-i+1)_{(\lb \backslash \mu)_i}
  (q-i+1)_{(\lb \backslash \mu)_i}}
   {(n-i+1)_{(\lb \backslash \mu)_i}},\\
   && \hspace{9cm}\mbox{using (3.2.4)}\\
\eean \bean &=&\frac{\ell
!}{p^{\ell}}\prod_1^p\frac{(q-i)!}{(n-i)!}\sum_{\kappa
\in \BY_{\ell-p(n-p)}  \atop \kappa_1^{\top}\leq
p}\frac{1}{(h^{\kappa})^2}\prod^p_{i=1}\frac{(p-i+1)_{\kappa_i}
(q-i+1)_{\kappa_i}}{(n-i+1)_{\kappa_i}}. \eean

%\newpage

\section{Expressing integrals as
mathematical expectation on partitions and on random
words}

\subsection{Expressing an integral on $Gr(p,\BF^n)$ as a
mathematical expectation on partitions}

Remembering the probability (3.1.3)
  \be P_{x,p}(\lb)=e^{-\beta
px}\frac{x^{|\lb|} J^{(1/\beta)}_{\lb}(1^p)}{
 j_{\lb}^{(1/\beta)}},~~~\lb \in \BY,\ee
 on partitions $ \BY$,
 with support on $\lb^{\top}_1 \leq p$,
 the following statement holds:

\begin{theorem}
For fixed $p\leq q\leq n/2$, the following holds
($\beta=1/2,1,2$)
  \bea
    \lefteqn{
    c^{-1}
 \int_{Gr(p,\BF^n)}e^{x\Tr(I+Z^{\dag}Z)^{-1}}\det(Z^{\dag
}Z)^{-\beta(q-p)} d\mu(Z) } %\nonumber
\\ &&\nonumber\\
 &=&
 ~
   {}_2F_1^{(1/\beta)}(\beta p,\beta q;\beta n;y)
    \Bigr|_{\sum_{\ell}y_{\ell}^i= \frac{x}{\beta}\delta_{1i}}
      %\nonumber
      \\  &&\nonumber\\
&=& e^{\beta
px}E_{x,p}\left(\frac{1}{\beta^{|\lb|}}\frac{J_{\lb}^{(1/\beta)}(1^q)}
{J_{\lb}^{(1/\beta)}(1^n)}%I_{\{\lb_1^{\top}\leq p\}}(\lb)
 \right)
 \eea
 where $c:=c^{(\beta)}_{n,q,p}$ is as in (0.0.20).
\end{theorem}

\proof
 For a symmetric function $f(z_1,\ldots,z_p)$, define
the integral, depending on $\beta$,
  \be  \Bigl\la f \Bigr\ra_{\beta}:=
 \int_{[0,1]^p}f(z_1,\ldots,z_p)~ |\Delta(z)|^{2\beta}
\prod_{i=1}^pz_i^{\beta(q-p+1) -1}
 (1-z_i)^{\beta(n-p-q+1)-1}
dz_i.
 \ee

 Kaneko's formula (2.3.13) will be used for
$a=\beta (n-q)$ and $b=\beta q$ in the sequence of
identities below; the inequalities $p\leq q\leq n/2$
imply $q-p+1,~n-p-q+1\geq1$, and so the integral
(4.1.5)
 above makes sense, and we first apply (1.0.3) in the
 following sequence of identities:
%guarantee $a-(\beta m+1-\beta),~b-(\beta
%, m+1-\beta)>-1$. So,

\bean
  \lefteqn{\int_{Gr(p,\BF^n)}e^{x\Tr(I+Z^{\dag}Z)^{-1}}\det(Z^{\dag
}Z)^{-\beta(q-p)} d\mu(Z)}\\
&=&
 \left\la \prod_{i=1}^pe^{xz_i}
\right\ra_{\beta} \\
&=&\sum_0^{\iy}x^{\ell}\left\la
\frac{\left(\displaystyle{\sum_1^p}z_i\right)^{\ell}}
{\ell !}\right\ra_{\beta}\\
&=&\sum_0^{\iy}x^{\ell}\left\la\frac{1}{\beta^{\ell}}
 \sum_{{{|\lb|=\ell}
 %\atop{\lb^{\top}_1\leq p}
 }}
\frac{J_{\lb}^{(1/\beta)}(z)}{j_{\lb}^{(1/\beta)}}
\right\ra _{\beta}, \mbox{~using (2.3.6)}\\
&\stackrel{\ast}{=}&\prod^p_{i=1}
\frac{\Gamma(i\beta+1)\Gamma(\beta(n-q-i+1))}
{\Gamma(\beta +1)} \\
 && \hspace{1cm} \sum^{\iy}_{\ell
=0}\left(\frac{x}{\beta}\right)^{\ell}
 \sum_{{{|\lb|=\ell}%\atop{%\lb^{\top}_1\leq p}
  }}
\frac{J_{\lb}^{(1/\beta)}(1^p)}{j_{\lb}^{(1/\beta)}}\prod^p_{i=1}\frac{\Gamma(\lb_i+\beta(q-
i+1))}{\Gamma (\lb_i+\beta(n-i+1))},\\
 && \hspace{7cm}
 \mbox{using Kaneko's formula (2.3.13),} \\
&=&\prod^p_{i=1}\frac{\Gamma(i\beta+1)\Gamma(\beta(n-q-i+1))}{\Gamma(\beta
  +1)}
 e^{\beta px}\\
 &&\hspace{1cm}\sum_{\lb \in \BY%\lb_1^{\top}\leq p
 }\frac{1}{\beta^{|\lb|}}P_{x,p}(\lb)\prod^p_{i=1}
\frac{\Gamma(\lb_i+\beta(q-i+1))}
 {\Gamma(\lb_i+\beta(n-i+1))},~~
\begin{array}{l}~\mbox{using the
definition }\\ \mbox{(4.1.1) of probability,}\\
 \mbox{with support on $\BY^{(p)}
 %\lb_1^{\top}\leq p
 $}
\end{array}
 \\ &=&c^{(\beta)}_{n,q,p}e^{\beta
px}\sum_{\lb \in \BY%\lb_1^{\top}\leq p
 }P_{x,p}(\lb)\frac{J_{\lb}^{(1/\beta)}(1^q)}
{\beta^{|\lb|}J_{\lb}^{(1/\beta)}(1^n)},~~
\begin{array}{l}~\mbox{
using
formula (2.3.4) %for $J_{\lb}^{(1/\beta)}(1^n)$
and }\\\mbox{~~the value (0.0.20) of $c^{(\beta)}_{n,q,p}$,}
\end{array}
\\ &=&c^{(\beta)}_{n,q,p}e^{\beta
px}E_{x,p}\left(\frac{1}{\beta^{|\lb|}}
\frac{J_{\lb}^{(1/\beta)}(1^q)}
{J_{\lb}^{(1/\beta)}(1^n)}
 %~I_{\{\lb_1^{\top}\leq p\}}(\lb)
  \right), \eean
  \vspace{-1.5cm}\be \ee
yielding (4.1.4).

Finally, looking at the expression to the right of
$\stackrel{\ast}{=}$, we find, using (3.1.2) and
(2.3.4), to be precise,

$$
 \left.J_{\lb}^{(1/\beta)}(y)
  \right|_{\sum_{\ell}y_{\ell}^i=
\frac{x}{\beta}\delta_{1i}}
 =\left( \frac{x}{\beta}  \right)^{|\lb|}
~~\mbox{and}~~
 J_{\lb}^{(1/\beta)}(1^p)=
  \left( \frac{1}{\beta}  \right)^{|\lb|}
  \left(  p\beta \right)_{\lb}^{(1/\beta)},$$
the following

 \bean  \left\la \prod_{i=1}^pe^{xz_i}
\right\ra_{\beta} &=&c^{(\beta)}_{n,q,p}\sum_{\lb
\in \BY }\frac{(\beta p)_{\lb}(\beta q)_{\lb}}{(\beta
n)_{\lb}} \frac{1}{\beta^{|\lb|}}
\left.\frac{J_{\lb}^{(1/\beta)}(y)}{j_{\lb}
^{(1/\beta)}}\right|_{\sum_{\ell}y_{\ell}^i=
\frac{x}{\beta}\delta_{1i}}
 ,~~~~\mbox{by (2.3.4),}\\ &=&
 ~
  c^{(\beta)}_{n,q,p}~ {}_2F_1^{(1/\beta)}(\beta p,\beta q;\beta n;y)
    \Bigr|_{\sum_{\ell}y_{\ell}^i= \frac{x}{\beta}\delta_{1i}}
     ,~~~~\mbox{by (2.3.10),} \nonumber
 \eean
 thus ending the proof of Theorem 4.1.\qed

  \remark Identity (4.1.3) is also an immediate
  consequence of Kaneko's formula.

%\newpage

\subsection{Expressing an integral over $Gr(p,\BC^n)$
 as a mathematical expectation on partitions and random
words}

We now specialize the previous section to $\BF=\BC$.
For fixed integer $p\geq 1$ and $\beta=1$, recall the
probability (3.1.5), with support in $\BY^{(p)}$,
 \be
P_{x,p}(\lb)=
\frac{x^{|\lb|}s_{\lb}(1^p)}{e^{px}h^{\lb}}
   = e^{-px}\frac{(xp)^{|\lb|}}{|\lb|!}
    P^{\ell,p}(\lb)
%   \frac{f^{\lb}s_{\lb}(1^p)}{p^{|\lb|}}
%
 ~,~~~\lb
\in \BY^{(p)},
 \ee
and the probability (3.1.6) on $\BY^{(p)}_{\ell}$  coming from the
  uniform distribution on
  $S_{\ell}^p$, via the RSK correspondence,
 \be
  P^{\ell,p}(\lb)=
   \frac{f^{\lb}s_{\lb}(1^p)}{p^{|\lb|}}
   ,~~~\lb \in \BY^{(p)}_{\ell}.\ee

   For integer $0\leq p\leq n/2$, consider
   the fixed rectangular Young diagram
   $\mu=(n-p)^p$.% of width $n-p$,
%$$\mu=(\overbrace{n-p,n-p,\ldots,n-p}^p)=(n-p)1^p.$$

\begin{theorem} For fixed $p\leq q\leq n/2$,
% ($U(p)\times U(n-p)$  and $U(q)\times U(n-q)$)
%\medbreak
\bea
 \lefteqn{\int_{Gr(p,\BC^n)}e^{x\Tr(I+Z^{\dag}Z)^{-1}}\det(Z^{\dag
}Z)^{-(q-p)} d\mu(Z)}\nonumber  \\
% &=& \int_{{Y\in\HR_p}\atop{0<Y<I}}e^{x\Tr
%  Y}\det Y^{q-p}\det(I-Y)^{n-q-p}dY \\
% &=&\tilde
% c_{n,q,p}\frac{e^{px}}{x^{(n-p)p}}E_{x,p}
% \left(I_{\{\lb\supseteq\mu \}}(\lb)\prod_{{(i,j)\in \lb}
% \atop {n-q<j\leq n-p}}
%  h_{(i,j)}^{\lb}\right),\\
  &=&
  \left\{ \begin{array}{l}
  \displaystyle{
 \frac{\tilde c_{n,q,p}}{x^{(n-p)p}}~e^{px}E_{x,p}
 \left(I_{\{\lb\supseteq\mu \}}(\lb)\prod_{{(i,j)\in \lb}
 \atop {n-q<j\leq n-p}}
  h_{(i,j)}^{\lb}\right)}\\
\\
 \displaystyle{
\frac{\tilde c_{n,q,p}}{x^{(n-p)p}}
 \sum_{\ell\geq p(n-p)}
  \frac{(px)^{\ell}}{\ell!}
 E^{\ell,p}\left( I_{\{\lambda \supseteq \mu \}}(\lb)
  \prod_{{(i,j)\in \lb}
 \atop {n-q<j\leq n-p}}
  h_{(i,j)}^{\lb}
   \right),}
   \\  \\
\displaystyle{   c_{n,q,p}^{(1)} \sum_{r\geq 0}x^r
\sum_{\kappa \in \BY_r \atop \kappa_1^{\top}\leq
p}\frac{1}{(h^{\kappa})^2}
 \frac{(p)_{\kappa}(q)_{\kappa}}{(n)_{\kappa}}
 }
\nonumber
\\ \displaystyle{
  c_{n,q,p}^{(1)}~{}_2F_1^{(1)}(p,q;n;y)
 \Bigr|_{\sum_{\ell}y_{\ell}^i= \delta_{1i}x}
}
\end{array} \right.
 \eea
  \vspace{-1.4cm}\be \ee
In particular, for $p=q$, the integral above has two
different formulations, as a probability or as a
generating function of probabilities,
  \bea
\lefteqn{\int_{Gr(p,\BC^n)}e^{x\Tr(I+Z^{\dag}Z)^{-1}}
 %\det(Z^{\dag }Z)^{\beta(p-q)}
   d\mu(Z)}\nonumber \\&&\nonumber \\
% &=&
%   \int_{{Y\in\HR_p}\atop{0<Y<I}}e^{x\Tr
%   Y}\det(I-Y)^{n-2p}dY  \\
%
%
 &=& \left\{ \begin{array}{l}
  \displaystyle{\frac{\tilde
 c_{n,p,p}}{x^{(n-p)p}}~e^{px}P_{x,p}
 \left(\lb\supseteq\mu  \right)}\\
\\
 \displaystyle{
\frac{\tilde c_{n,p,p}}{x^{(n-p)p}}
 \sum_{\ell\geq p(n-p)}
  \frac{(px)^{\ell}}{\ell!}
 P^{\ell,p}\left( \lb\supseteq\mu
   \right),}
   \\
 \displaystyle{
\frac{\tilde c_{n,p,p}}{x^{(n-p)p}}
 \sum_{\ell\geq p(n-p)}
  \frac{(px)^{\ell}}{\ell!}
 P^{\ell,p}\left(\pi \in S_{\ell}^p~\Bigr|~{
 {d_1(\pi)=p
 \mbox{ and} }\atop {  i_{p-1}(\pi)\leq
 \ell-n+p} } \right),}  \\
\displaystyle{   c_{n,p,p}^{(1)} \sum_{r\geq 0}x^r
\sum_{\kappa \in \BY_r \atop \kappa_1^{\top}\leq
p}\frac{1}{(h^{\kappa})^2}
 \frac{((p)_{\kappa})^2  }{(n)_{\kappa}}
%
%
%\prod_1^p \frac{(p-i+1)_{\kappa_i}(q-i+1)_{\kappa_i}}
% {(n-i+1)_{\kappa_i}}
 }
\nonumber
 \\ \displaystyle{
  c_{n,p,p}^{(1)}~{}_2F_1^{(1)}(p,p;n;y)
 \Bigr|_{\sum_{\ell}y_{\ell}^i= \delta_{1i}x}
}
     \end{array} \right.
    \eea
    \vspace{-1.2cm}\be \ee
  where $P^{\ell,p}$ denotes the uniform distribution on the
  set $S_{\ell}^p$ of words of length $\ell$ from an
  alphabet $1,\ldots,p$, and
  \bean
  d_1(\pi)&=&\mbox{length of longest strictly decreasing
 subsequence} \\
 i_k(\pi)&=&
  \mbox{length of the longest union of $k$ disjoint weakly increasing
 subsequences.}
 \eean
\end{theorem}

%\remark The constant in () and () equal
%$$\tilde C_{n,q,p}:=(\sqrt{\pi})^{p(p-1)}%\linebreak
%\displaystyle{\prod_1^p}(n-q-i)!~.$$

\remark
The constant $c^{(1)}_{n,q,p}$ in (4.2.2) and (4.2.3)
 is the same as
 $c^{(\beta)}_{n,q,p}$ for $\beta =1$ (see (0.0.20)) and $\tilde
 c_{n,q,p}$ is a new constant:
 \be
c^{(1)}_{n,q,p}=
  \prod_1^p\left({n-j}\atop
    {n-q-j,q-j,j}\right)^{-1}
      ~~~~~ \tilde c_{n,q,p}=\prod_{i=1}^p
i!~(n-q-i)!.
 \ee
 From the fourth expression of (4.2.2), it follows readily that, near $x=0$,
\bea \lefteqn{
 ( c^{(1)}_{n,q,p})^{-1}\int_{Gr(p,\BC^n)}e^{x\Tr(I+Z^{\dag}Z)^{-1}}\det(Z^{\dag
}Z)^{-(q-p)} d\mu(Z)} \nonumber\\
 &=&
    1+\frac{p~q}{n}x+
     \frac{pq}{4n}
     \left(\frac{(p+1)(q+1)}{n+1}+
    \frac{(p-1)(q-1)}{n-1} \right)x^2+\ldots
\eea

  \proof From the last identity in (2.1.2) with $m=p$, it
follows that, since $\lb_1^{\top}\leq p$,
  \be
\frac{h^{\lb+(q-p)^p}}{h^{\lb+(n-p)^p}}=\prod_{i=1}^p\frac{(q+\lb_i-i)!
}{(n+\lb_i-i)!} \ee
For a partition $\lb$ such that $\lb_1^{\top}\leq p$
and for $\lb'=\lb +k^p$, with arbitrary integer
$k\geq 0$, we have $s_{\lb}(1^p)= s_{\lb'}(1^p)$,
using the last identity (2.1.3).
 Using these facts, we have, continuing from (4.1.6),

 \bean
\lefteqn{\int_{Gr(p,\BC^n)}e^{x\Tr(I+Z^{\dag}Z)^{-1}}\det(Z^{\dag
}Z)^{p-q} d\mu(Z)}\nonumber
\\
&=&\int_{[0,1]^p}\Delta(z)^2\prod^p_{i=1}
 e^{xz_i}z_i^{q-p}(1-z_i)^{n-p-q}dz_i\\
&=&\tilde c_{n,q,p}\sum_{{\lb \in \BY
}\atop{\lb^{\top}_1\leq
p}}x^{|\lb|}\frac{s_{\lb}(1^p)}{h^{\lb}}\prod^p_{i=1}
\frac{(q+\lb_i-i)!}{(n+\lb_i-i)!}
,~~  \begin{array}{l}\mbox{\small by the 4th identity
of (4.1.6) and }\\ \mbox{\small
$J^{(1)}_{\lb}=h^{\lb}s_{\lb}$ and
 $j^{(1)}_{\lb}=(h^{\lb})^2$ }
\end{array}
\\
 &=&\tilde
c_{n,q,p}\sum_{\lb\in \BY
}x^{|\lb|}\frac{s_{\lb}^p}{h^{\lb}}
\frac{h^{\lb+(q-p)^p}}{h^{\lb+(n-p)^p}}
 ,~~
  \mbox{using (4.2.7)}
\eean \bean   %\\
   &\stackrel{\ast}{=}&\tilde
c_{n,q,p}\sum_{{\lb'\in \BY }\atop {\lb'\supseteq \mu}
}x^{|\lb'|-(n-p)p}\frac{s_{\lb'}(1^p)}{h^{\lb'}}
\frac{h^{\lb'-(n-q)^p}}{h^{\lb'-(n-p)^p}}
,~~
\begin{array}{l}\mbox{\small setting
 $\lb'=\lb+(n-p)^p$}\\
\mbox{\small  and using $s_{\lb}(1^p)= s_{\lb'}(1^p)$
}
\end{array}
 \\
&=&\frac{\tilde
c_{n,q,p}e^{px}}{x^{(n-p)p}}\sum_{{\lb' \in \BY
}\atop{\lb'\supseteq \mu} }
 P_{x,p}(\lb')
\frac{h^{\lb'-(n-q)^p}}{h^{\lb'-(n-p)^p}} ,~~ ,~~
\begin{array}{l}\mbox{\small using $P_{x,p}$ defined in (4.2.1), }\\
 \mbox{\small and $|\lb^{\prime}|=
  |\lb|+(n-p)p,$ }
\end{array}
 \\
&=&\frac{\tilde c_{n,q,p}e^{px}}{x^{(n-p)p}}E_{x,p}
 \left(I_{\lb\supseteq\mu}(\lb)\prod_{n-q
 <j\leq n-p}h^{\lb}_{(i,j)}\right).\\
 \eean
\vspace{-2cm}\be  \ee
  The last equality in (4.2.3)
follows from (4.1.3), while the second to the last one
follows from (2.3.12).

%\newpage

Finally, to prove the second formula on the right hand
side of (4.2.3), start with equality
$\stackrel{\ast}{=}$ in (4.2.8), omitting
 $\tilde c_{n,q,p}$, and replacing $\lb'$ by $\lb$,

\bean
 \lefteqn{ \frac{1}{x^{p(n-p)}}\sum_{\lb\in
\BY\atop
 \lb\supseteq \mu} x^{|\lb|}\frac{s_{\lb}(1^p)}{h^{\lb}}
\frac{h^{\lb-(n-q)^p}}{h^{\lb-(n-p)^p}},~~\mbox{for
fixed }~\mu=(n-p)^p }\\
 \\
&=&\frac{1}{x^{p(n-p)}}\sum_{\lb\in \BY\atop
 \lb\supseteq \mu}\frac{(px)^{|\lb|}}{|\lb|!}\frac{|\lb|!}{h^{\lb}}
\frac{s_{\lb}(1^p)}{p^{|\lb|}} \prod_{(i,j)\in\lb\atop
n-q<j\leq n-p}h^{\lb}_{(i,j)}\\
 \\
&=&\frac{1}{x^{p(n-p)}}\sum_{\ell\geq
p(n-p)}\frac{(px)^{\ell}}{\ell !} \sum_{\lb\in
\BY_{\ell}\atop\lb\supseteq
\mu}\frac{f^{\lb}s_{\lb}(1^p)}{p^{\ell}}
\prod_{(i,j)\in\lb\atop n-q<j\leq n-p}
h^{\lb}_{(i,j)},~~\mbox{using (2.1.2),}\\
 \\
 &=&
 \frac{1}{x^{p(n-p)}}\sum_{\ell\geq
p(n-p)}\frac{(px)^{\ell}}{\ell !}E^{\ell,p}
\left(I_{\{\lb\supseteq\mu\}}(\lb)
\prod_{(i,j)\in\lb\atop n-q<j\leq n-p}
h^{\lb}_{(i,j)}\right). \eean
In particular, setting $q=p$, the latter equals \bean
\lefteqn{\frac{1}{x^{p(n-p)}}\sum_{\ell\geq
p(n-p)}\frac{(px)^{\ell}}{\ell !}P^{\ell,p}(
\lb\supseteq \mu)}\\
&=&
 \frac{1}{x^{p(n-p)}}\sum_{\ell\geq
p(n-p)}\frac{(px)^{\ell}}{\ell !}
 P^{\ell,p}( \lb\supseteq \mu,\lb^{\top}_1=p)
  \\
&=&
 \frac{1}{x^{p(n-p)}}\sum_{\ell\geq p(n-p)}
  \frac{(px)^{\ell}}{\ell !}P^{\ell,p} (\pi \in
S_{\ell}^p~\Big\vert~d_1(\pi)=p, ~~i_{p-1}(\pi)\leq\ell
-n+p).
 \eean

To see the last two equalities, one proceeds as follows.
From $\lb \supseteq \mu=(n-p)^p$ and $P^{\ell,p}(\lb)=0$ for $\lb_1^{\top}>p$, it follows that
$\lb_1^{\top}=p$ and by Greene's theorem
 (see Sagan \cite{Sagan}, p. 110), $d_1(\pi)
=\lb_1^{\top}=p.$ Since also by Greene, $i_k(\pi)=\sum_1^k \lb_i$,
and, in particular $i_p(\pi)=\sum_1^{p} \lb_i=\ell$,
we have
 $$
 \lb_p=i_p(\pi)-i_{p-1}(\pi)=\ell-i_{p-1}(\pi)
 .
 $$
 From $\lb \supseteq \mu=(n-p)^p$, it also follows that
   $\lb_p\geq n-p$, and thus
 $i_{p-1}(\pi)\leq \ell-n+p$. Conversely, if
 $d_1(\pi)=\lb_1^{\top}=p$, and $i_{p-1}(\pi)\leq
 \ell-n+p$, then
  $\sum_1^{p-1}\lb_i=i_{p-1}(\pi)\leq
  \sum_1^{p}\lb_i-n+p$; hence $\lb_p\geq n-p$, and so
  $\lb \supseteq \mu$. \qed
%\newpage

\section{Testing Statistical Independence
 of Gaussian Populations}

 The statistical facts, used in this paper and summarized in this section, are due to
 James \cite{James} and Constantine
 \cite{Constantine}; see also Muirhead
 \cite{Muirhead}.

 \subsection{The Wishart distribution}

Let the $p\times n$ matrix $X$, with $n\geq p$ and $n$
identically distributed independent columns, have the
normal distribution
 \be (\det
2\pi\Sigma)^{-n/2}e^{-\frac{1}{2}\Tr
\Sigma^{-1}(X-M)(X-M)^{\top}} .
 \ee
 Then the $p\times p$
matrix $S=X\,X^{\top}$ has the {\em non-central Wishart
distibution with $n$ degrees of freedom}, with
$p\times p$
covariance matrix $\Sigma$ and non-centrality matrix
$\Omega=\frac{1}{2}M\,M^{\top}\Sigma^{-1}$, namely
 $$
 \Gamma_p(n/2)^{-1}(\det 2\Sigma)^{-n/2}e^{-\Tr(\Omega
+\frac{1}{2}\Sigma^{-1}S)}(\det
S)^{\frac{1}{2}(n-p-1)}\,_0F_{1}
\left(\frac{n}{2};\frac{1}{2}\Sigma^{-1}\Omega
S\right), $$
 where $S>0$, where $\Gamma_m$ is the multivariate Gamma function
 and where\footnote{Here $(a)_{\lb}:=
\prod_i(a+\beta(1-i))_{\lb_i}$ for $\beta=1/2$, with
$(x)_n:=x(x+1)\ldots (x+n-1),~x_0=1$.}(for the definition of $C_{\lb}$, see
(2.3.7))
 $$ %(\det S)^{\frac{n}{2}-(m+1)}
  \,_0F_1\left(\frac{n}{2};
\frac{1}{2}\Sigma^{-1}\Omega
S\right)=%\sum^{\iy}_{k=0}
 \sum_{\lb\in
\BY}\frac{C_{\lb} \left(\frac{1}{2}\Sigma^{-1}\Omega
S\right)} {\left(\frac{n}{2}\right)_{\lb}|\lb|!}
 =
  \sum_{\lb\in
\BY}\frac{2^{|\lb|}J^{(2)}_{\lb}
\left(\frac{1}{2}\Sigma^{-1}\Omega S\right)}
{\left(\frac{n}{2}\right)_{\lb}j_k^{(2)}}
 . $$
When $M=0$, we find the {\em (central) Wishart
distribution} $W_p(n, \Sigma)$, with $p\leq n$,
for the $p\times p$ matrix $S=X\,X^{\top}$:
\be
 \Gamma_p(n/2)^{-1}(\det 2\Sigma)^{-n/2}
 e^{-\frac{1}{2}\Tr\Sigma^{-1}S}
(\det S)^{\frac{1}{2}(n-p-1)}
 \prod_{{1\leq i\leq j \leq p}} dS_{ij}.
\ee

\subsection{The canonical correlation
 coefficients}

In testing the statistical independence of two
Gaussian populations, one needs to know the
distribution of {\em canonical correlation
 coefficients}.  To set up the problem, consider $p+q$ normally
distributed random variables $(X_1,...,X_p)^{\top}$
and $(Y_1,...,Y_q)^{\top}$
 ($p\leq q$) with mean zero and
covariance matrix

\smallbreak

$\hspace*{66mm}\stackrel{p}{\longleftrightarrow}
\quad\stackrel{q}{\longleftrightarrow}$
$$\mbox{cov}\MAT{1}X\\Y\mat :=\Sigma=
 \left(\begin{array}{cc}\Sigma_{11}&\Sigma_{12}\\
\Sigma_{12}^{\top}&\Sigma_{22}\end{array}\right)
 \begin{array}{l}\updownarrow p\\ \updownarrow q\end{array}
$$ The method proposed by Hotelling \cite{Hotelling}
is to find linear transformations $U=L_{1}X$ and
$V=L_2Y$ of $X$ and $Y$ having the property that the
correlation between the first components $U_1$ and
$V_1$ of the vectors $U$ and $V$ is maximal subject to
the condition that Var $U_1=$ Var $V_1=1$; moreover,
one requires the second components $U_2$ and $V_2$ to
have maximal correlation subjected to $$
\left\{\begin{array}{ll} {\rm (i)}\quad \mbox{Var\,}
U_2=\mbox{Var\,}V_2=1\\ {\rm (ii)}\quad U_2
\mbox{\,\,and\,\,} V_2 \mbox{\,\,are uncorrelated with
both\,\,} U_1 \mbox{\,\,and\,\,} V_1,
\end{array}\right.
$$ etc\ldots

Then there exist $O_p\in O(p)$, $O_q\in O(q)$ such
that $$
\Sigma_{11}^{-1/2}\Sigma_{12}\Sigma_{22}^{-1/2}=O^{\top}_{p}P
\,O_q $$ where $P$ has the following form:
\newpage $${{q}\atop{\longleftarrow\longrightarrow}}$$ $$
\left\{\begin{array}{lll} P=\left(
\begin{array}{lllllll|l}
\rho_1& & & & & & & \\
 &\ddots& & O& & & &\\
 & &\rho_k& & & & &O\\
& & &\rho_{k+1}& & & &  \\ & O& & & &\ddots& &\\ & & &
& & & \rho_{p}&
\end{array}
\right){\Bigg\updownarrow} p,\quad
k=\mbox{\,rank\,}\Sigma_{12},\\
\hspace*{3cm}{{\longleftarrow\longrightarrow}\atop{p}}
\\
 \\
1\geq\rho_1\geq\rho_2\geq\ldots\geq\rho_k
>0,~\rho_{k+1}= ...=\rho_p=0\quad\mbox{(canonical
correlation coefficients),}\\  \\ \rho_i \mbox{\,\,are
solutions ($\geq 0$) of\,}
\det(\Sigma_{11}^{-1}\Sigma_{12}\Sigma_{22}^{-1}\Sigma_{12}^{\top}-\rho^2I)=0.
\end{array}
\right. $$

\bigbreak

\noindent Then the covariance matrix of the vectors $$
U=L_1X:=O_p\Sigma_{11}^{-1/2}X\quad\mbox{and}\quad
V=L_2Y:=O_q\Sigma_{22}^{-1/2}Y $$ has the canonical
form ($\det \Sigma_{can}=\prod_1^p (1-\rho_i^2)$) $$
\mbox{cov}\MAT{1}U\\V\mat=\Sigma_{{can}}=
\MAT{2}I_{p}&P\\P^{\top}&I_q\mat ,
$$ with
 $$\mbox{spectrum}~\Sigma_{{can}}=
 \underbrace{1,\ldots ,1}_{q-p},1-\rho_1,1+\rho_1,\ldots,
  1-\rho_p,1+\rho_p
  $$
 and inverse $$
\Sigma_{can}^{-1}=\frac{1}{\prod_1^p(1-\rho_i^2)^2}
\MAT{2}I_{p}&-P\\-P^{\top}&I_q\mat.
 $$

\subsection{Distribution of the sample canonical correlations}

From here on, we may take $\Sigma=\Sigma_{can}$.
 The $n$ ($n\geq p+q$) independent samples
$(x_{11},\ldots,x_{1p},y_{11},\ldots,y_{1q})^{\top},\ldots$,
$(x_{n1},\ldots,x_{np},y_{n1},\ldots,y_{nq})^{\top}$,
arising from observing $\MAT{1}X\\Y\mat$ lead to a
matrix $\MAT{1}x\\y\mat$ of size $(p+q,n)$, having the
normal distribution \cite{Muirhead} (p. 79 and p. 539)

\bean &&(2\pi)^{-n(p+q)/2}(\det \Sigma)^{-n/2}
\exp ~{-\frac{1}{2}\Tr~ (x^{\top}~y^{\top})
\left(\begin{array}{cc}\Sigma_{11}&\Sigma_{12}\\
\Sigma_{12}^{\top}&\Sigma_{22}\end{array}\right)^{-1}
\MAT{1}x\\y\mat}\\ &&~~~= (2\pi)^{-n(p+q)/2}(\det
\Sigma)^{-n/2}
e^{ -\frac{1}{2}\Tr~\left(x^{\top}(\Sigma^{-1})_{11}x
+
 y^{\top}(\Sigma^{-1})_{22}y
 +2y^{\top}(\Sigma^{-1})_{12}^{\top}x\right)}
\eean The conditional distribution of $p\times n$
matrix $x$ given the $q\times n$ matrix $y$ is also
normal: \be (\det 2\pi\Omega)^{-n/2}
 e^{-\frac{1}{2}\Tr\Omega^{-1}
  (x-Py)(x-Py)^{\top}} \ee with
\bean \Omega&=&\Sigma_{11}-\Sigma_{12}\Sigma_{22}^{-1}
\Sigma_{21} =\diag(1-\rho_1^2,\ldots,1-\rho_p^2)\\
P&=& \Sigma_{12}\Sigma_{22}^{-1}. \eean
 Then the
 maximum likelihood estimates $r_i$ of the $\rho_i$
  satisfy the determinantal equation
\be
  \det
(S_{11}^{-1}S_{12}S_{22}^{-1}S_{12}^{\top}-r^2I)=
0,\ee
 corresponding to
  $$
  S=\MAT{2}S_{11}&S_{12}\\S_{12}^{\top}&S_{22}\mat
  :=
  \displaystyle{\left(\begin{tabular}{ll}$xx^{\top}$&$xy^{\top}$ \\
 $yx^{\top}$&$yy^{\top}$
   \end{tabular}\right)},
   $$
where $S_{ij}$ are the associated submatrices of the
{\em sample} covariance matrix $S$.

 \remark  The $r_i$
can also be viewed as $r_i=\cos \theta_i$
, where the $\theta_1,...,\theta_{p}$
 are the {\em critical} angles
 between two planes in $\BR^n$:

 (i) a
 $p$-dimensional plane $=$ span $\{(x_{11},...,x_{n1}),...,
 (x_{1p},...,x_{np})\}$

 (ii) a $q$-dimensional plane $=$ span $\{
 (y_{11},...,y_{n1})^{\top},...,~
 (y_{1{q}},...,y_{n{q}})\}$.

 \noindent As we shall see, $z_i=r_i^2=\cos^2 \theta_i$ are
 the precise variables $z_i$ appearing in section 1.

\bigbreak

Since the $(q,n)$-matrix $y$ has rank$(y)=q$, there
exists a matrix $H_n\in O(n)$ such that
$yH_n=(y_1\,\,\Big |\,\, O)$; therefore acting on $x$
with $H_n$ leads to

$\hspace*{35mm}\stackrel{q}{\leftrightarrow}\,\,\,\,\,\,
\stackrel{n-q}{\leftrightarrow}\hspace{3cm}
\stackrel{q}{\leftrightarrow}\,\,\,\,\,\,\stackrel{n-q}{\leftrightarrow}$

\vspace{-.8cm}

\be yH= (y_1\,\,\,\big |\,\,\,O)\updownarrow q, \qquad
xH_n=(u\,\,\,\big | \,\,\,v)\updownarrow p~~.
 \ee
 With this in mind,
 \bean
  \lefteqn{S_{12}
S_{22}^{-1}S_{12}^{\top}-r^2S_{11}}\\
 &=& xy^{\top}(yy^{\top})^{-1}
 yx^{\top}-r^2xx^{\top}  \\
&=&xH(yH)^{\top}(yH(yH)^{\top})^{-1}
 yH(xH)^{\top}-r^2(xH)(xH)^{\top}\\
 &=&(u~\big |
~v)\MAT{1}y^{\top}_1\\O\mat
 \left((y_1~\big |~ O)
  \MAT{1}y^{\top}_1\\O\mat\right)^{-1}
(y_1~\big | ~O) \MAT{1}u^{\top}\\v^{\top}\mat
-r^2(u~\big | ~v)\MAT{1}u^{\top}\\v^{\top}\mat\\
%
%
%&=&\det\left(w\MAT{2}I_q&O\\O&0_{n-q}\mat
%w^{\top}-r^2ww^{\top}\right)\\
&=&(u~\big |~v)\MAT{2}I_q&O\\O&0_{n-q}\mat
%\left(\frac{u^{\top}}{v^{\top}}\right)
 \MAT{1}u^{\top}\\v^{\top}\mat
-r^2(u~\big |~v)\MAT{1}u^{\top}\\v^{\top}\mat
 \\
 &=&uu^{\top}-r^2(uu^{\top}+vv^{\top}),
\eean
 and so the equation (5.3.2) for the $r_i$ can be
rewritten
 \be\det ( uu^{\top}-r^2(uu^{\top}+vv^{\top})=0.
 \ee
 Then setting the forms (5.3.3) of $x$ and $y$
 in the conditional distribution
(5.3.1) of $x$ given $y$, one computes the following,
setting $H:=H_n$,
%$$
%(\det
%2\pi\Omega)^{-n/2}e^{-\frac{1}{2}\Tr\Omega^{-1}(x-\Sigma_{12}\Sigma_{22}^{-1}y)
%(x-\Sigma_{12}\Sigma_{22}^{-1}y)^{\top}}
%$$

%\hfill with $\Omega=\Sigma_{11}-\Sigma_{12}\Sigma_{22}^{-1}\Sigma_{21}$

\medbreak

$\Tr\Omega^{-1}(x-Py)(x-Py)^ {\top}$ \bean
&=&\Tr\Omega^{-1}(xH-PyH)(xH-PyH )^{\top}\\
&=&\Tr\Omega^{-1}\left((u\,\,\big
|\,\,v)-P(y_1\,\,\big |\,\,O)\right)\left((u\,\,\big
|\,\,v)-P(y_1\,\,\big |\,\,O)\right)^{\top}\\
&=&\Tr\Omega^{-1}(u-
%\Sigma_{12}\Sigma_{22}^{-1}
 Py_1)(u-
 %\Sigma_{12}\Sigma_{22}
 Py_1)
^{\top} +\Tr\Omega^{-1}vv^{\top};  ~~\Omega
=\diag(1-\rho_1^2,\ldots,1-\rho_p^2); \eean this
establishes the independence of
 the normal distributions $u$ and $v$, given the
 matrix $y$,
 with
$$ u\equiv N(P
%\Sigma_{12}\Sigma_{22}^{-1}
y_1,\Omega),\quad v\equiv
N(O,\Omega).~~P=\diag(\rho_1,\ldots,\rho_p). $$ Hence
$uu^{\top}$ and $vv^{\top}$ are conditionally
independent and both Wishart distributed; to be
precise:

\smallbreak
\begin{itemize}
  \item The $p\times p$ matrices $vv^{\top}$ are Wishart
   distributed, given $y$, with $n-q$ degrees of freedom and
   covariance $\Omega$;

  \item The $p\times p$ matrices $uu^{\top}$ are
  non-centrally Wishart distributed, given $y$, with $q$ degrees
  of freedom, with covariance $\Omega$ and with
  non-centrality matrix
$$\frac{1}{2}Py_1y_1^{\top} P^{\top}\Omega^{-1}
%\frac{1}{2}\diag(\frac{\rho_i}{1-\rho_i^2})y_1y_1^{\top}
%\diag(\frac{\rho_i}{1-\rho_i^2})
 .$$

  \item The marginal distribution of the $q\times q$ matrices $yy^{\top}$
   are Wishart distributed,
 with $n$ degrees of freedom and covariance $I_q$,
 because the marginal distribution of $y$ is normal
 with covariance $I_q$.

\end{itemize}

 %Hotelling \cite{Hotelling} shows that after non-singular linear
%transformations on the $X$'s and
 %$Y$'s, the covariance matrix above has the canonical form

% In testing the null-hypothesis that the two sets of
%  variables $X_i$ and $Y_i$ are independent,
 % an interval of confidence will be found, based on
%  the distribution above, given that
 % the $\rho_i=0$.
%If the null-hypothesis is that $\rho_i$ is non-zero,
%then one has an additional factor in (2.1.1),
%involving sums of {\em zonal} polynomials in the
%$r_i^2=z_i$ (Jack polynomials).

To summarize, given the matrix $y$, the sample
canonical correlation coefficients $r_1^2>\ldots
>r_p^2$ are the roots of
\bean (r_1^2>\ldots >r_p^2)&=&\mbox{~roots
of~}\det(xy^{\top}(yy^{\top})^{-1}yx^{\top}
-r^2xx^{\top})=0\\ &=&\mbox{~roots
of~}\det(uu^{\top}-r^2(uu^{\top}+vv^{\top}))=0\\
&=&\mbox{~roots
of~}\det(uu^{\top}(uu^{\top}+vv^{\top})^{-1}-r^2I)=0.
\eean
%From (), it follows that

\medbreak

Then one shows that, knowing $uu^{\top}$ and $
vv^{\top}$ are Wishart and conditiionally independent,
the conditional distribution of $r_1^2>\ldots
>r^2_p$, given the
matrix $y$ is given by \bean
\pi^{p^2/2}c_{n,p,q}e^{-\frac{1}{2}\Tr
 P yy^{\top}
  P^{\top}\Omega^{-1}}
\Delta_p(r^2) \prod^p_1(r_i^2)^{\frac{1}{2}(q-p-1)}
(1-r_i^2)^{\frac{1}{2}(n-q-p-1)}.\\
%
%\sum^{\iy}_{k=0}
 \sum_{\lb \in \BY}
\frac{(n/2)_{\lb}C_{\lb}(\frac{1}{2}
 P yy^{\top}
  P^{\top} \Omega^{-1})}{(q/2)_{\lb}C_{\lb}(I_p)~|\lb|!}C_{\lb}(R^2),
\eean where\footnote{$c_{n,p,q}$ is a different
constant from (0.0.20).} $$
R^2=\diag(r_1^2,\ldots,r_p^2),\quad c_{n,p,q}=
\frac{\Gamma_p(n/2)}{\Gamma_p(q/2)
\Gamma_p((n-q)/2)\Gamma_p(p/2)}. $$ By taking the
expectation with regard to $y$ or, what is the same,
by integrating over the matrix $yy^{\top}$, which is
Wishart distributed, we obtain:

\begin{theorem} Let $X_1,\ldots,X_p,Y_1,\ldots,Y_q$ ($p\leq q$) be normally
distributed random variables with zero means and
covariance matrix $\Sigma
=\displaystyle{\MAT{2}\Sigma_{11}&\Sigma_{12}\\
\Sigma_{21}&\Sigma_{22}\mat}$. If
$\rho^2_1,\ldots,\rho^2_p$ are the roots of
$\det(\Sigma_{11}^{-1}\Sigma_{12}\Sigma_{22}^{-1}\Sigma_{12}^{\top}-\rho^2I)=0$,
then the maximum likelihood estimates
$r_1^2,\ldots,r^2_p$ from a sample of size $n$ ($n\geq
p+q$) are given by the roots of $$
\det(xy^{\top}(yy^{\top})^{-1}yx^{\top}-r^2xx^{\top})=
0 .$$
Setting
$$Z:=\diag(z_1,\ldots,z_p)=\diag(r^2_1,\ldots,r^2_p)
~~\mbox{and}~~P^2:
=\diag(\rho_1^2,\ldots,\rho_p^2),$$
the $z_i=r_i^2$ have the following density
\bea \lefteqn{
\pi^{-p^2}c_{n,p,q}\Delta_p(z) \prod_1^p
% \prod_{1\leq i<j\leq p}(r_i^2-r_j^2)
   z_i^{(q-p-1)/2}
  (1-z_i)^{%\frac{1}{2}
(n-q-p-1)/2}dz_i
 }
 \nonumber\\
&&. %\sum^{\iy}_{n=0}
 \Bigl(\prod_1^p
(1-\rho_i^2)^{n/2}\Bigr)\sum_{\lb \in \BY}
 \frac{
 \left(\frac{n}{2}\right)_{\lb}
\left(\frac{n}{2}\right)_{\lb} }
{\left(\frac{q}{2}\right)_{\lb}}
 \frac{C_{\lb}(Z)C_{\lb}(P^2)}{C_{\lb}(1^p)~|\lb|!}.
\eea
\end{theorem}

\begin{corollary} If $\rho_1^2=\ldots =\rho_p^2=0$, then the joint density
of the $z_i=r^2_i$ is given by the density appearing
in the integral of Theorem 1.1, namely
 \be \pi^{p^2/2}c_{n,p,q}
% \prod_{1\leq i<j\leq p}(r_i^2-r_j^2)
  \Delta_p(z)
  \prod^p_{i=1}z_i^{(q-p-1)/2}
   (1-z_i)^{(n-q-p-1)/2}
dz_i.
 \ee
\end{corollary}

\remark As was shown here, the normal distribution
over $\BR$ leads to the density (5.3.6) for the $z_i$,
which corresponds to the case $\beta=1/2$ for (1.0.3).
Starting with normal distributions over $\BC$ and
$\BH$ leads, in a similar way, to integrals (1.0.3)
for the cases $\beta=1$ and $2$.

%\newpage

\section{Differential equations for the Grassmannian integrals
 and the hypergeometric functions}

\subsection{Differential equations for the Grassmannian
integral }

Theorem 6.1 shows that the integral over the
Grassmannian $Gr(p,\BF^n)$ satisfies Painlev\'e-like
 differential equations; for $\beta=1$, this equation
 is the Painlev\'e V equation with a specific boundary
 condition (Theorem 6.2).

\begin{theorem}
 The following holds for the integral
  \be
 I_p(y)%=I_{n,q,p}^{(\beta)}
  = \int_{Gr(p,\BF^n)}e^{y\Tr(I+Z^{\dag}Z)^{-1}}\det(Z^{\dag
}Z)^{-\beta(q-p)} d\mu(Z)=
  c^{\prime}\exp \int_0^x H(y)dy
   ,
 \ee
  where  $H_{}(y)=\frac{d}{d y}\log I_p(y)$
   satisfies the differential equation ($
 H':=\frac{dH}{dy}$ and remember
 $\delta_1^{\beta}=1$ for $\beta=1$ and $=0$ otherwise):
\bean && \hspace{-1cm}
 4\left(y^3
H^{\prime\prime\prime}+6y^3{H^{\prime}
}^2+(1+\delta^{\beta}_{1})(2 y^2
 H^{\prime\prime}
 +4 y^2HH^{\prime}
  + yH^2)\right)-yP_0H^{\prime}+P_1H
  +P_2 \\
%  &&
% -yP_0H^{\prime}+P_1H
%  +P_2\\
&&\\ &&%-yP_0H^{\prime}+P_1H
 % +P_2
   =
 \left\{ \begin{array}{ll}
    0,~~& \mbox{for}~~ \beta =1, ~~~~~(\mbox{\bf Painlev\'e V}) \\ \\
   \displaystyle{  \frac{3}{16}\frac{p(p-1)}{(p+1)(p+2)}y^3~\frac{I_{p-2} I_{p+2}}{I_{p}^2}
} ,~~& \mbox{for}~~ \beta =1/2, \\
        \\
  \displaystyle{  \frac{3}{16^2}\frac{p}{p+1}y^3~\frac{I_{p-1}I_{p+1}}{I_{p}^2}
},~~& \mbox{for}~~ \beta =2,
 \end{array}
  \right.
  \eean
  \vspace{-.9cm} \be \ee
where $$
  \begin{array}{c|ccc}
   & P_0&P_1 &P_2\\
   \hline \\
  \beta=1& 4y^2-8sy+4n^2-8& 4(sy-n^2)    &
   4r(y+n)    \\
  \beta=1/2 &  4y^2-4sy+(n+2)^2-8   &
  2sy-n(n-2)   & r(2y+n-2) \\
  \beta=2 & y^2-4sy+4((n-1)^2-2)   &
  2sy-4n(n+1)   & r(2y+4(n+1)) \\
  \end{array}
  $$
 in terms of
  \be
  r=pq,~~~~s=n-2p-2q.
   \ee

\end{theorem}

%\newpage

\begin{theorem}
For $\beta =1$, we have
  \bea
 \lefteqn{
%\prod_1^p\left({n-j}\atop
%   {n-q-j,q-j,j}\right)^{}
 (c^{(1)}_{n,q,p})^{-1} \int_{Gr(p,\BC^n)}e^{x\Tr(I+Z^{\dag}Z)^{-1}}\det(Z^{\dag
}Z)^{-(q-p)} d\mu(Z)}\nonumber  \\
 &=&
  \displaystyle{ \prod_1^p\frac{(n-j)!}{(q-j)!}%\tilde c_{n,q,p}
 \sum_{\ell\geq p(n-p)}
  \frac{p^{\ell} x^{\ell -p(n-p)}}{\ell!}
 E^{\ell,p}\left( I_{\{\lambda \supseteq \mu \}}(\lb)
  \prod_{{(i,j)\in \lb}
 \atop {n-q<j\leq n-p}}
  h_{(i,j)}^{\lb}
   \right),}
  \nonumber \\
&=&    \displaystyle{ \sum_{r\geq 0}x^r \sum_{\kappa
\in \BY_r %\atop \kappa_1^{\top}\leq
%p
}
\frac{1}{(h^{\kappa})^2}\frac{(p)_{\kappa}(q)_{\kappa}}{(n)_{\kappa}}
%
%\prod_1^p \frac{(p-i+1)_{\kappa_i}(q-i+1)_{\kappa_i}}
% {(n-i+1)_{\kappa_i}}
    }
 \nonumber \\
 &=&
  ~{}_2F_1^{(1)}(p,q;n;y)
 \Bigr|_{\sum_{\ell}y_{\ell}^i= \delta_{1i}x}
  \nonumber \\
 &=&
    \exp {\displaystyle{\int_0^x\frac{u(y)-p(n-p)+py}{y}
    dy}
    }
 \eea
 where $u(x)$
  is the unique solution to the initial value problem:
 \be
 \left\{\begin{array}{l} \displaystyle{ x^2
u^{\prime\prime\prime}+xu^{\prime\prime}
 +6x{u^{\prime} }^2-4uu^{\prime}+4Qu^{\prime}
 -2Q^{\prime}u+2R=0 } \\\hspace{10cm}
  \mbox{\bf (Painlev\'e V)} \\
\displaystyle{\mbox{with}~~  u_{}(x)=p(n-p)-
 \frac{p(n-q)}{n}x+\ldots+
+  a_{n+1}x^{n+1}+O(x^{n+1})+\ldots ,~\mbox{near}~ x=0
.}
\end{array}
 \right. \ee
  with $a_{n+1}$ specified (see remark).
$Q$ and $R$ are polynomials in $x$:
  \bea
 4Q&=&
 -x^2+2(n+2(p-q))x-(n-2p)^2 \nonumber\\
 2R&=&
  p(p-q)(x+n-2p).
  \eea

 The third order
equation (6.1.5)
 has a first integral, which is second order in $u$
and quadratic in $u^{\prime\prime}$,
 \be
  u^{
\prime\prime 2} +\frac{4}{x^2}
 \left( (xu^{\prime 2}+Q u^{\prime}+R)u^{\prime}
 - ( u^{\prime 2}+\frac{}{}Q' u^{\prime}+R')u^{}
 +\frac{1}{2}Q^{\prime\prime} u^2
   -  \frac{p^2(q-p)^2}{4} \right)=0
\ee
\end{theorem}

\remark Note that the Painlev\'e equation (6.1.5)
admits a solution
 \be
  u(x)=p(n-p)- \frac{p(n-q)}{n}x +
  \sum_{i\geq 2} a_i  x^i
 \ee
 with the $a_i$ given by the indicial equation
$$ i(i-1-n)(i-1+n)a_i = g_i(a_0,\ldots , a_{i-1}),~~
 i=1,2,\ldots ~,$$
showing the existence of a free parameter at $i=n+1$.
However the fact that, according to (4.2.3),
 $$ u(x)= p(n-p)-px
 + x \frac{d}{dx} \log
 \sum_{r\geq 0}x^r \sum_{\kappa
\in \BY_r \atop \kappa_1^{\top}\leq
p}\frac{1}{(h^{\kappa})^2} \prod_1^p
\frac{(p-i+1)_{\kappa_i}(q-i+1)_{\kappa_i}}
 {(n-i+1)_{\kappa_i}}
 $$
 leads to an explicitly known value for $a_{n+1}$.

% \newpage

{\medskip\noindent{\it Proof of Theorem 6.1:\/} }
Define
 \be
  \tilde I_p(t):= \int_{[-1,1]^p}
 |\Delta_p(y)|^{2\beta}\prod_1^p e^{\sum_{i=1}^{\iy}
 t_iy_k^i}
  (1-y_k)^{a}
 (1+y_k)^{b}dy_k,
\ee
 and the locus \be \LR:=\{t_1=x
 \neq 0,
~\mbox{all other}~t_i=0 \}.\ee
 Using (1.0.3), and setting
  $a=\beta(n-q-p+1)-1, b=\beta(q-p+1)-1$ in (6.1.9), the linear
   change of variables $y_k:= 2z_k-1$ leads to:

\bea \tilde I_p(t)\Bigr|_{\LR}
 &=&
 \int_{[-1,1]^p}
 |\Delta_p(y)|^{2\beta}\prod_1^p e^{
 xy_k}
  (1-y_k)^{\beta(n-q-p+1)-1}
 (1+y_k)^{\beta(q-p+1)-1}dy_k  \nonumber \\
 &=&  c_p^{(\beta)} e^{-px}
 \int_{[0,1]^p}
  |\Delta_p(z)|^{2\beta}
  \prod_1^p e^{
 2xz_k}
  z_k^{\beta(q-p+1)-1}
 (1-z_k)^{\beta(n-q-p+1)-1}dz_k
  \nonumber\\
 &=& c_p^{(\beta)}
e^{-px}\int_{Gr(p,\BF^n)}e^{2x\Tr(I+Z^{\dag}Z)^{-1}}\det(Z^{\dag
}Z)^{-\beta (q-p)} d\mu(Z)  \nonumber\\
  &=&  c_p^{(\beta)}e^{-px}I_p(x), \eea with
\be
 c_p^{(\beta)}=2^{p(\beta n -\beta p+\beta-1)}.
 \ee

 %   \vspace{5cm}

According to the appendix, the integral $\tilde
I_p(t)$ satisfies the Virasoro constraints (8.0.10),
with
 $a=\beta(n-q-p+1)-1, b=\beta(q-p+1)-1$, and thus
  $$ b_0=a-b=\beta(n-2q),~b_1=a+b=\beta(n-2p+2)-2,$$
and $$ \sigma_1=\beta n ,~~\sigma_2=\beta (n-1)+1.$$

These expressions and their first $t_1$- and $t_2$-
derivatives, evaluated along the locus $ \LR$ read as
follows:
 ($F(t):=F_p(t):=\log \tilde I_p(t)$)
{\footnotesize \bean
 0&=&\left.\frac{{\cal
J}^{(2)}_{-1}\tilde I_{p}}{ \tilde
I_{p}}\right|_{\LR}=\left.\left(t_{1}\frac{\pl}{\pl
t_{2}}+ \sigma_1\frac{\pl}{\pl
t_{1}}\right)F_p+p(b_0-t_{1})\right|_{\LR}\\
0&=&\left.\frac{{\cal J}^{(2)}_{0}
 \tilde  I_{p}}{\tilde I_{p}}\right|_{\LR}=
\left(t_{1}\frac{\pl}{\pl t_{3}}+\sigma_2
\frac{\pl}{\pl t_2}+(b_0-t_1)
 \frac{\pl}{\pl
t_{1}}+\beta\frac{\pl^2}{\pl t^2_1}\right)F_p\\ &
&\hspace{5cm}+\left.\beta\left(\frac{\pl F_p}{\pl
t_1}\right)^2-\frac{p}{2}(\sigma_1-b_1)\right|_{\LR}\\
& &
\\
0&=&
 \left.\frac{\pl}{\pl t_1} \frac{{\cal
J}^{(2)}_{-1}\tilde I_{p}}{\tilde I_{p}}\right|_{\LR}
% =\left(\sum_{i\geq
% 1}it_{i}\frac{\pl^{2}}{\pl t_{i+1}\pl t_{1}}
% +\frac{\pl}{\pl t_{2}}-
%\sum_{i\geq 2}it_{i} \frac{\pl^{2}}{\pl t_{i-1}\pl
%t_{1}}\right.\\ &
%&\left.\left.\hspace{4cm}+(2n+a)\frac{\pl^2}{\pl
%t_1^2}\right)F_n\right|_{\LR}-n
%\\ && \hspace{2.3cm}
=\left.\left(t_1\frac{\pl^2}{\pl t_2\pl
t_1}+\frac{\pl}{\pl t_{2}}+\sigma_1\frac{\pl^2}{\pl
t_1^2}\right)F_p
 \right|_{\LR}-p\\
 & & \\
\eean \bean
 0&=&\left.\frac{\pl}{\pl t_1}\frac{{\cal
J}^{(2)}_{0}\tilde I_{p}}{\tilde I_{p}}\right|_{\LR}
 =\left(t_1\frac{\pl^2}{\pl t_3\pl
t_1}+\sigma_2\frac{\pl^2}{\pl t_2\pl
t_1}+(b_0-t_1)\frac{\pl^2}{\pl t_1^2}+\frac{\pl}{\pl
t_3} -\frac{\pl}{\pl t_1}\right.\\ &
&\hspace{3cm}\left.+\beta \frac{\pl^3}{\pl
t_1^3}\right)F_p+2\beta \left.\frac{\pl F_p}{\pl
t_1}\frac{\pl^2 F_p}{\pl t_1^2}\right|_{\LR}\\ & & \\
 0&=&\left.\frac{\pl}{\pl t_2}\frac{{\cal
J}^{(2)}_{-1} \tilde I_{p}}{\tilde I_{p}}\right|_{\LR}
 =\left.\left(t_1\frac{\pl^2}{\pl
t_2^2}+ \sigma_1\frac{\pl^2}{\pl t_1\pl
t_2}+2(\frac{\pl}{\pl t_3}-\frac{\pl}{\pl
t_1})\right)F_p
 \right|_{\LR}. \eean
 }
The five equations above form a (triangular) linear
system in five unknowns $$ \left.\frac{\pl F_p}{\pl
t_2}\right|_{\LR},
 \quad\left.\frac{\pl F_p}{\pl t_3}\right|_{\LR},\quad\left.
 \frac{\pl^2F_p}{\pl
t_1\pl
t_2}\right|_{\LR},\quad\left.\frac{\pl^2F_p}{\pl
t_1\pl t_3}\right|_{\LR},\quad
\left.\frac{\pl^2F_p}{\pl t_2^2}\right|_{\LR},\quad $$
which can be expressed in terms of
 $$
 \frac{\pl F_p}{\pl t_1^{}}~,~~
   \frac{\pl^{2} F_p}{\pl t_1^{2}}~,~~
    \frac{\pl^{3} F_p}{\pl t_1^{3}}~.
 $$
Setting $t_1=x$ and $F'_n=\pl F_n/ \pl x$,
 these expressions
are

{\footnotesize
 \bean
  \left.\frac{\pl F_n}{\pl
t_2}\right|_{\LR}
 &=&-\frac{1}{t_1}\Bigl(\sigma_1 F^{\prime}_n-p(t_1-b_0)\Bigr)\\
& & \\ \left.\frac{\pl F_n}{\pl
t_3}\right|_{\LR}&=&\frac{1}{t_1^2}\Bigl(-2\beta
t_1(F^{\prime\prime}_n+ F_n^{\prime
2})+2(t_1^2-b_0t_1+\sigma_1\sigma_2)F^{\prime}_n\\
&&\hspace{3cm}+p((\sg_1-2\sg_2-b_1)t_1+2b_0\sg _2)
\Bigr)\\ & &
\\ \left.\frac{\pl^2 F_n}{\pl t_1\pl
t_2}\right|_{\LR}&=&\frac{1}{t_1^2}
 \left(-\sigma_1t_1F_n^{\prime\prime}+\sigma_1
 F_n^{\prime}+b_0 p\right)\\
\left.\frac{\pl^2 F_n}{\pl t_1\pl t_3}\right|_{\LR}
 &=&\frac{1}{2t_1^3}\Bigl(-2\beta t_1^2(
 F_n^{\prime\prime\prime}+2
F_n^{\prime}F_n^{\prime\prime}-F_n^{\prime 2
})+2t_1(t_1^2-b_0t_1+\sg_1\sg_2+\beta)
F_n^{\prime\prime}\\
 & &\hspace*{3cm}+2(b_0t_1-2\sg_1\sg_2)F_n^{\prime
}+p((2\sg_2-\sg_1+b_1)t_1-4\sg_2b_0)\Bigr)\\ & & \\
\left.\frac{\pl^2 F_n}{\pl t_2^2}\right|_{\LR}
&=&\frac{1}{t_1^3}\Bigl( (\sg_1^2+2\beta)t_1F^{\prime
\prime}_n+2\beta t_1 F^{\prime
2}_n+(2b_0t_1-2\sg_1\sg_2-\sg_1^2)F^{\prime}_n\\ &
&\hspace*{3cm}+p((-\sg_1+2\sg_2
+b_1)t_1-b_0(\sg_1+2\sg_2))\Bigr). \eean }
\vspace{-.8cm}\be \ee
%
%
%\newpage
%
From (8.0.2), it follows that

 \bea
\tilde I_p(t)&=&
   \left\{
   \begin{array}{lll}
   p!\tau_p(t) &\mbox{$p$ even,}& \beta=1/2\\
    p!\tau_p(t)& \mbox{$p$ arbitrary,} & \beta=1\\
   p!\tau_{2p}(t/2) &\mbox{$p$ arbitrary,}& \beta=2,
   \end{array} \right.
 \eea
where in all three cases $\tau_p(t)$ is a
$\tau$-function satisfying the KP and Pfaff-KP
equations (8.0.4). Substituting (6.1.14)
 in the equation (8.0.4) and evaluating along the locus
$\LR$ leads to the following equations:
 {\footnotesize
\bean \left(\left(\frac{\pl}{\pl t_1}
\right)^4+3\left(\frac{\pl}{\pl
t_2}\right)^2-4\frac{\pl^2}{\pl t_1 \pl
t_3}\right)F_p+6\left(\frac{\pl^2}{\pl t^2_1}F_p
\right)^2
 &&  \\
 &&  \hspace{-1cm}
  =12\frac{p(p-1)}{(p+1)(p+2)}
  \frac{\tilde I_{p-2} \tilde I_{p+2}}{\tilde I_{p}^2}
 (1-\delta^{\beta}_1)\\
 &&   \hspace{1cm}~~\mbox{for}~~ \beta =1/2,1,\\
  \left(\left(\frac{\pl}{\pl t_1}
\right)^4+\frac{3}{4}\left(\frac{\pl}{\pl
t_2}\right)^2-\frac{\pl^2}{\pl t_1 \pl
t_3}\right)F_p+6\left(\frac{\pl^2}{\pl t^2_1}F_p
\right)^2
 &=&
  \frac{3}{4}\frac{p}{p+1}
   \frac{\tilde I_{p-1}\tilde I_{p+1}}{\tilde I_{p}^2}
 %(1-\delta^{\beta}_1)
 \\
 &&   \hspace{1cm}~~\mbox{for}~~ \beta =2.
\eean } \vspace{-1.8cm} \be\ee
%
%
%\newpage
Then, using $F(t) = \log \tilde I_p(t)$, with $\tilde
I_p$ as in (6.1.9) and (6.1.11),
  \bean
 H(y)&:=&
 %\frac{d}{d y}\log I_{p,q}
 % \\
 %
 %
 %&=&
 \frac{d}{d y}\log
 \int_{Gr(p,\BF^n)}e^{y\Tr(I+Z^{\dag}Z)^{-1}}\det(Z^{\dag
}Z)^{-\beta(q-p)} d\mu(Z)\\
&=& \frac{d}{d y}\log (e^{py/2}\\
 &&
   \int_{[-1,1]^p}
 |\Delta_p(z)|^{2\beta}\prod_1^p e^{
 yz_k/2}
  (1-z_k)^{\beta(n-q-p+1)-1}
 (1+z_k)^{\beta(q-p+1)-1}dz_k)
 \\
 &=&
  \frac{d}{dy}F_p(\frac{y}{2},0,0,\ldots)+\frac{p}{2}\\
 &=&
  \frac{1}{2}\left(\frac{\pl F_p}{\pl t_1}
   (\frac{y}{2},0,0,\ldots)+p\right)
 \eean \vspace{-1.5cm} \be \ee
 satisfies the differential equation, upon substituting
 the derivative (6.1.13) into (6.1.15) ($
 H':=\frac{dH}{dy}$):

\bean && \hspace{-1cm}
 4\left(y^3
H^{\prime\prime\prime}+6y^3{H^{\prime}
}^2+(1+\delta^{\beta}_{1})(2 y^2
 H^{\prime\prime}
 +4 y^2HH^{\prime}
  + yH^2)\right)-yP_0H^{\prime}+P_1H
  +P_2 \\
&&\\ &&=
 \left\{ \begin{array}{ll}
    0,~~& \mbox{for}~~ \beta =1,  \\ \\
   \displaystyle{ 3\frac{p(p-1)}{(p+1)(p+2)}y^3~
    \left.\frac{\tilde I_{p-2}\tilde I_{p+2}}{\tilde I_{p}^2}
\right|_{\LR}} ,~~& \mbox{for}~~ \beta =1/2, \\
        \\
  \displaystyle{ \left. \frac{3}{16}\frac{p}{p+1}y^3~
   \frac{\tilde I_{p-1}\tilde I_{p+1}}{\tilde I_{p}^2}
\right|_{\LR}},~~& \mbox{for}~~ \beta =2,
 \end{array}
  \right.
  \eean
 where $P_0,P_1,P_2$ are polynomials in $y$, with
 coefficients depending on
 $r=pq$ and $s=n-2p-2q$, given by table (6.1.3).

From
 $$ \tilde I_p(t)\bigr| _{\LR}
  = c_p^{(\beta)}e^{-px}I_p(x)=2^{p(\beta n -\beta p+\beta-1)}
  e^{-px}I_p(x),$$
it follows that for $\beta=1/2$ and $2$,
 $$
\left. \frac{\tilde I_{p-2} \tilde I_{p+2}}{\tilde
I_{p}^2}\right|_{\LR}=
  2^{-4}
   \frac{ I_{p-2}  I_{p+2}(x)}{I_{p}^2(x)}
  ~~ \mbox{and}~~
 \left. \frac{\tilde I_{p-1}\tilde I_{p+1}}{\tilde I_{p}^2}
   \right|_{\LR}=2^{-4}
   \frac{ I_{p-1}(x) I_{p+1}(x)}{ I_{p}^2(x)},
    $$
     thus establishing (6.1.2).
     \qed

{\medskip\noindent{\it Proof of Theorem 6.2:\/} }
 In particular, for $\beta =1$, from Theorem 4.2,
 Theorem 6.1 and (6.1.5),
   \bean
u(x)&:=&
 x \frac{\pl}{\pl x}\log E_{x,p}
   \left(I_{\{\lb\supseteq\mu\}}(\lb)\prod_{(i,j) \in \lb
    \atop n-q<j<n-p}
  h_{(i,j)}^{\lb}\right)\\
 &=& x \frac{\pl}{\pl x}\log \left(
 e^{-px}x^{p(n-p)}
%I_p(t)\Bigr|_{\LR}
\int_{Gr(p,\BC^n)}e^{x\Tr(I+Z^{\dag}Z)^{-1}}\det(Z^{\dag
}Z)^{-(q-p)} d\mu(Z)
  \right)\\
%  &=&
%  -px+p(n-p)+x\frac{d}{dx} \log I_{p}(x)\\
  &=&
  -px+p(n-p)+xH(x)
\eean
  satisfies the differential equation (6.1.5). From
  $I_{p}(x)=1+\frac{pq}{n}x+\ldots$, as in (4.2.6), it
  has the behavior near $x=0$, spelled out in (6.1.5),
  namely
  \be
 u(x)=p(n-p)-\frac{p(n-q)}{n} x+\ldots .
 \ee

Equation (6.1.7) follows from Cosgrove and Scoufis
(\cite{CosgroveScoufis}) (see \cite{AvM4}) , and the
constant is obtained by setting $x=0$ in the equation
and using the Taylor expansion (6.1.16) of $u(x)$
about $x=0$. \qed

% \newpage

\subsection{Hypergeometric functions, KP hierarchy and
          \\ Painlev\'e equations}

\begin{theorem}  The hypergeometric function for $\beta =1$ and $p,q,n\in\BC$,
   expressed in $t_i$-variables,
 \be
\tau(t):=\,_2F_1^{(1)}(p, q; n;u)\Biggl|_{\sum^m_{k=1}
u^i_k=it_i}, \ee satisfies the Hirota bilinear
relation as a function of $t_1,t_2,...,$ namely for
all $t,t'\in\BC^{\iy}$,
\be
\oint_{z=\iy}e^{\sum_1^{\iy}(t_i-t'_i)z^i}\tau(t-[z^{-1}])\tau(t'+[z^{-1}])dz=0.
\ee In particular, $\tau(t)$ satisfies the KP
hierarchy\footnote{Given a polynomial
$p(t_1,t_2,...)$, define the customary Hirota symbol
$p(\pl_t)f\circ g:= p(\frac{\pl}{\pl
y_1},\frac{\pl}{\pl y_2},...)f(t+y)g(t-y)
\Bigl|_{y=0}$. The ${ s}_{\ell}$'s are the elementary
Schur polynomials
$e^{\sum^{\iy}_{1}t_iz^i}:=\sum_{i\geq 0} {s}_i(t)z^i$
and for later use, set ${ s}_{\ell}(\tilde \pl):={
s}_{\ell}(\frac{\pl}{\pl
t_1},\frac{1}{2}\frac{\pl}{\pl t_2},\ldots).$}
($k=0,1,2,...$):
\be
\left(s_{k+4}\left(\frac{\pl}{\pl
t_1},\frac{1}{2}\frac{\pl}{\pl t_2},\frac{1}{3}
\frac{\pl}{\pl
t_3},\ldots\right)-\frac{1}{2}\frac{\pl^2}{\pl t_1\pl
t_{k+3}} \right)\tau\circ\tau =0. \ee
\end{theorem}

\proof Using $$ 1-a=e^{-\sum_1^{\iy}\frac{a^i}{i}}, $$
 one has the
following formula
 \bea  \prod_{{1\leq j\leq
p}\atop{1\leq k\leq m}}(1-z_ju_k)^{-\beta}
 \Bigl|_{t_i:=\frac{1}{i}\sum_{k=1}^mu_k^i}&=&
\prod_{{1\leq j\leq p}\atop{1\leq k\leq
m}}e^{\beta\sum_{i=1}^{\iy}\frac{u_k^i z_j^i}{i}}\nonumber\\
&=&\prod_{1\leq j\leq
p}e^{\beta\sum_{i=1}^{\iy}z_j^i\frac{1}{i}
 \sum^m_{k=1}u^i_k} \nonumber\\
&=&\prod_{1\leq j\leq
p}e^{\beta\sum_{i=1}^{\iy}t_iz_j^i}. \eea

On the one hand, setting $a=\beta(n-q)$ and $b=\beta
q$ in (2.3.14) using the new variables $t_i$,
using the constant (0.0.20) and (6.2.4), the hypergeometric function
 equals an integral for
integer $p\geq 1$,

\medbreak

$\displaystyle{c^{(\beta)}_{n,q,p}~\,_2F_1^{(1/\beta)}(\beta
p,\beta q, \beta n;u_1,...,u_m)}$
 \bea
&=&\int_{[0,1]^p}\prod_{{1\leq k\leq p}\atop{1\leq
j\leq m}}(1-z_ku_j)^{-\beta}
|\Delta_p(z)|^{2\beta}\prod^p_{\ell=1}
 z_{\ell}^{\beta(q-p+1)-1}
(1-z_{\ell})^{\beta(n-p-q+1)-1}dz_{\ell}\nonumber\\
&=&\int_{[0,1]^p}|\Delta_p(z)|^{2\beta}
 \prod^p_{\ell=1}e^{\beta\sum_{i=1}^{\iy}t_iz
_{\ell}^i}
z_{\ell}^{\beta(q-p+1)-1}
 (1-z_{\ell})^{\beta(n-p-q+1)-1}
 dz_{\ell}.
\eea According to (8.0.4), this integral is a solution
of the Hirota bilinear relation for $\beta =1$,
$q,n\in\BR$, integer $p>1$, such that $q-p+1$ and
$n-p-q+1>0$.

On the other hand, the hypergeometric function, is
also defined by (0.0.16), for $p,q,n\in\BR$ and
$\beta=1$:
 \be
\left.  \,_2F_1^{(1}(
p, q,  n;u_1,...,u_m)\right|_{t_i:=\frac{1}{i}\sum_{k=1}^mu_k^i}=\sum_{\lb\in\BY}\frac{(p)_{\lb}
(q)_{\lb}}{(n)_{\lb}h^{\lb}}{\bf s}_{\lb}(t)=:
 \sum_{\lb\in\BY}c_{\lb}{\bf s}_{\lb}(t).
\ee
 For integer $p\geq 1$, the integral (6.2.5) was shown to
be a solution of the Hirota bilinear relations and so
the coefficients $$ c_{\lb}\Bigl|_{\mbox{integer~}
p\geq 1} $$ satisfy Pl\"ucker relations.
They are homogeneous quadratic relations in a finite
number of the
$c_{\lb}$'s for $\lb\in\BY$ and they characterize the
KP $\tau$-functions; so, we have for every integer
$p\geq 1$, \bea
0&=&\sum_{\lb,\mu}\al_{\lb,\mu}c_{\lb}c_{\mu},\hspace*{2cm}\al_{\lb,\mu}=\pm
1\nonumber\\
&=&\sum_{\lb,\mu}\al_{\lb,\mu}\frac{(p)_{\lb}(q)_{\lb}
(p)_{\mu}(q)_{\mu}}{(n)_{\lb}(n)_{\mu}}\nonumber\\
&=:&\frac{1}{Y(n)}\sum_iX_i(n,q)p^i,
  \eea
   where $X_i(n,q)$ are
polynomials in $n,q$ and $\sum_i$ is a finite sum.
Setting $p=$ distinct integers $p_k$ in sufficient
number, one solves the homogeneous
 linear system (6.2.7) in
the $X_i(n,q)$ with coefficients $p^i_k$, whose
determinant is a Vandermonde, therefore non-zero. This
implies that $X_i(n,q)=0$ for all $n,q\in\BR$ and so,
we also have $$ \sum_iX_i(n,q)p^i=0,\qquad\mbox{for
all~}p\in\BR, $$ implying the Pl\"ucker relations (6.2.6)
for $p,q,n\in\BR$. One can extend the argument further
by analytic continuation to $p,q,n\in\BC$. \qed

% \newpage

\begin{theorem}
The hypergeometric function
\be
 H(x)=\frac{d}{dx} \log~{}_2F_1^{(1)}(p,q;n;y)
 \Bigr|_{\sum_{\ell}y_{\ell}^i= \delta_{1i}x}
 \ee
satisfies the Painlev\'e V equation:
\be
 4\left(x^3
H^{\prime\prime\prime}+6x^3{H^{\prime} }^2+2(2 x^2
 H^{\prime\prime}
 +4 x^2HH^{\prime}
  + xH^2)\right)-xP_0H^{\prime}+P_1H
  +P_2 =0\ee
  with $P_0,P_1,P_2$ as in (6.1.3).

  \end{theorem}

  \proof For integer $p,q$, with $0\leq p,q\leq n/2$,
  we have by (2.3.12),
\bean
_2F_1^{(1)}(p,q;n;u)\Biggl|_{\sum_{\ell}u_{\ell}^i=\delta_{1i}x}&=&1+\sum_{r\geq
1}b_rx^r\\ &=&\sum_{r\geq
0}x^r\sum_{\lb\in\BY_r}\frac{1}{(h^{\lb})^2}\frac{(p)_{\lb}(q)_{\lb}}{(n)_{\lb}}
,\quad p,q,n\in\BC \eean with $b_r$ rational in
$p,q,n$, and so $$
H(x):=\frac{d}{dx}\log\,_2F_1^{(1)}(p,q,n;y)\Biggl|_{\sum
y_{\ell}^i= \delta_{1i}x}=\sum_{r\geq 0}c_rx^r, $$
also with $c_r$ rational in $p,q,n$. Putting this
expression into the left hand side of the Painlev\'e
$V$ equation leads to $$ \sum^{\iy}_{r=0}d_rx^r,\quad
d_r=\mbox{~rational in $p,q,n$, with universal
coefficients}. $$ For integer $p,q,n$ with $n/2\geq
p,q\geq 0$, $H(x)$ satisfies the Painlev\'e equation
by the proof of Theorem 6.1
and so $d_r=0$. Using the same argument as in
 Theorem 6.3,
we have $d_r=0$ for all $p,q,n\in\BC$, whenever the
series $_2F_1^{(1)}$ makes sense.\qed

%\newpage

\section{Differential equations for
the Wishart and canonical correlations distributions
 %
% Differential equations for random ensembles
%in $T~Gr(m,\BF^n)$ and statistics
}

Remember from section 1, the tangent space to the
Grassmannians at $Id=\MAT{1}I_p\\O\mat$ is given by $$
T_{Id}Gr(p,\BF^n)=\{Z~\Big|~\mbox{arbitrary
$(n-p)\times p$ matrix with values in $\BF$}\}. $$ The
subgroup $$ K=\left\{B=\MAT{2} B_{11}&O\\O&B_{22}\mat,
B_{11}\in K_1, B_{22}\in K_2\right\} $$ acts on
$T_{Id}~Gr(m,\BF^n)$ as $$ Ad_B(Z)=B_{22}ZB^{-1}_{11},
$$ for which the spectrum $(\lb_1,\ldots,\lb_m)$ of
$Z^{\dag}Z$ is a invariant under the action of $K$ and
$\lb_i\geq 0$, since the matrix $Z^{\dag}Z$ is
positive definite.

Given $E\subseteq [0,\iy)$, consider the following
probability: \bea
P^{(\beta)}_{n,p}(E)&:=&P(\mbox{all~}\lb_i\in
E)\nonumber\\ &=&c_{p,n}\int_{{{Z\in
T_{Id}Gr(p,\BF^n)}\atop{\mbox{spectrum~}(Z^{\dag}Z)\subset
E}}}e^{-\Tr V(Z^{\dag}Z)}d\mu(Z)\nonumber\\
&=&\frac{\int_{E^p}|\Dt_p(z)|^{2\beta}\prod_{i=1}^n
e^{-V(z_i)}z_i^{\beta (n-2p+1)-1}
dz_i}{\int_{F^p}|\Dt_p(z)|^{2\beta}\prod_{i=1}^n
e^{-V(z_i)}z_i^{\beta (n-2p+1)-1}dz_i}, \eea where
$d\mu(Z)$ is Haar measure (1.0.11) on
$T_{Id}Gr(p,\BF^n)$ and where $\beta =1/2,1,2$
correspond to $\BR,\BC,\BH$ respectively. The function
$e^{-V(z)}$ will be either

%\newpage

\begin{description}

  \item[(i)] Laguerre: $e^{-V(z)}=e^{-bz}z^a$, $F=\BR^+$,
$\beta =1/2,1,2$
  \item[(ii)] Jacobi:  $e^{-V(z)}=(1-z)^az^b$, $F=[0,1]$,
$\beta =1/2,1,2$
 \item[(iii)] Gaussian: $e^{-V(z)}=e^{-z^{2}}$, $F=\BR^+$,
$\beta =1$, $n=2p$
\end{description}

We shall only consider (i) and (ii), the Gaussian
distribution being as in \cite{AvM3}.

%\newpage

\subsection{Laguerre probability on the tangent space to
 $Gr(p,\BF^n)$ and the Wishart distribution}

\begin{theorem}
%Here we consider the following probability, on the
%ensemble $T_{Id}Gr(p,\BF^n)$ with the Laguerre
%distribution $e^{-by}y^a$:

For $e^{-V(y)}:= e^{-by}y^a$, the probability (7.0.1)
defined on $Z \in$ the tangent space (at the identity)
to the symmetric space $ Gr(p,\BF^n)$, leads to the
following probability on the (positive) spectrum
$(\lb_1,\ldots, \lb_p)$ of $Z^{\dag}Z$,

\bea P_{n,p}^{(\beta)}(\max_{i} ~
 \lb_i \leq x%,~ \lb_i >0
 )
 &:=&c_{p,n}\int_{{{Z\in
T_{Id}Gr(p,\BF^n)}\atop{\mbox{spectrum~}
 (Z^{\dag}Z)\subset [0,x]
 }}}e^{-b\Tr Z^{\dag}Z }(\det Z^{\dag}Z)^a d\mu(Z)\nonumber\\
 %
% &=&P(\mbox{spectrum~}(Z^{\dag}Z)\in
%E, Z\in T_{Id}Gr(p,\BF^n))\nonumber\\
&=&\frac{\int_{[0,x]^p}|\Dt_p(y)|^{2\beta}\prod_{1}^p
e^{-by_i}y_i^{a+\beta(n-2p)+\beta
-1}dy_i}{\int_{[0,\iy)^p}|\Dt_p(y)|^{2\beta}\prod_{1}^p
e^{-by_i}y_i^{a+\beta(n-2p)+\beta-1}dy_i} \eea
 Then
 $$f(x)=x\frac{d}{dx}\log
   P_{p,n}^{(\beta)}(\max_{i} ~\lb_i \leq x )$$
satisfies %($f:=f_n(x)$)
\begin{itemize}
  \item
   for \underline{$\beta=1$}:
$$
x^2f^{\prime\prime\prime}+xf^{\prime\prime}
 +6xf^{\prime
2}-4 f~f^{\prime}-(2hf^{\prime}-h^{\prime}f)=0\quad
(\mbox{{\bf Painlev\'e V}}) $$ with $$
2h=(a+n-2p-bx)^2-4nbx. $$

  \item
   for \underline{$\beta=\left\{{1/2\atop 2}\right.$}:

   \bea
 \lefteqn{
 Q   \left(\frac{P^{(\beta)}_{n-{4\atop 2},p-{2\atop 1}}
  P^{(\beta)}_{n+{4\atop 2},p+{2\atop 1}}}
  {(P^{(\beta)}_{n,p})^2}-1 \right)
 -\left(3  f
+\frac{ b^2x^2 }{2\beta} -bQ_0x-3
  Q_1\right)f }\nonumber\\
 &=&
 x^3f^{\prime\prime\prime}-
  x^2f^{\prime\prime}+6x^2f^{\prime 2}
%  \nonumber\\
%&&~~
 -x\left(8 f+
\frac{b^2x^2}{2\beta}
-2bQ_0x-Q_2 \right) f^{\prime},
\nonumber\\ \eea
\end{itemize}
where
$$ \al :=a+\beta(n-2p)+\beta -1,
$$ and
 \bean
  Q&=&\left\{\begin{array}{ll}
\displaystyle{\frac{3}{4} p(p-1)( p+2\al)( p+2\al
+1)},~~~\mbox{for}~~\beta=1/2\\  \\ \displaystyle{
\frac{3}{2}p(2p+1)(2p+\al)(2p+\al-1)}
,~~~\mbox{for}~~\beta=2
\end{array}\right.\\
Q_2&=&
  6\beta p^{2}-\frac{\al^2}{2\beta}+6\al p+
   4(1- \beta)
  \al
  \\
  Q_1&=&
  2\beta p^{2}+2\al p+(1-\beta) \al,~~~~
  Q_0=   p+
 \frac{\al}{2\beta}.
  \eean  \vspace{-1cm} \be \ee

\end{theorem}

\begin{corollary} Consider a matrix $A$, which is Wishart
$W_p(n,\frac{1}{2b}I_p)$-distributed with eigenvalues
$\lb_1,\ldots,\lb_p$  (see Muirhead, p. 107). Then $$
P^W_{n,p}(x):=P^W_{n,p}(\max_i~\lb_i\leq x),\qquad
f:=x\frac{d}{dx}\log P^W_{n,p}(x) $$ satisfy the equation

\medbreak

$\displaystyle{\frac{3}{4}p(p-1)n(n-1)
 \left(\frac{P^W_{n+2,p+2}P^W_{n-2,p-2}}{(P^
W_{n,p})^2} -1\right)-(3f+b^2x^2-bxQ_0-3Q_1)f}$ $$
=x^3f^{\prime\prime\prime}-x^2f^{\prime\prime}+6x^2f^{\prime
2}-x(8f+(bx-Q_0)^2-4Q_1)f^{\prime}, $$ with
 $$
Q_0=\frac{1}{2}(n+p-1) ~~\mbox{and}~~4Q_1=
(n-1)(4p+1)-p,\qquad
%Q_2=-\frac{1}{4}(n-p-1)^2+(n-p-1) +  3p(n-1)
.$$

\end{corollary}

\proof According to (5.1.2), the Wishart distribution
$W_p(n,\frac{1}{2b}I_p)$ of the $p\times p$
matrix $A$ is given by
 \be
P^W_{n,p}(dA)= \Gamma_p(n/2)^{-1}b^{np/2}
 e^{-b\Tr A}
(\det A)^{\frac{1}{2}(n-p-1)}
 \prod_{{1\leq i\leq j \leq p}} dA_{ij}.
\ee
and so the joint probability $P^W_{n,p}(\max_i~\lb_i\leq x)$
 is precisely formula (7.1.1)
, with
 $a=p/2$ and $\beta =1/2$,
 \bean
P^W_{n,p}(\max_i~\lb_i\leq x)&=&
 c_{n,p,b}\int_{[0,x]^p}|\Delta_p(y)|\prod_1^p
e^{-by_i}y_i^{\frac{1}{2}(n-p-1)}dy_i  \\
&=&P_{n,p}^{(\beta)}(\max_i~\lb_i\leq x)\Bigr|_
 {{a=p/2}\atop
{\beta=1/2} }. \eean
Therefore $P^W_{n,p}(x)$ also satisfies the inductive
differential equation (7.1.2); we only need to check
that
 \bean P^W_{n\pm
2,p\pm 2}(x)&=&c\int_{[0,x]^{p\pm 2}}|\Delta_{p\pm
2}(y)|\prod_1^{p\pm
2}e^{-by_i}y_i^{\frac{1}{2}(n-p-1)}dy_i\\
&=&c\int_{[0,x]^{p\pm 2}}|\Delta_{p\pm
2}(y)|\prod_1^{p \pm 2}e^{-by_i}y_i^{p/2}
 y_i^{\frac{1}{2}((n\pm
4)-2(p\pm 2))- \frac{1}{2}}\\ &=& P^{(1/2)}_{n\pm
4,p\pm 2}(\max_i~\lb_i\leq x)\mbox{~for~}a=p/2. \eean
$Q,Q_0,Q_1,Q_2$ can immediately be computed by setting
$\beta =1/2$, $\delta^{\beta}=1$ and $2\al =n-p-1$ in
(7.1.3).\qed

%\newpage

\subsection{Jacobi probability on the tangent space to
 $Gr(p,\BF^n)$ and the sample canonical correlation
 distribution }

\begin{theorem}

For $e^{-V(y)}:= (1-y)^ay^b$, the probability (7.0.1)
defined on $Z \in$ the tangent space (at the identity)
to the symmetric space $ Gr(p,\BF^n)$, leads to the
following probability on the (positive) spectrum of
$Z^{\dag}Z$,
\bea
P^{(\beta)}_{n,p}(\max_i \lb_i\leq \frac{x+1}{2})
&:=&c_{p,n}\int_{{{Z\in
T_{Id}Gr(p,\BF^n)}\atop{\mbox{spectrum~}
 (Z^{\dag}Z)\subset [0,\frac{x+1}{2}]
 }}}e^{-\Tr V(Z^{\dag}Z)}d\mu(Z)\nonumber\\
&=&\frac{\int_{[0,\frac{x+1}{2}]^p}
 |\Dt_p(y)|^{2\beta}\prod_{i=1}^p
(1-y_i)^ay_i^{b+\beta (n-2p+1)-1}
dy_i}{\int_{[0,1]^p}|\Dt_p(y)|^{2\beta}\prod_{i=1}^p
(1-y_i)^ay_i^{b+\beta (n-2p+1)-1}dy_i};\nonumber\\
\eea
it satisfies the following differential equations,
upon setting $$ f(x):=(1-x^2)\frac{d}{dx}  \log
P^{(\beta)}_{n,p}(\max_{i} \lb_i \leq
\frac{x+1}{2}),$$

\begin{itemize}
  \item

 for \underline{$\beta=1$}: \bea
 &&\hspace{-1.3cm} 2(x^2-1)^2f^{\prime\prime\prime}
 +4(x^2-1)\left(xf^{\prime\prime}
 -3f^{\prime 2}\right)
 +\left(16
xf-u(x^{2}-1)-2sx-r \right)f^{\prime}
 \nonumber\\&&~~~~~-f\left(4f-ux-s \right)=0
 \eea

\item
for \underline{$\beta=\left\{{1/2\atop 2}\right.$}:

$\displaystyle{Q\left(\frac{P^{(\beta)}_{n+{4\atop
2},p+{2 \atop 1}}
 P^{(\beta)}_{n-{4\atop
2},p-{2 \atop 1}}} {(P^{(\beta)}_{n,p})^2}-1\right)}$
 {\footnotesize
\bea &=& 4(u+1)(x^2-1)^{2}\Bigl(-u(x^2-1)
f^{\prime\prime\prime}+
 (12f -ux-3s)
  f^{\prime\prime}+6 u(u-1)
 f^{\prime 2}\Bigr)\nonumber\\
  &&-(x^2-1)f^{\prime}\Bigl(24f(u+3)(2f-s)+8fu
  (5u-1)x-
   u(u+1)(ux^{2}+2sx+8)+
 Q_2
 \Bigr)
 \nonumber\\
 && +f\Bigl(48f^{3}+48f^{2}(ux+2x-s)+2f\left(8u^{2}x^{2}
 +2ux^{2}-12usx-24s
 x+Q_4\right)\nonumber\\
  &&~~~~~-u(u+1)x(3ux^{2}+sx-2
 ux-3u)+Q_3x- Q_1s\Bigr), \eea
 }

\end{itemize}

\vspace{-1cm}
where

$$b_0=a-b-\beta(n-2p+1)+1
,\quad b_1=a +b+\beta(n-2p+1)-1 $$
$$
  r=\frac{2}{\beta}(b_0^2+(b_1+2-2\beta)^2)
%  = \frac{4}{\beta}( (a-b)^2+(a+b +2-\beta)^2)
%
~~~~~~s=\frac{2}{\beta}b_0(b_1+2-2\beta)
%=\frac{4}{\beta}(a-b)(a+b +2-\beta)
% \\
 % m&=& \beta n+\gamma+\delta
%\\
$$ $$ u= \frac{2}{\beta}(2\beta p+b_1
+2-2\beta)(2\beta p
 +b_1),
%=\frac{4}{\beta}m(m+2-\beta)
$$
 and %the following {\em invariant} polynomials in $u,r,s$:
 \bean
 Q&=&\frac{3}{16}\left(
(s^{2}-ur+u^{2})^{2}-4(rs^{2}-4u
 s^{2}-4s^{2}+u^{2}r)\right)
\\
  Q_1&=& 3s
 ^{2}-3ur-6r+2u^{2}+23u+24\\
 Q_2&=&3us^{2}
 +9s^{2}-4u^{2}\,r+2ur+4u^{3}+10u^{2}
 \\
 Q_3&=&3us^{2}+6s^{2}-3u^{2}
 r+u^{3}+4u^{2}
 \\
 Q_4&=& 9s
 ^{2}-3ur-6\,r+u^{2}+22u+24=Q_1+(6s^2-u^2-u).\\
 \eean
\vspace{-2cm} \be \ee

\end{theorem}

\begin{corollary} Let $A=ZZ^{\top}$ have the Wishart
distribution
$W_{p+q}(n,\Sigma)$-distribution.
 Break up the
matrices $\Sigma$ and $A$, as follows:

  $\hspace*{23mm}\stackrel{p}{\longleftrightarrow}
\quad\stackrel{q}{\longleftrightarrow} $
\vspace{-.3cm}
 $$
  \Sigma=
 \left(\begin{array}{cc}\Sigma_{11}&\Sigma_{12}\\
\Sigma_{12}^{\top}&\Sigma_{22}\end{array}\right)
 \begin{array}{l}\updownarrow p\\ \updownarrow q\end{array}
~~~~~~~~A=
 \left(\begin{array}{cc}A_{11}&A_{12}\\
A_{12}^{\top}&A_{22}\end{array}\right)
 \begin{array}{l}\updownarrow p\\ \updownarrow q\end{array}
$$
 Assume the eigenvalues $\rho^2_1,\ldots,\rho^2_p$
of $\Sigma_{11}^{-1}\Sigma_{12}\Sigma_{22}^{-1}
\Sigma_{12}^{\top}$ all zero. Then the probability
distribution of the eigenvalues
 $r_1^2,\ldots,r_p^2$ of $A_{11}^{-1}A_{12}A_{22}^{-1}
   A_{12}^{\top}$ (sample canonical correlation
    coefficients) is given by
  \bean
 P^C_{n,p,q}(x)&=&
 P(0\leq  r_i^2\leq \frac{1+x}{2}~~\mbox{for } 1\leq i\leq p)
\\ &=&\frac
 {\int_{[0,\frac{1+x}{2}]^p}
|\Dt_p(z)|\prod_{i=1}^{p}z_i^{\frac{1}{2}(q-p-1)}
(1-z_i)^{\frac{1}{2}(n-p-q-1)}dz_i}
 {
 \int_{[0,1]^p}
|\Dt_p(z)|\prod_{i=1}^{p}z_i^{\frac{1}{2}(q-p-1)}
(1-z_i)^{\frac{1}{2}(n-p-q-1)}dz_i }\eean
 and
  $$
  f(x):=(1-x^2)\frac{d}{dx}  \log P^C_{n,p,q}(x)$$
  satisfy the inductive PDE:
%
%
%\bigbreak
 $$
 Q\left(\frac{P^C_{n+4,p+2,q+2}
P^C_{n-4,p-2,q-2}} {{P^C_{n,p,q}}^2}-1\right)
 =\left\{\begin{array}{l}\mbox{same expression as the}\\
 \mbox{right hand side of (7.2.3)}
 \end{array}\right\},$$
with $Q_1,Q_2,Q_3$ being symmetric polynomials of
$p,q$, given by (7.2.4), where
 \bean u&=&n(n-2)\\
   r/4&=& n^{2}/2-n(p+q)+p^{2}+q^{2}
   \\
  s&=&  (n-2 p) (n-2 q )\\
%&&\\
  Q&=&48pq\left(p-1\right)\left(q-1\right)
  \left(n-p\right)
 \left(n-q\right)\left(n-p-1\right)\left(n-q-1\right)
 .\eean

 \end{corollary}

 \remark
 For instance
 \bean
   Q_1&=&24\left(n-1\right)\left(p^2+q^{2}-n\left(p+q\right)
\right)
 +48pq\left(n^{2}-n\left(p+q\right)+pq\right)\\
  &&-
  \left(n-6
 \right)\left(n-1\right)^{2}\left(n+4\right).\eean

\proof
% Remembering the
% $$
%2^{-n(p+q)/2}\Gamma_{p+q}(n/2)^{-1}(\det\Sigma)^{-n/2}
% e^{-\frac{1}{2}\Tr\Sigma^{-1}A}
%(\det A)^{\frac{1}{2}(n-p-q-1)}, $$
%
The proof follows immediately from Corollary 5.2 and Theorem 7.3.\qed

\section{Appendix: The Pfaff-KP hierarchy and Virasoro constraints}

Consider weights of the form $\rho(z)dz:=e^{-V(z)}dz$
on an interval $F=[A,B]\subseteq\BR$, with rational
logarithmic derivative and subjected to the following
boundary conditions:
\be
-\frac{\rho'}{\rho}=V^{\prime}=\frac{g}{f}=\frac{\sum_0^{\iy}b_iz^i}{\sum_0^{\iy}a_iz^i},
\quad \lim_{z\rightarrow
A,B}f(z)\rho(z)z^k=0\mbox{\,\,for all\,\,}k\geq 0,
  \ee

\begin{theorem}

 The multiple integrals
 %$I(t,c;\beta)=(I_0 =1,I_1(t,c;\beta),...)$, with
 \bea
I_p(t;\beta)&=&\int_{F^p}|\Dt_p(z)|^{2\beta}\prod_{k=1}^p
\left(e^{\sum_1^{\iy}t_i z_k^i}\rho(z_k)dz_k\right)
,~\mbox{for} ~~p>0
\nonumber\\
 &=&
   \left\{
   \begin{array}{lll}
   p!\tau_p(t) &\mbox{$p$ even,}& \beta=1/2\\
    p!\tau_p(t)& \mbox{$p$ arbitrary,} & \beta=1\\
   p!\tau_{2p}(t/2) &\mbox{$p$ arbitrary,}& \beta=2
   \end{array} \right.
 \eea
 with $I_0=1$, satisfy

 \item[(i)] the following
 Virasoro constraints
 %\footnote{
% When $E$ equals the whole range $F$, then the
%  ${\cal B}_k$'s are absent in the formulae (2.1.7).}
  for all $k\geq -1$:
\bea && \hspace{-1cm}\sum_{i\geq 0}\left( a_i~
{}^{\beta}\BJ_{k+i,p}^{(2)}(t,p)-b_i ~ {}^{\beta}
\BJ_{k+i+1,p}^{(1)}(t,p)\right) I_p(t;\beta) = 0,
%\nonumber
 \eea
  in terms of the coefficients
$a_i,~b_i$ of the rational function $(-\log \rho)'$

\item[(ii)] The \underline{Pfaff-KP hierarchy}: (see
 footnote 14 for notation)
\be
\left({\bf
s}_{k+4}(\tilde\pl)-\frac{1}{2}\frac{\pl^2}{\pl t_1\pl
t_{k+3}}\right)\tau_{p}\circ\tau_{p}=
 (1-\delta^{\beta}_{1}){\bf
s}_k(\tilde \pl)~\tau_{p+2}\circ\tau_{p-2} \ee \hfill
$p ~\mbox{even},~ k=0,1,2,...~.$

 \noindent of which the first equation reads ($p$
 even)
 \bea
  \lefteqn{\left(\left(\frac{\pl}{\pl t_1}
\right)^4+3\left(\frac{\pl}{\pl
t_2}\right)^2-4\frac{\pl^2}{\pl t_1 \pl
t_3}\right)\log\tau_p+6\left(\frac{\pl^2}{\pl
t^2_1}\log\tau_p
\right)^2}\nonumber\\
&&\hspace{6cm}=12\frac{\tau_{p-2}\tau_{p+2}}{\tau_{p}^2}
(1- \delta_{\beta,1}).
\eea

%\newpage

\item[(iii)] More generally, the functions
$\tau(t)$ satisfy the Hirota bilinear
relations for all
$t,t'\in\BC^{\iy}$ and $m,p$ positive integers (see
footnote 10 for notation)
 \begin{itemize}
  \item $\beta=1$
\be
\oint_{z=\iy}
\tau_p(t-[z^{-1}])\tau_p(t'+[z^{-1}])
 e^{\sum_1^{\iy}(t_i-t'_i)z^i}dz=0,
% \mbox{for $\beta=1$}
\ee
 \item $\beta=1/2$ and $2$
\begin{multline}
\oint_{z=\iy}\tau_{2p}(t-[z^{-1}])\tau_{2m+2}(t'+[z^{-1}])
e^{\sum_{0}^\iy(t_i-t'_i)z^i} z^{2p-2m-2}dz
 %\mbox{for $\beta=1/2$ and $2$}
  \\
{}+\oint_{z=0}\tau_{2p+2}(t+[z])\tau_{2m}(t'-[z])
e^{\sum_{0}^\iy(t'_i-t_i)z^{-i}}z^{2p-2m}dz=0\,,
%\label{PF3.1}
\end{multline}

\end{itemize}
\end{theorem}

\subsubsection*{Example   (Jacobi $\beta$-integral)}
 This case is particularly important, because it covers
the integrals in Theorems 0.2 and 0.3. The weight and
the $a_i$ and $b_i$, as in (8.0.1), are given by
 $$
\rho (z):=e^{-V}=(1-z)^{a}(1+z)^{b} ,
V'=\frac{g}{f}=\frac{a-b+(a+b)z}{1-z^2}
 $$
 $$
 a_0=1,a_1=0,a_2=-1,b_0=a-b,b_1=a+b
 ,~\mbox{and all
 other}~a_i,b_i=0 .$$
 The integrals
 \be
I_p= \int_{E^p}|\Dt_p(z)|^{2\beta}\prod_{k=1}^p
(1-z_k)^{a}(1+z_k)^{b}e^{\sum_{i=1}^{\iy}t_iz_k^i}dz_k
 \ee
satisfy the Virasoro constraints $(k\geq-1)$:
 \bea
 {\cal J}^{(2)}_{k}I_p=
\left(~{}^{\beta}\BJ_{k+2,p}^{(2)}-~{}^{\beta}\BJ_{k,p}^{(2)}+
b_0~{}^{\beta}\BJ_{k+1,p}^{(1)}
+b_1~{}^{\beta}\BJ_{k+2,p}^{(1)}\right)I_p=0. \eea
 Introducing
  $ \sigma_i=(2p-i-1)\beta+i+1+b_1,
$

Then introducing the function $F_p:=\log \tau_p(t)$,
the two first Virasoro constraints for $m=1,2$ divided
by $\tau_p$ are given by \bea \frac{{\cal
J}^{(2)}_{-1}\tau_p}{\tau_p}&=&\left(\sum_{i\geq
1}it_i\frac{\pl}{\pl t_{i+1}}- \sum_{i\geq
2}it_i\frac{\pl}{\pl t_{i-1}} +\sigma_1\frac{\pl}{\pl
t_1}\right) F_p+p(b_0-t_1)=0 \nonumber\\
\frac{{\cal J}^{(2)}_{0}\tau_p}{\tau_p}
  &=& \left(\sum_{i\geq 1}it_i
  \Bigl(\frac{\pl}{\pl t_{i+2}}
  -\frac{\pl}{\pl t_{i}}\Bigr) +
 b_0\frac{\pl}{\pl t_1}
 +\beta\frac{\pl^2}{\pl t_1^2}
+\sigma_2\frac{\pl}{\pl t_2}\right) F_p \nonumber\\
 &&\hspace{5cm}+\beta\left(\frac{\pl F_p}{\pl t_1}
 \right)^2-\frac{p}{2}(\sigma_1-b_1)=0.\nonumber\\
\eea

%\newpage

\end{document}